\documentclass[leqno,11pt]{amsart}
\usepackage{amsmath,amscd,amsthm,amsxtra}
\usepackage{epsfig,graphics,color,colortbl}
\usepackage{amssymb,latexsym}
\usepackage{mathabx}
\usepackage{mathrsfs,eucal,upgreek}
\usepackage[poly,all]{xy}
\usepackage{hyperref}
\usepackage{tikz-cd}
\usepackage{arydshln}

\newtheorem{thm}{\bf Theorem}[section]
\newtheorem{df}[thm]{\bf Definition}
\newtheorem{prop}[thm]{\bf Proposition}
\newtheorem{cor}[thm]{\bf Corollary}
\newtheorem{lem}[thm]{\bf Lemma}
\newtheorem{rem}[thm]{\bf Remark}
\newtheorem{ex}[thm]{\bf Example}

\numberwithin{equation}{section}

\newcommand{\bs}{\boldsymbol}
\newcommand{\A}{\mathcal{A}}

\newcommand{\B}{\mathbf{B}}

\newcommand{\W}{\mathcal{W}}
\newcommand{\cP}{\mathscr{P}}

\newcommand{\pf}{\noindent{\bfseries Proof. }}
\newcommand{\ov}{\overline}

\newcommand{\bi}{\bs{\rm i}}

\newcommand{\N}{\mathbb{N}}
\newcommand{\gl}{\mathfrak{gl}}
\newcommand{\Z}{\mathbb{Z}}

\newcommand{\te}{\widetilde{e}}
\newcommand{\tf}{\widetilde{f}}
\newcommand{\g}{\mathfrak{g}}
\newcommand{\td}{\widetilde}

\newcommand{\mc}{\mathcal}
\newcommand{\mf}{\mathfrak}

\newcommand{\J}{\mathbb{J}}

\newcommand{\la}{\lambda}

\newcommand{\red}[1]{{\color{red}#1}}

\begin{document}
\title[Lusztig data of KN tableaux]
{Lusztig data of Kashiwara-Nakashima tableaux in types $B$ and $C$}
\author{JAE-HOON KWON}

\address{Department of Mathematical Sciences, Seoul National University, Seoul 08826, Korea}
\email{jaehoonkw@snu.ac.kr}


\thanks{This work was supported by Samsung Science and Technology Foundation under Project Number SSTF-BA1501-01.}

\begin{abstract}
We provide an explicit combinatorial description of the embedding of the crystal of Kashiwara-Nakashima tableaux in types $B$ and $C$ into that of $\bi$-Lusztig data associated to a family of reduced expressions $\bi$ of the longest element $w_0$. The description of the embedding is simple and elementary using only the Sch\"{u}tzenberger's jeu de taquin and RSK algorithm.
A spinor model for classical crystals plays an important role as an intermediate object connecting Kashiwara-Nakashima tableaux and Lusztig data. 
\end{abstract}

\maketitle
\setcounter{tocdepth}{2}


\section{Introduction}
Suppose that $\g$ is a complex finite-dimensional simple Lie algebra. Let $B(\infty)$ be the crystal of the negative part of the quantum group $U_q(\g)$ and let $B(\la)$ be the crystal of an integrable highest weight $U_q(\g)$-module with highest weight $\la$ \cite{Kas91}.

Let $\bi$ be a reduced expression of the longest element $w_0$ of length $N$ in the Weyl group of $\g$. 
There is a realization of $B(\infty)$ given by the set of $N$-tuples of non-negative integers called the {\em $\bi$-Lusztig data}, which parametrizes a PBW basis associated to $\bi$ and hence the canonical basis of the negative part of $U_q(\g)$ \cite{Lu90, Lu93}. 
The crystal of this Lusztig parametrization has a rich combinatorial structure including a connection with the geometry of the totally positive varieties \cite{BZ-1}, and is also closely related with other important realizations of $B(\infty)$, for example the Kashiwara's string parametrizaion and MV polytopes \cite{BZ-1,Kam}.

On the other hand, when $\g$ is of classical type, there is a well-known combinatorial model for $B(\la)$ called {\em Kashiwara-Nakashima tableaux} ({\em KN tableaux} for short) \cite{KashNaka}, which generalizes the notion of Young tableaux in type $A$, and has been an important tool for studying crystals of classical type \cite{Na} and the Kirillov-Reshetikhin crystals of non-exceptional affine type \cite{FOS,LOS}.

The purpose of this paper is to give an explicit description of the relationship between these two primarily important combinatorial models for $B(\infty)$ and $B(\la)$ in types $B$ and $C$. 

Let $\B_{\bi}$ be the set of $\bi$-Lusztig data and ${\bf KN}_\la$ the set of KN tableaux of shape $\la$. We assume that the Weyl group of type $B_n$ or $C_n$ is embedded into that of type $A_{2n-1}$ as the set of fixed elements under the Dynkin symmetry. 
Let $\bi$ be a reduced expression of $w_0$ such that the image of $\bi$ in the Weyl group of type $A$ is adapted to a Dynkin quiver of type $A_{2n-1}$, which has a single sink in the middle vertex.

The main result in this paper is to give a combinatorial description of the crystal embedding
\begin{equation}\label{eq:embedding KN to Bi}
\xymatrixcolsep{4pc}\xymatrixrowsep{3pc}\xymatrix{
{\bf KN}_\la \ \ar@{^{(}->}[r] &\ \B_{\bi}\otimes T_{\omega_\la}},
\end{equation}
where $\omega_\la$ is the highest weight corresponding to $\la$, and $T_{\omega_\la}=\{t_{\omega_{\la}}\}$ is an abstract crystal with ${\rm wt}(t_{\omega_\la})=\omega_\la$ and
$\varphi_i(t_{\omega_{\la}})=-\infty$ for each simple root index $i$ (Theorem \ref{thm:main}). 

Let us give some remarks on previous works related to this problem.
It is well-known that there is a concrete description of the embedding in case of type $A$ with $\bi=(1,2,1,3,2,1,\ldots)$, where the $\bi$-Lusztig datum of a Young tableau is determined by counting the number of its entries in each row (see \cite{K16-2} for a generalization to the case when ${\bf i}$ is adapted to a quiver of type $A$ with a single sink). 
For other types, recall that there is another realization of $B(\infty)$ of type $ABCD$ and type $G_2$ in terms of marginally large tableaux \cite{HL}.
Recently, 
 a simple combinatorial description of the isomorphism from the crystal of marginally large tableaux of type $BCD$ to that of $\bi$-Lusztig data is obtained in \cite{SST16,SST17} for some $\bi$ which is different from our choice. Using the insertion scheme of $BCD$ types, one may describe an embedding of ${\bf KN}_\la$ into the crystal of marginally large tableaux to have an embedding from ${\bf KN}_\la$ into $\B_{\bi}$, but it seems  more involved than in type $A$ (see Remark \ref{rem:related works}).

The description of \eqref{eq:embedding KN to Bi} presented here is simple and elementary using only the Sch\"{u}tzenberger's jeu de taquin sliding and RSK algorithm (see Examples in Section \ref{sec:examples}).
Let us explain our result in more details.
A key ingredient of our approach is to 
consider intermediate combinatorial models of crystals between ${\bf KN}_\la$ and $\B_{\bi}$, and factorize the embedding into a composite of the following three morphisms; 
\begin{equation}\label{eq:factorization of embedding}
\xymatrixcolsep{4pc}\xymatrixrowsep{3pc}\xymatrix{
{\bf KN}_\la \ar@{->}^{\Psi_\la }[r] &
{\bf T}_\la  \ar@{->}^{\!\!\!\!\!\!\!\!\Theta_\la }[r] &
{\bf V}_\la\otimes T_{r\omega_n} \ar@{->}^{\Phi_\la }[r] &
\ \B_{\bi}\otimes T_{\omega_\la}},
\end{equation}
where $r=\langle \omega_\la,h_n\rangle$, ${\bf T}_\la$ is another tableau model of $B(\omega_\la)$ introduced in \cite{K15,K16}, and ${\bf V}_\la$ is the crystal of a parabolic Verma module with respect to a maximal Levi subalgebra $\mf l$ of type $A_{n-1}$ in $\g$ \cite{K09,K12}.

The crystal ${\bf T}_\la$ is given by a set of sequences of Young tableaux of two-columned skew shapes satisfying certain conditions of admissibility. It is based on the fact that a non-spinor fundamental weight module can be embedded in a two-fold tensor power of the spin representation of a $q$-deformed Clifford algebra \cite{Ha}. Hence the isomorphism $\Psi_\la$ maps each column of a non-spinor weight in KN tableaux to a two-columned skew Young tableaux (Theorems \ref{thm:KN_C_to_psst} and \ref{thm:KN_B_to_psst}), where each single column on the side of ${\bf T}_\la$ corresponds to a monomial in a spin representation. We may call ${\bf T}_\la$ a {\em spinor model} of $B(\la)$ for this reason (see \cite{K15-1,K15,K16} for other applications). 

The crystal ${\bf V}_\la$ is given by a set of pairs of Young tableaux, which was introduced in the study of irreducible characters of classical type as a truncation of Cauchy or Littlewood identity (multiplied by a Schur polynomial) and their stable branching formulas. 
It is the crystal of a parabolic Verma module in the sense that ${\bf V}_\la$ is a limit of the crystals of finite-dimensional $U_q(\g)$-modules, whose highest weight with respect to the Levi subalgebra $\mf l$ is $\la$.
Using a combinatorial description of $B(\infty)$ in terms of quiver representations \cite{Re} and a crystal theoretic interpretation of RSK algorithm \cite{K09}, the embedding $\Phi_\la$ maps a pair of Young tableaux to a pair of upper triangular matrix which are the $\bi$-Lusztig data corresponding to the Levi subalgebra $\mf l$ and the associated nilradical of $\g$, respectively (Theorem \ref{eq:Phi isomorphism}).

Now the classical methods in type $A$ combinatorics naturally enters into this picture of types $B$ and $C$ since both of the crystals ${\bf T}_\la$ and ${\bf V}_\la$ are given in terms of Young tableaux. 
Indeed, we show that the embedding $\Theta_\la$ is obtained by successive application of the Sch\"{u}tzenberger's jeu de taquin sliding to a tableau in ${\bf T}_\la$ (Theorem \ref{thm:Theta isomorphism}). The description of $\Theta_\la$ is a crucial part of \eqref{eq:factorization of embedding}. 
We conjecture a similar result for type $D$ using a spinor model of type $D$ \cite{K16}, while we need an approach different from types $B$ and $C$ since we do not know yet a direct analogue of $\Phi_\la$ in \eqref{eq:factorization of embedding} for type $D$. 


The paper is organized as follows. In Section \ref{sec:prel}, we briefly recall some background on necessary materials. In Section \ref{sec:spinor}, we review the notion of spinor model ${\bf T}_\la$ and describe the isomorphism $\Psi_\la$ from ${\bf KN}_\la$ to ${\bf T}_\la$. In Section \ref{sec:parabolic}, we discuss the crystal structure on $\B_{\bi}$ and
describe the embedding $\Phi_\la$ of ${\bf V}_\la\otimes T_{r\omega_n}$ into $\B_{\bi}\otimes T_{\omega_\la}$. Finally in Section \ref{sec:main}, we give a description of the embedding $\Theta_\la$ of ${\bf T}_\la$ into ${\bf V}_\la\otimes T_{r\omega_n}$.\vskip 2mm

{\bf Acknowledgement} The author would like to thank the referee for helpful comments and Il-Seung Jang for very careful reading of the manuscript.

\section{Preliminaries}\label{sec:prel}
\subsection{Notations}
Let $\Z_+$ be the set of non-negative integers.
Let $\mathscr{P}$ be the set of partitions. 
We identify $\lambda=(\la_1,\la_2,\ldots)\in \cP$ with a {Young diagram}. We denote by $\la'=(\la'_1,\la'_2,\ldots)$ the conjugate of $\la$, and by $\lambda^\pi$ the skew Young diagram obtained by $180^{\circ}$-rotation of $\lambda$. For $n\geq 1$ let $\cP_n=\{\,\la=(\la_1,\ldots,\la_n)\,|\,\la_i\in \Z_+, \la_i\geq \la_{i+1}\ \text{for all $i$}\,\}\subset \cP$.

Let $\mc{A}$ be a linearly ordered  set with a $\mathbb{Z}_2$-grading $\mc{A}=\mc{A}_0\sqcup\mc{A}_1$. For a skew Young diagram
$\lambda/\mu$, let $SST_\A(\lambda/\mu)$ be the set of all
semistandard tableaux of shape $\lambda/\mu$ with entries in $\A$, that is, (1) the entries in each row (resp. column) are
weakly increasing from left to right (resp. from top to bottom), (2)
the entries in $\mc{A}_0$ (resp. $\mc{A}_1$) are strictly increasing in each
column (resp. row).
For $T\in SST_\A(\lambda/\mu)$, let $w(T)$ be a word obtained by reading the entries of $T$ row
by row from top to bottom, and from right to left in each row (cf.~\cite{Ful}).
We denote by $\W_\A$ the set of words of finite length with letters in $\A$.

\subsection{Crystals}
Let us give a brief review on crystals (see \cite{HK,Kas95} for more details).
Let $\g$ be the Kac-Moody algebra associated to a symmetrizable generalized Cartan matrix
$A =(a_{ij})_{i,j\in I}$ indexed by a set $I$. Let $P^\vee$ be the dual weight lattice,
$P = {\rm Hom}_\Z( P^\vee,\Z)$ the weight lattice,
$\Pi^\vee=\{\,h_i\,|\,i\in I\,\}\subset P^\vee$ the set of simple coroots, and
$\Pi=\{\,\alpha_i\,|\,i\in I\,\}\subset P$ the set of simple roots of $\g$ such that $\langle \alpha_j,h_i\rangle=a_{ij}$ for $i,j\in I$. Let $(d_i)_{i\in I}$ be a sequence of positive integers such that $(d_ia_{ij})$ is symmetric.

Let $P^+$ be the set of integral dominant weights.
Let $W$ be the Weyl group of $\g$ generated by the simple reflection $s_i$ for $i\in I$.

A {\it $\g$-crystal} (or simply a {\it crystal} if there is no confusion on $\g$) is a set
$B$ together with the maps ${\rm wt} : B \rightarrow P$,
$\varepsilon_i, \varphi_i: B \rightarrow \mathbb{Z}\cup\{-\infty\}$ and
$\te_i, \tf_i: B \rightarrow B\cup\{{\bf 0}\}$ for $i\in I$ satisfying the following:
\begin{itemize}
\item[(1)]  
$\varphi_i(b) =\langle {\rm wt}(b),h_i \rangle +
\varepsilon_i(b)$,

\item[(2)] $\varepsilon_i(\te_i b) = \varepsilon_i(b) - 1,\ \varphi_i(\te_i b) =
\varphi_i(b) + 1,\ {\rm wt}(\te_ib)={\rm wt}(b)+\alpha_i$ if $\te_i b \in B$,

\item[(3)] $\varepsilon_i(\tf_i b) = \varepsilon_i(b) + 1,\ \varphi_i(\tf_i b) =
\varphi_i(b) - 1,\ {\rm wt}({\tf_i}b)={\rm wt}(b)-\alpha_i$ if $\tf_i b \in B$,

\item[(4)] $\tf_i b = b'$ if and only if $b = \te_i b'$ for $b' \in B$,

\item[(5)] $\te_ib=\tf_ib={\bf 0}$ when $\varphi_i(b)=-\infty$,
\end{itemize}
for $b\in B$ and $i\in I$.
Here ${\bf 0}$ is a formal symbol and $-\infty$ is the smallest
element in $\Z\cup\{-\infty\}$ such that $-\infty+n=-\infty$
for all $n\in\Z$. 
For $b, b'\in B$, we write $b\stackrel{i}{\rightarrow}b'$ if $b'=\tf_i b$ for $i\in I$.

A  tensor product $B_1\otimes B_2$  of crystals $B_1$ and $B_2$
is defined to be a crystal, which is $B_1\times B_2$  as a set with elements  denoted by
$b_1\otimes b_2$, satisfying  {\allowdisplaybreaks
\begin{align*}
{\rm wt}(b_1\otimes b_2)&={\rm wt}(b_1)+{\rm wt}(b_2), \\
\varepsilon_i(b_1\otimes b_2)&= {\rm
max}\{\varepsilon_i(b_1),\varepsilon_i(b_2)-\langle {\rm
wt}(b_1),h_i\rangle\}, \\
\varphi_i(b_1\otimes b_2)&= {\rm max}\{\varphi_i(b_1)+\langle {\rm
wt}(b_2),h_i\rangle,\varphi_i(b_2)\},\\
{\te}_i(b_1\otimes b_2)&=
\begin{cases}
{\te}_i b_1 \otimes b_2, & \text{if $\varphi_i(b_1)\geq \varepsilon_i(b_2)$}, \\
b_1\otimes {\te}_i b_2, & \text{if
$\varphi_i(b_1)<\varepsilon_i(b_2)$},
\end{cases}\\
{\tf}_i(b_1\otimes b_2)&=
\begin{cases}
{\tf}_i b_1 \otimes b_2, & \text{if  $\varphi_i(b_1)>\varepsilon_i(b_2)$}, \\
b_1\otimes {\tf}_i b_2, & \text{if $\varphi_i(b_1)\leq
\varepsilon_i(b_2)$},
\end{cases}
\end{align*}
\noindent for $i\in I$. Here we assume that ${\bf 0}\otimes
b_2=b_1\otimes {\bf 0}={\bf 0}$.}

A morphism
$\psi : B_1 \rightarrow B_2$ is a map from $B_1\cup\{{\bf 0}\}$ to
$B_2\cup\{{\bf 0}\}$ such that
\begin{itemize}
\item[(1)] $\psi(\bf{0})=\bf{0}$,

\item[(2)] ${\rm wt}(\psi(b))={\rm wt}(b)$,
$\varepsilon_i(\psi(b))=\varepsilon_i(b)$, and
$\varphi_i(\psi(b))=\varphi_i(b)$ when $\psi(b)\neq \bf{0}$,

\item[(3)] $\psi(\te_i b)=\te_i\psi(b)$ for $b\in B_1$ such that $\psi(b)\neq \bf{0}$ and
$\psi(\te_i b)\neq \bf{0}$,

\item[(4)] $\psi(\tf_i
b)=\tf_i\psi(b)$ for $b\in B_1$ such that $\psi(b)\neq \bf{0}$ and
$\psi(\tf_i b)\neq \bf{0}$.
\end{itemize}
We call $\psi$ an embedding and $B_1$ a subcrystal of
$B_2$ when $\psi$ is injective.

The dual crystal $B^\vee$ of a crystal $B$ is defined
to be the set $\{\,b^\vee\,|\,b\in B\,\}$ with ${\rm
wt}(b^\vee)=-{\rm wt}(b)$, $\varepsilon_i(b^\vee)=\varphi_i(b)$,
$\varphi_i(b^\vee)=\varepsilon_i(b)$, $\te_i(b^\vee)=(\tf_i b)^\vee$ and $\tf_i(b^\vee)=\left(\te_i b \right)^\vee$ for
$b\in B$ and $i\in I$. (Here we assume that ${\bf 0}^\vee={\bf 0}$.)
For $\mu\in P$, let $T_\mu=\{t_\mu\}$ be a crystal, where ${\rm wt}(t_\mu)=\mu$, and $\varphi_i(t_\mu)=-\infty$ for all $i\in I$.
Given $b_i $ in crystals $B_i $ ($i=1,2$), we write $b_1 \equiv b_2$ if there is an isomorphism of crystals $C(b_1) \rightarrow C(b_2)$ mapping $b_1$ to $b_2$, where $C(b_i)$ denotes the connected component of $b_i$ in $B_i$.

Let $q$ be an indeterminate. Let $U_q(\g)$ be the quantized enveloping algebra of $\g$, which is an associative $\mathbb{Q}(q)$-algebra with $1$ generated by $e_i$, $f_i$, and $q^h$ for $i\in I$ and $h\in P^\vee$. 
Let $U^-_q(\g)$ be the negative part of $U_q(\g)$, the subalgebra generated by $f_i$ for $i\in I$. We denote by $B(\infty)$ the crystal associated to $U^-_q(\g)$, and by $B(\Lambda)$ the crystal associated to an irreducible highest weight $U_q(\g)$-module with highest weight
$\Lambda\in P^+$.

Let $\ast$ be the $\mathbb{Q}(q)$-linear anti-automorphism of $U_q(\g)$ such that 
$e_i^\ast =e_i$, $f_i^\ast=f_i$, and $(q^h)^\ast =q^{-h}$ for $i\in I$ and $h\in P$. It is shown in \cite{Kas91,Kas93} that $\ast$ induces a bijection on $B(\infty)$. For $i\in I$, we define $\te_i^\ast =\ast \circ \te_i \circ \ast$ and $\tf_i^\ast =\ast \circ \tf_i \circ \ast$ on $B(\infty)$.

\subsection{PBW basis and Lusztig data}\label{sec:PBW crystal}
Suppose that $\g$ is a complex finite-dimensional simple Lie algebra.  
Let us review the PBW basis of $U^-_q(\g)$ \cite{Lu90,Lu90-2,S94}.
Let $w_0$ be the longest element in $W$ of length $N$, and let $R(w_0)=\{\,(i_1,\ldots,i_N)\,|\,w_0=s_{i_1}\ldots s_{i_N}\,\}$ be the set of reduced expressions of $w_0$.


We put  $[m]_q=\frac{q^m-q^{-m}}{q-q^{-1}}$, $[m]_q!=[1]_q[2]_q\cdots [m]_q$ for $m\in\N$, and $[0]_q!=1$. Let $t_i=q^{h_i}$, $e_i^{(m)}=e_i^m/[m]_q!$, and $f_i^{(m)}=f_i^m/[m]_q!$  for $m\in\Z_+$ and $i\in I$.

For $i\in I$, let $T_i$ be the $\mathbb{Q}(q)$-algebra automorphism of $U$, which is equal to {$T''_{i,1}$} in \cite{Lu93}, given by
{\allowdisplaybreaks
\begin{align*}
T_i(t_j)&=t_j t_i^{-a_{ij}}, \\
T_i(e_j)&=  
\begin{cases}
{-f_it_i}, & \text{if $j=i$}, \\
\sum_{k+l=-a_{ij}} (-1)^kq_i^{-l}e_i^{(k)}e_je_i^{(l)}, & \text{if $j\neq i$},
\end{cases}\\
T_i(f_j)&=  
\begin{cases}
{-t_i^{-1}e_i}, & \text{if $j=i$}, \\
\sum_{k+l=-a_{ij}} (-1)^kq_i^{k}f_i^{(l)}f_jf_i^{(k)}, & \text{if $j\neq i$},
\end{cases}
\end{align*}}
for $j\in I$, where $a_{ij}=\langle \alpha_j,h_i \rangle$ and $q_i=q^{d_i}$.

For ${\bi}=(i_1,\ldots, i_N)\in R(w_0)$ and ${\bf c}=(c_1,\ldots,c_N)\in\Z_+^N$, let
\begin{equation}\label{eq:PBW vector}
\begin{split}
b _{\bi}(\bf c)=&
f^{(c_{1})}_{i_1}T_{i_1}(f^{(c_{2})}_{i_2})\cdots T_{i_{1}}T_{i_2}\cdots T_{i_{N-1}}(f^{(c_{N})}_{i_{N}}).
\end{split}
\end{equation}
Then the set $B_{\bf i}:=\{\,b _{\bi}({\bf c})\,|\,{\bf c}\in\Z_+^{N}\,\}$ is a $\mathbb{Q}(q)$-basis of $U^-_q(\g)$ called a {\it PBW basis}.

Let $A_0$ denote the subring of $\mathbb{Q}(q)$ consisting
of rational functions regular at $q=0$.
The $A_0$-lattice of $U^-_q(\g)$ generated by $B_{\bf i}$ is independent of the choice of ${\bf i}$, which we denote by $L(\infty)$. 
For $i\in I$, let $\te_i$ and $\tf_i$ denote the modified Kashiwara operators on $U^-_q(\g)$ given by
\begin{equation*}
\te_i x =\sum_{k\geq 1}f_i^{(k-1)}x_k,\quad\quad \tf_i x =\sum_{k\geq 0}f_i^{(k+1)}x_k,
\end{equation*}
for $x=\sum_{k\geq 0}f_i^{(k)}x_k$, where $x_k\in T_i(U^-_q(\g))\cap U^-_q(\g)$ for $k\geq 0$. 
Then $L(\infty)$ is invariant under $\te_i$, $\tf_i$.
If $\pi : L(\infty) \rightarrow L(\infty)/q L(\infty)$ is a canonical projection, then $\pi(B_{\bf i})$ is a $\mathbb{Q}$-basis of  $L(\infty)/qL(\infty)$, which is also independent of the choice of ${\bf i}$, and $\pi(B_{\bf i})\cup \{0\}$ is invariant under $\te_i$, $\tf_i$ for $i\in I$. Hence $\pi(B_{\bf i})$ becomes a $\g$-crystal with respect to $\te_i$ and $\tf_i$, which is isomorphic to $B(\infty)$. We identify ${\bf B}_{\bi}:=\Z_+^N$ with a crystal $\pi(B_{\bi})$ under the map ${\bf c}\mapsto b_{\bi}({\bf c})$, and call ${\bf c}\in \B_{\bi}$ an {\it $\bi$-Lusztig datum}.

\subsection{Kashiwara-Nakashima tableaux}\label{subsec:KN}
Let us review the notion of Kashiwara-Nakashima tableaux (KN tableaux for short) \cite{KashNaka}. 

Fix a positive integer $n\geq 2$ throughout the paper.
From now on, we assume that $\g={\mf b}_n$ of type $B_n$ or $\g={\mf c}_n$ of type $C_n$ with $I=\{\,1,\ldots,n\,\}$. 
Let $\mf l$ be the Levi subalgebra of $\mf g$ corresponding to $\{\,\alpha_1,\ldots,\alpha_{n-1}\,\}$.

We assume that the weight lattice of $\g$ is $P=\bigoplus_{i=1}^n\Z\epsilon_i$, where $\{\,\epsilon_i\,|\,1\leq i\leq n\,\}$ is an orthonormal basis with respect to a symmetric bilinear form $(\,,\,)$. We denote by $\omega_i$ the $i$th fundamental weight for $i\in I$. For $\la\in\cP_n$, we put
\begin{equation*}
\omega_\la=\la_1\epsilon_1+\cdots+\la_n\epsilon_n.
\end{equation*}

We consider a linearly ordered set 
\begin{equation*}
\begin{split}
\J_n&=\{\,1<2<\cdots<n<0<\ov{n}<\cdots<\ov{2}<\ov{1}\,\}. \\
\end{split}
\end{equation*}
We put $\J^\times_n=\J_n\setminus\{0\}$, $[n]=\{\,1,\ldots,n\,\}$, and $[\ov{n}]=\{\,\ov{n},\ldots,\ov{1}\,\}$. We assume that $\J_n$ is $\Z_2$-graded where $(\J_n)_0=\J_n^\times$ and $(\J_n)_1=\{0\}$.

\subsubsection{Type $C_n$} Suppose that $\g={\mf c}_n$ with Dynkin diagram
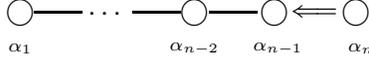
\begin{center} 
\setlength{\unitlength}{0.19in}
\begin{picture}(20,3.5)
\put(5.6,2){\makebox(0,0)[c]{$\bigcirc$}}
\put(12.6,2){\makebox(0,0)[c]{$\bigcirc$}}
\put(10.4,2){\makebox(0,0)[c]{$\bigcirc$}}
\put(14.85,2){\makebox(0,0)[c]{$\bigcirc$}}
\put(6,2){\line(1,0){1.3}} \put(8.7,2){\line(1,0){1.3}} \put(10.82,2){\line(1,0){1.3}}
%
\put(13.7,2){\makebox(0,0)[c]{$\Longleftarrow$}}

\put(8,1.95){\makebox(0,0)[c]{$\cdots$}}
\put(5.6,1){\makebox(0,0)[c]{\tiny ${\alpha}_1$}}
\put(12.7,1){\makebox(0,0)[c]{\tiny ${\alpha}_{n-1}$}}
\put(10.4,1){\makebox(0,0)[c]{\tiny ${\alpha}_{n-2}$}}
\put(15,1){\makebox(0,0)[c]{\tiny ${\alpha}_n$}}

\end{picture}
\end{center}\vskip -3mm 
where $\alpha_i=\epsilon_i-\epsilon_{i+1}$ for $1\leq i\leq n-1$, and $\alpha_n=2\epsilon_n$. We have $\omega_i=\epsilon_1+\cdots+\epsilon_i$ for $1\leq i\leq n$, and $P^+=\{\,\omega_\la\,|\,\la\in \cP_n\,\}$. 
Note that $\omega_\la=\sum_{i\geq 1}\omega_{\la'_i}$ for $\la\in \cP_n$.

\begin{df}\label{def:KN-C}{\rm
For $\la\in \cP_n$, let ${\bf KN}^{\mf c_n}_{\la}$ be the set of $T\in SST_{\J_n^\times}(\la^\pi)$ satisfying  
\begin{itemize}
\item[(${\mf c}$-1)] if $T({i_1,j})=\ov{a}$ and $T({i_2,j})=a$ for some $a$ and $1\leq i_1< i_2\leq \la'_j$, then we have $i_1+(\la'_j-i_2+1)\leq a$,

\item[(${\mf c}$-2)] if either $T({p,j})=\ov{a}$, $T({q,j})=\ov{b}$, $T({r,j})=b$, $T({s,j+1})=a$ or 
 $T({p,j})=\ov{a}$, $T({q,j+1})=\ov{b}$, $T({r,j+1})=b$, $T({s,j+1})=a$ for some $1\leq a\leq b\leq n$, and $p\leq q<r\leq s$, then we have $(q-p)+(s-r)<b-a$,

\end{itemize}
where $T({i,j})$ denotes the entry of $T$ in the $i$th row from the bottom and the $j$th column from the right. We call ${\bf KN}^{\mf c_n}_{\la}$ the set of {\em KN tableaux of type $C_n$ of shape $\la$.}}
\end{df}

The set ${\bf KN}^{\mf c_n}_{(1)}$ has a $\mf c_n$-crystal structure such that ${\bf KN}^{\mf c_n}_{(1)}\cong B(\omega_1)$, where
\begin{equation*}
\xymatrixcolsep{2pc}\xymatrixrowsep{0pc}\xymatrix{
\boxed{1}
\, \ar@{->}[r]^{1} & \, \boxed{2} \, \ar@{->}[r]^{2} & \, \cdots \, \ar@{->}[r]^{n-1}  
& \, \boxed{n}\,  \ar@{->}[r]^{n} & \, \boxed{\ov{n}}\, \ar@{->}[r]^{n-1} &   \cdots   \ar@{->}[r]^{2} & \boxed{\ov{2}}\, \ar@{->}[r]^{1} & \, \boxed{\ov{1}}}
\end{equation*}
with ${\rm wt}(\,
\resizebox{.02\hsize}{!}{\def\lr#1{\multicolumn{1}{|@{\hspace{.6ex}}c@{\hspace{.6ex}}|}{\raisebox{-.25ex}{$#1$}}}\raisebox{-.65ex}
{$\begin{array}[b]{c}
\cline{1-1} 
\lr{{i}} \\
\cline{1-1} 
\end{array}$}}
\,)=\epsilon_i$ and 
${\rm wt}(\,
\resizebox{.02\hsize}{!}{\def\lr#1{\multicolumn{1}{|@{\hspace{.6ex}}c@{\hspace{.6ex}}|}{\raisebox{-.25ex}{$#1$}}}\raisebox{-.65ex}
{$\begin{array}[b]{c}
\cline{1-1} 
\lr{\ov{i}} \\
\cline{1-1} 
\end{array}$}}
\,)=-\epsilon_i$ for $1\leq i\leq n$. 
For $\la\in \cP_n$ and $i\in I$, we define $\te_i$ and $\tf_i$ on ${\bf KN}^{\mf c_n}_{\la}$ under the identification of $T\in {\bf KN}^{\mf c_n}_{\la}$ with 
$\resizebox{.04\hsize}{!}{\def\lr#1{\multicolumn{1}{|@{\hspace{.6ex}}c@{\hspace{.6ex}}|}{\raisebox{-.1ex}{$#1$}}}\raisebox{-.9ex}
{$\begin{array}[b]{c}
\cline{1-1} 
\lr{w_1} \\
\cline{1-1} 
\end{array}$}}\otimes \cdots \otimes 
\resizebox{.04\hsize}{!}{\def\lr#1{\multicolumn{1}{|@{\hspace{.6ex}}c@{\hspace{.6ex}}|}{\raisebox{-.1ex}{$#1$}}}\raisebox{-.9ex}
{$\begin{array}[b]{c}
\cline{1-1} 
\lr{w_r} \\
\cline{1-1} 
\end{array}$}}\in ({\bf KN}^{\mf c_n}_{(1)})^{\otimes r}$ when $w(T)=w_1\cdots w_r$. 
Then ${\bf KN}^{\mf c_n}_{\la}$ is a $\mf c_n$-crystal with respect to $\te_i$ and $\tf_i$ for $i\in I$, and
\begin{equation}\label{eq:KN=B(la):c}
{\bf KN}^{\mf c_n}_{\la}\cong B(\omega_\la).
\end{equation}

\begin{rem}{\rm
Note that a KN tableau in \cite{KashNaka} is of shape $\la$. One may consider ${\bf KN}^{\mf c_n}_{\la}$ here as the dual crystal of the KN tableaux of shape $\la$ in \cite{KashNaka}. This immediately implies \eqref{eq:KN=B(la):c}. The same holds for the case of type $B_n$ in the next subsection.
}
\end{rem}

\subsubsection{Type $B_n$}
Suppose that $\g={\mf b}_n$ with Dynkin diagram
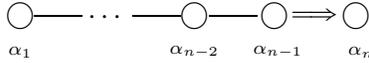
\begin{center} 
\setlength{\unitlength}{0.19in}
\begin{picture}(20,3.5)
\put(5.6,2){\makebox(0,0)[c]{$\bigcirc$}}
\put(12.6,2){\makebox(0,0)[c]{$\bigcirc$}}
\put(10.4,2){\makebox(0,0)[c]{$\bigcirc$}}
\put(14.85,2){\makebox(0,0)[c]{$\bigcirc$}}
\put(6,2){\line(1,0){1.3}} \put(8.7,2){\line(1,0){1.3}} \put(10.82,2){\line(1,0){1.3}}
%
\put(13.7,2){\makebox(0,0)[c]{$\Longrightarrow$}}

\put(8,1.95){\makebox(0,0)[c]{$\cdots$}}
\put(5.6,1){\makebox(0,0)[c]{\tiny ${\alpha}_1$}}
\put(12.7,1){\makebox(0,0)[c]{\tiny ${\alpha}_{n-1}$}}
\put(10.4,1){\makebox(0,0)[c]{\tiny ${\alpha}_{n-2}$}}
\put(15,1){\makebox(0,0)[c]{\tiny ${\alpha}_n$}}

\end{picture}
\end{center}\vskip -3mm 
where $\alpha_i=\epsilon_i-\epsilon_{i+1}$ for $1\leq i\leq n-1$, and $\alpha_n=\epsilon_n$. 
We have $\omega_i=\epsilon_1+\cdots+\epsilon_i$ for $1\leq i\leq n-1$, $\omega_n=(\epsilon_1+\cdots+\epsilon_n)/2$, and $P^+=\{\,\omega_\la, \omega_\la+\omega_n\,|\,\la\in \cP_n\,\}$.

\begin{df}{\rm
For $\la\in \cP_n$, let ${\bf KN}^{\mf b_n}_{\la}$ be the set of $T\in SST_{\J_n}(\la^\pi)$ satisfying
\begin{itemize}
\item[(${\mf b}$-1)] if $T({i_1,j})=\ov{a}$ and $T({i_2,j})=a$ for some $a$ and $i_1< i_2$, then we have $i_1+(\la'_j-i_2+1)\leq a$,

\item[(${\mf b}$-2)] if either $T({p,j})=\ov{a}$, $T({q,j})=\ov{b}$, $T({r,j})=b$, $T({s,j+1})=a$ or 
 $T({p,j})=\ov{a}$, $T({q,j+1})=\ov{b}$, $T({r,j+1})=b$, $T({s,j+1})=a$ for some $1\leq a\leq b< n$, and $p\leq q<r\leq s$, then we have $(q-p)+(s-r)<b-a$,
 
\item[(${\mf b}$-3)] if either $T({p,j})=\ov{a}$, $T({s,j+1})=a$, $T({q,j}), T({r,j})\in \{n,0,\ov{n}\}$ or $T({p,j})=\ov{a}$, $T({s,j+1})=a$, $T({q,j+1}), T({r,j+1})\in \{n,0,\ov{n}\}$ for some $1\leq a< n$ and $p\leq q<r=q+1\leq s$, then we have $(q-p)+(s-r)<b-a$,
 
\item[(${\mf b}$-4)] there is no $p<q$ such that $T(p,j)\in\{0,\ov{n}\}$, $T(q,j+1) \in \{n,0\}$ for some $j$.
   
\end{itemize}}
\end{df}

The set ${\bf KN}^{\mf b_n}_{(1)}$ has a $\mf b_n$-crystal structure such that ${\bf KN}^{\mf b_n}_{(1)}\cong B(\omega_1)$, where 
\begin{equation*}
\xymatrixcolsep{1.8pc}\xymatrixrowsep{0pc}\xymatrix{
\boxed{1}\, \ar@{->}[r]^{1} & \, \boxed{2} \, \ar@{->}[r]^{2} & \, \cdots \, \ar@{->}[r]^{n-1}  
& \, \boxed{n}\,  \ar@{->}[r]^{n} & \, \boxed{0}\,  \ar@{->}[r]^{n} & \, \boxed{\ov{n}}\, \ar@{->}[r]^{n-1} &   \cdots   \ar@{->}[r]^{2} & \boxed{\ov{2}}\, \ar@{->}[r]^{1} & \, \boxed{\ov{1}}}
\end{equation*}
with  ${\rm wt}(\,
\resizebox{.02\hsize}{!}{\def\lr#1{\multicolumn{1}{|@{\hspace{.6ex}}c@{\hspace{.6ex}}|}{\raisebox{-.25ex}{$#1$}}}\raisebox{-.65ex}
{$\begin{array}[b]{c}
\cline{1-1} 
\lr{{i}} \\
\cline{1-1} 
\end{array}$}}
\,)=\epsilon_i$,
${\rm wt}(\,
\resizebox{.02\hsize}{!}{\def\lr#1{\multicolumn{1}{|@{\hspace{.6ex}}c@{\hspace{.6ex}}|}{\raisebox{-.25ex}{$#1$}}}\raisebox{-.65ex}
{$\begin{array}[b]{c}
\cline{1-1} 
\lr{\ov{i}} \\
\cline{1-1} 
\end{array}$}}
\,)=-\epsilon_i$ for $1\leq i\leq n$, and ${\rm wt}(0)=0$.  
As in ${\bf KN}^{\mf c_n}_{\la}$, one can define $\te_i$ and $\tf_i$ ($i\in I$) on ${\bf KN}^{\mf b_n}_{\la}$ for $\la\in \cP_n$. Then ${\bf KN}^{\mf b_n}_{\la}$ is a $\mf b_n$-crystal with respect to $\te_i$ and $\tf_i$ for $i\in I$, and
\begin{equation}\label{eq:KN=B(la):b}
{\bf KN}^{\mf b_n}_{\la}\cong B(\omega_\la).
\end{equation}

\begin{df}{\rm 
Let ${\bf KN}^{\rm sp}$ be the set of $\J_n^\times$-semistandard tableaux $T$ of single-columned shape with height $n$ and half width such that $T$ does not contain $a$ and $\ov{a}$ at the same time for any $1\leq a\leq n$. 

Let $\cP_{n}^{\rm sp}=\{\,\mu+\sigma_n\,|\,\mu\in\cP_n\,\}$, where $\sigma_n=(1/2,\ldots,1/2)$ ($n$ times).
For $\la\in \cP_{n}^{\rm sp}$ with $\la=\mu+\sigma_n$, let ${\bf KN}^{\mf b_n}_{\la}$ to be the set of $\J_n^\times$-semistandard tableaux such that 
\begin{itemize}
\item[(1)] the first column of $T$ (from the right) is a tableau in ${\bf KN}^{\rm sp}$, 

\item[(2)] the subtableau consisting of the other columns is a tableau in ${\bf KN}^{\mf b_n}_{\mu}$, 

\item[(3)] the first two columns of $T$ (from the right) satisfy (${\mf b}$-2)-(${\mf b}$-4).
\end{itemize}}
\end{df}

For $T\in {\bf KN}^{\rm sp}$, put ${\rm wt}(T)=\sum_{i=1}^n c_i(\epsilon_i/2)$, where $c_i=1$ (resp. $c_i=-1$) if $T$ has $i$ (resp. $\ov{i}$). 
For $i\in I$, define $\te_iT$ to be the unique tableau in ${\bf KN}^{\rm sp}$ such that ${\rm wt}(\te_i T)=s_i({\rm wt}(T))$ if $c_i=-1$, and ${\bf 0}$ if $c_i =1$.  
Then ${\bf KN}^{\rm sp}$ is a crystal isomorphic to $B(\omega_n)$. 
For $\la\in \cP_{n}^{\rm sp}$ with $\la=\mu+\sigma_n$, we define a crystal structure on ${\bf KN}^{\mf b_n}_{\la}$ by regarding ${\bf KN}^{\mf b_n}_{\la}\subset {\bf KN}^{\rm sp}\otimes{\bf KN}^{\mf b_n}_{\mu}$. 
Then we have
\begin{equation}\label{eq:KN=B(la):b:sp}
{\bf KN}^{\mf b_n}_{\la}\cong B(\omega_\la),
\end{equation}
where $\omega_\la=\omega_\mu + \omega_n$.
For $\la\in \cP_n\cup\cP_n^{\rm sp}$, we call ${\bf KN}^{\mf b_n}_{\la}$ the set of {\em KN tableaux of type $B_n$ of shape $\la$}.

\section{Spinor model for classical crystals}\label{sec:spinor}

In this section, we recall the spinor model of classical crystals \cite{K15,K16} and then construct an isomorphism $\Psi_\la$ from the crystal of KN tableaux to the spinor model. 
\subsection{Combinatorics of two-columned  skew tableaux}\label{subsec:PSST-1} 
For a single-columned tableau $U$, let us denote by ${\rm ht}(U)$ the height of $U$, and by $U(i)$ the $i$th entry of $U$ from the bottom for $i\geq 1$.

For $a,b,c\in\Z_+$, let $\lambda(a,b,c)=(2^{b+c},1^a)/(1^b)$, which is a skew Young diagram of two-columned shape.  
Let $T$ be a tableau of shape $\lambda(a,b,c)$. We denote by $T^{\tt L}$ and $T^{\tt R}$ the left and right columns of $T$, respectively.  We also denote by $T^{\tt tail}$ the subtableau of $T$ corresponding to the tail of $\la(a,b,c)$, a lower single column of height $a$, and by $T^{\tt body}$ the subtableau of $T$ above $T^{\tt tail}$. For example,

$$T=\resizebox{.055\hsize}{!}
{\def\lr#1{\multicolumn{1}{|@{\hspace{.75ex}}c@{\hspace{.75ex}}|}{\raisebox{-.04ex}{$#1$}}}\raisebox{-.6ex}
{$\begin{array}{cc}
\cline{2-2}
&\lr{1}\\
\cline{1-1}\cline{2-2}
\lr{2} &\lr{4}\\
\cline{1-1}\cline{2-2}
\lr{4}&\lr{6}\\
\cline{1-1}\cline{2-2}
\lr{5} & \lr{8}\\
\cline{1-1}\cline{2-2}
\lr{6} \\
\cline{1-1}
\lr{7} \\
\cline{1-1}
\end{array}$}}\ \in SST_{[8]}(\la(2,1,3))
$$\vskip 2mm
$$
T^{\tt L}=\resizebox{.048\hsize}{!}
{\def\lr#1{\multicolumn{1}{|@{\hspace{.75ex}}c@{\hspace{.75ex}}|}{\raisebox{-.04ex}{$#1$}}}\raisebox{-.6ex}
{$\begin{array}{cc}
 & \\
\cline{1-1} 
\lr{2} & \\
\cline{1-1} 
\lr{4}& \\
\cline{1-1} 
\lr{5} &  \\
\cline{1-1} 
\lr{6} \\
\cline{1-1}
\lr{7} \\
\cline{1-1}
\end{array}$}}\ \ \ 
T^{\tt R}=\resizebox{.027\hsize}{!}
{\def\lr#1{\multicolumn{1}{|@{\hspace{.75ex}}c@{\hspace{.75ex}}|}{\raisebox{-.04ex}{$#1$}}}\raisebox{-.6ex}
{$\begin{array}{cc}
\cline{1-1}
\lr{1}\\
\cline{1-1}
\lr{4}\\
\cline{1-1}
\lr{6}\\
\cline{1-1}
\lr{8}\\
\cline{1-1}
 \\
 \\ 
\end{array}$}}\ \ \ \ \ \  
T^{\tt body}=\resizebox{.055\hsize}{!}
{\def\lr#1{\multicolumn{1}{|@{\hspace{.75ex}}c@{\hspace{.75ex}}|}{\raisebox{-.04ex}{$#1$}}}\raisebox{-.6ex}
{$\begin{array}{cc}
\cline{2-2}
&\lr{1}\\
\cline{1-1}\cline{2-2}
\lr{2} &\lr{4}\\
\cline{1-1}\cline{2-2}
\lr{4}&\lr{6}\\
\cline{1-1}\cline{2-2}
\lr{5} & \lr{8}\\
\cline{1-1}\cline{2-2}
\\
\\
\end{array}$}}\ \ \ \ \
T^{\tt tail}=\ \resizebox{.026\hsize}{!}
{\def\lr#1{\multicolumn{1}{|@{\hspace{.75ex}}c@{\hspace{.75ex}}|}{\raisebox{-.04ex}{$#1$}}}\raisebox{-.6ex}
{$\begin{array}{c}
\\
\\
\\
\\
\cline{1-1}
\lr{6} \\
\cline{1-1}
\lr{7} \\
\cline{1-1}
\end{array}$}}\ \ \ \ \
$$\vskip 3mm

Let $\A$ be a linearly ordered set.
Let $T\in SST_\A(\la(a,b,c))$ be given. One may slide down $T^{\tt R}$ by $k$ positions for $0\leq k\leq \min\{a,b\}$ to have a (not necessarily semistandard) tableau $T'$ of shape $\lambda(a-k,b-k,c+k)$. We define ${\mf r}_T$ to be the maximal $k$ such that $T'$ is semistandard. Note that ${\mf r}_T=0$ when $a=0$ or $b=0$.

\begin{df}\label{def:mc X}
{\rm
For $T\in SST_\A(\la(a,b,c))$ with ${\mf r}_T=0$, we define  
\begin{itemize}
\item[(1)] $\mc E T$ to be tableau in $SST_{\A}(\la(a-1,b+1,c))$ obtained from $T$ by applying jeu de taquin sliding to the position below the bottom of $T^{\tt R}$, when $a>0$,

\item[(2)] $\mc F T$ to be tableau in $SST_{\A}(\la(a+1,b-1,c))$ obtained from $T$ by applying jeu de taquin sliding to the position above the top of $T^{\tt L}$, when $b>0$.

\end{itemize}
}
\end{df}

We assume that ${\mc E}T={\bf 0}$ when $a=0$, and ${\mc F}T={\bf 0}$ when $b=0$. Here ${\bf 0}$ is a formal symbol. By definition one can check the following (cf.\cite{La}).
\begin{lem}\label{lem:cal E and F}
Under the above hypothesis,
\begin{itemize}
\item[(1)] ${\mf r}_{\mc E T}=0$ and ${\mf r}_{\mc F T}=0$ whenever ${\mc E}T$ and ${\mc F}T$ are defined

\item[(2)] $\{\,{\mc E}^k T\,|\,0\leq k\leq a\,\}\cup \{\,{\mc F}^l T\,|\,0\leq l\leq b\,\}$ forms a regular $\mf{sl}_2$-crystal with respect to $\mc E$ and $\mc F$.
\end{itemize}
\end{lem}

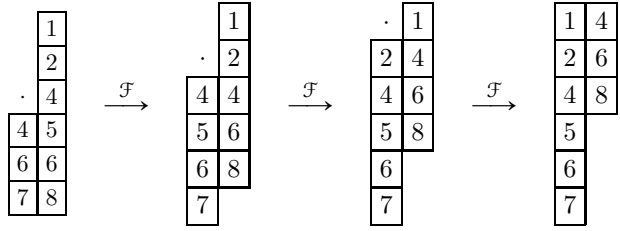
\begin{ex}{\rm

$$
\resizebox{.055\hsize}{!}
{\def\lr#1{\multicolumn{1}{|@{\hspace{.75ex}}c@{\hspace{.75ex}}|}{\raisebox{-.04ex}{$#1$}}}\raisebox{-.6ex}
{$\begin{array}{cc}
\cline{2-2}
&\lr{1}\\
\cline{2-2}
 &\lr{2}\\
 \cline{2-2}
\cdot &\lr{4}\\
\cline{1-1}\cline{2-2}
\lr{4} & \lr{5}\\
\cline{1-1}\cline{2-2}
\lr{6}& \lr{6} \\
\cline{1-2}
\lr{7} & \lr{8} \\
\cline{1-2}
\end{array}$}}
\quad  \stackrel{\mc{F}}{\longrightarrow} \quad 
\resizebox{.06\hsize}{!}
{\def\lr#1{\multicolumn{1}{|@{\hspace{.75ex}}c@{\hspace{.75ex}}|}{\raisebox{-.04ex}{$#1$}}}\raisebox{-.6ex}
{$\begin{array}{cc}
\cline{2-2}
&\lr{1}\\
\cline{2-2}
\cdot &\lr{2}\\
\cline{1-1}\cline{2-2}
\lr{4}&\lr{4}\\
\cline{1-1}\cline{2-2}
\lr{5} & \lr{6}\\
\cline{1-1}\cline{2-2}
\lr{6}& \lr{8} \\
\cline{1-2}
\lr{7} \\
\cline{1-1}
\end{array}$}}
\quad  \stackrel{\mc{F}}{\longrightarrow} \quad 
\resizebox{.06\hsize}{!}
{\def\lr#1{\multicolumn{1}{|@{\hspace{.75ex}}c@{\hspace{.75ex}}|}{\raisebox{-.04ex}{$#1$}}}\raisebox{-.6ex}
{$\begin{array}{cc}
\cline{2-2}
\cdot &\lr{1}\\
\cline{1-1}\cline{2-2}
\lr{2} &\lr{4}\\
\cline{1-1}\cline{2-2}
\lr{4}&\lr{6}\\
\cline{1-1}\cline{2-2}
\lr{5} & \lr{8}\\
\cline{1-1}\cline{2-2}
\lr{6} \\
\cline{1-1}
\lr{7} \\
\cline{1-1}
\end{array}$}}
\quad  \stackrel{\mc{F}}{\longrightarrow} \quad 
\resizebox{.06\hsize}{!}
{\def\lr#1{\multicolumn{1}{|@{\hspace{.75ex}}c@{\hspace{.75ex}}|}{\raisebox{-.04ex}{$#1$}}}\raisebox{-.6ex}
{$\begin{array}{cc}
\cline{1-2}
\lr{1} & \lr{4}\\
\cline{1-1}\cline{2-2}
\lr{2} &\lr{6}\\
\cline{1-1}\cline{2-2}
\lr{4}&\lr{8}\\
\cline{1-1}\cline{2-2}
\lr{5} &  \\
\cline{1-1} 
\lr{6} \\
\cline{1-1}
\lr{7} \\
\cline{1-1}
\end{array}$}}
$$\vskip 2mm
}
\end{ex}

We define ${\mc E}$ and ${\mc F}$ on a pair of single columned tableaux in general.
For a pair $(U,V)\in SST_{\A}(1^u)\times SST_{\A}(1^v)$ ($u,v\geq 0$), let $T$ be a unique tableau in $SST_{\A}(\la(u-k,v-k,k))$ for some $0\leq k\leq \min\{u,v\}$ such that $(T^{\tt L},T^{\tt R})=(U,V)$ and ${\mf r}_T=0$. Then we define 
\begin{equation}\label{eq:E and F via jdt}
\begin{split}
\mc{X}(U,V)&=
\begin{cases}
((\mc{X}T)^{\tt L},(\mc{X}T)^{\tt R}), &\text{if $\mc{X}T\neq {\bf 0}$},\\
{\bf 0}, &\text{if $\mc{X}T={\bf 0}$},
\end{cases} \quad\quad (\mc{X}=\mc{E}, \mc{F}),
\end{split}
\end{equation}
where ${\mc X}T$ is defined in Definition \ref{def:mc X}.
In particular, when $U$ is empty, we have $\mc F^k(U,V)=(U_k,V_k)$ for $0\leq k\leq {\rm ht}(V)$, where
$U_k$ is the subtableau of $V$ consisting of the first $k$ entries from the bottom and $V_k$ is the complement of $U_k$ in $V$.

The above $\mf{sl}_2$-crystal structure will be used in Section \ref{sec:main}.

\subsection{Spinor model for classical crystals}\label{subsec:PSST}
In this subsection, we review another combinatorial model of $B(\la)$ in types $B_n$ and $C_n$ introduced in \cite{K15}.
We will follow the convention in \cite{K15-1},
where the definition of this model is slightly different from \cite{K15}.

For $0\leq a< n$, let
\begin{equation*}
{\bf T}^{{\mf g}}(a)=
\begin{cases}
\{\,T\,|\,T\in SST_{[\ov{n}]}(\la(a,0,c)), \ c\in \Z_+ \,\}, & \text{if $\mf g=\mf c_n$},\\
\{\,T\,|\,T\in SST_{[\ov{n}]}(\la(a,b,c)), \  (b,c)\in \Z_+\times  \Z_+, \  {\mf r}_T=0\,\},
 & \text{if $\mf g=\mf b_n$},
\end{cases}
\end{equation*}
and 
\begin{equation*}
{\bf T}^{\rm sp}= \bigsqcup_{0\leq a\leq n} SST_{[\ov{n}]}((1^a)).
\end{equation*}

Suppose that $\B$ is either ${\bf T}^{{\mf g}}(a)$ or ${\bf T}^{\rm sp}$.
First we regard $\B$ as an $\mf l$-crystal since $[n]$ is the crystal of the natural representation of $\mf l$, and  $[\ov{n}]$ is its dual (cf.\cite{KashNaka}). 
We denote by ${\rm wt}_{\mf l}(T)$ the ${\mf l}$-weight of $T\in \B$.
Next we define $\te_n$ on $\B$ as follows ($\tf_n$ is defined in a similar way):
\begin{itemize}
\item[(1)] if $T\in {\bf T}^{\mf c_n}(a)$, then $\te_n T$ is the tableau obtained from $T$ by removing a domino
$\resizebox{.05\hsize}{!}{\def\lr#1{\multicolumn{1}{|@{\hspace{.6ex}}c@{\hspace{.6ex}}|}{\raisebox{-.25ex}{$#1$}}}\raisebox{-.65ex}
{$\begin{array}[b]{cc}
\cline{1-1}\cline{2-2}
\lr{\ov{n}}&\lr{\ov{n}}\\
\cline{1-1}\cline{2-2}
\end{array}$}}$\ if $T$ has $\resizebox{.05\hsize}{!}{\def\lr#1{\multicolumn{1}{|@{\hspace{.6ex}}c@{\hspace{.6ex}}|}{\raisebox{-.25ex}{$#1$}}}\raisebox{-.65ex}
{$\begin{array}[b]{cc}
\cline{1-1}\cline{2-2}
\lr{\ov{n}}&\lr{\ov{n}}\\
\cline{1-1}\cline{2-2}
\end{array}$}}$ on its top, and ${\bf 0}$, otherwise,

\item[(2)] if $T\in {\bf T}^{\rm sp}$, then
$\te_n T$ is the tableau obtained from $T$ by removing  $\boxed{\ov{n}}$ if $T$ has $\boxed{\ov{n}}$ on its top, and ${\bf 0}$, otherwise,

\item[(3)] if $T\in {\bf T}^{\mf b_n}(a)$, then $\te_n T$ is the tableau given by $\te_n(T^{\tt R}\otimes T^{\tt L})$, where we regard ${\bf T}^{\mf b_n}(a)\subset ({\bf T}^{\rm sp})^{\otimes 2}$.
\end{itemize}
Then $\B$ is a $\g$-crystal with respect to $\te_i$, $\tf_i$ for $i\in I$ and
\begin{equation*}
{\rm wt}(T)=
\begin{cases}
2\omega_{n} + {\rm wt}_{\mf l}(T), & \text{if $T\in {\bf T}^{\mf b_n}(a)$},\\
\omega_{n} + {\rm wt}_{\mf l}(T), & \text{otherwise},
\end{cases}
\end{equation*} 
for $T\in \B$.
By \cite[Theorem 7.1]{K15} we have 
\begin{equation}
\begin{split}
&{\bf T}^{\mf c_n}(a)\cong B(\omega_{n-a})\quad (0\leq a<n),\\
&{\bf T}^{\mf b_n}(a)\cong B(\omega_{n-a})\quad (1\leq a<n),\quad
{\bf T}^{\mf b_n}(0)\cong B(2\omega_{n}),\quad  {\bf T}^{\rm sp}\cong B(\omega_n).
\end{split}
\end{equation}
Note that the highest weight element in ${\bf T}^{\g}(a)$ is given by $H_a\in SST_{[\ov{n}]}(\la(a,0,0))$ where $H^{\tt L}_a(i)=\ov{n-a+i}$ for $1\leq i\leq a$ and $H^{\tt R}_a$ is empty, while the highest weight element $H_{\rm sp}$ in ${\bf T}^{\rm sp}$ is an empty tableau.  

Suppose that $T\in {\bf T}^{\mf g}(a)$ ($0\leq a<n$) is given. We define 
\begin{equation}\label{eq:LT RT}
{}^{\tt L}T=({\mc E}^a T)^{\tt L}, \quad   {}^{\tt R}T= ({\mc E}^a T)^{\tt R}.
\end{equation}

\begin{ex}\label{ex:R matrix}
{\rm
$$
{\bf T}^{\mf b_9}(3)\ni\ T=  \
\resizebox{.48\hsize}{!}{$
{\def\lr#1{\multicolumn{1}{|@{\hspace{.75ex}}c@{\hspace{.75ex}}|}{\raisebox{-.04ex}{$#1$}}}\raisebox{-.6ex}
{$\begin{array}{cc}
\cline{2-2}
 \cline{2-2}
 &\lr{\ov{9}}\\
\cline{1-1}\cline{2-2}
\lr{\ov{8}}&\lr{\ov{7}}\\
\cline{1-1}\cline{2-2}
\lr{\ov{6}}&\lr{\ov{2}}\\
\cline{1-1}\cline{2-2}
\lr{\ov{4}}&\cdot \\
\cline{1-1}
\lr{\ov{2}} \\
\cline{1-1}
\lr{\ov{1}} \\
\cline{1-1} \\
\!\!\!\! T^{\tt L}\!\! & \!\! T^{\tt R}\!\!\!\!\!\!
\end{array}$}}
\ \ \  \stackrel{\mc E}{\longrightarrow}  \ \ \
{\def\lr#1{\multicolumn{1}{|@{\hspace{.75ex}}c@{\hspace{.75ex}}|}{\raisebox{-.04ex}{$#1$}}}\raisebox{-.6ex}
{$\begin{array}{cc}
\cline{2-2}
&\lr{\ov{9}}\\
\cline{2-2}
& \lr{\ov{7}}\\
\cline{1-1}\cline{2-2}
\lr{\ov{8}}&\lr{\ov{6}}\\
\cline{1-2} 
\lr{\ov{4}}&\lr{\ov{2}}  \\
\cline{1-1}\cline{2-2}
\lr{\ov{2}}&\cdot  \\
\cline{1-1} 
\lr{\ov{1}}& \\
\cline{1-1} \\ \\
\end{array}$}}
\ \ \  \stackrel{\mc E}{\longrightarrow}  \ \ \
{\def\lr#1{\multicolumn{1}{|@{\hspace{.75ex}}c@{\hspace{.75ex}}|}{\raisebox{-.04ex}{$#1$}}}\raisebox{-.6ex}
{$\begin{array}{cc}
\cline{2-2}
&\lr{\ov{9}}\\
\cline{2-2}
& \lr{\ov{7}}\\
\cline{2-2}
&\lr{\ov{6}}\\
\cline{1-2} 
\lr{\ov{8}}&\lr{\ov{4}}  \\
\cline{1-1}\cline{2-2}
\lr{\ov{2}}&\lr{\ov{2}}  \\
\cline{1-2} 
\lr{\ov{1}}&\cdot \\
\cline{1-1} \\ \\
\end{array}$}}
\ \ \  \stackrel{\mc E}{\longrightarrow}  \ \ \
{\def\lr#1{\multicolumn{1}{|@{\hspace{.75ex}}c@{\hspace{.75ex}}|}{\raisebox{-.04ex}{$#1$}}}\raisebox{-.6ex}
{$\begin{array}{cc}
\cline{2-2}
&\lr{\ov{9}}\\
\cline{2-2}
& \lr{\ov{7}}\\
\cline{2-2}
& \lr{\ov{6}}\\
\cline{2-2}
&\lr{\ov{4}} \\
\cline{1-2}
\lr{\ov{8}}&\lr{\ov{2}} \\
\cline{1-2}
\lr{\ov{2}}&\lr{\ov{1}} \\
\cline{1-2}\\ 
\!\!\!\! {}^{\tt L}T\!\! & \!\! {}^{\tt R}T\!\!\!\!\!\!
\end{array}$}}$}
\ = {\mc E}^3 T$$
}
\end{ex}

\begin{rem}\label{rem:equivalent algorithm}
{\rm
One can also determine 
$({}^{\tt L}T, {}^{\tt R}T)$ by the following algorithm \cite[(3.6)]{K16}, which will be useful in later arguments:
\begin{itemize}
\item[($\mf s$-1)] Let $y_i=T^{\tt R}(i)$ for $1\leq i\leq {\rm ht}(T^{\tt R})$. First, slide down the box $\boxed{y_1}$ in $T^{\tt R}$ as far as the entry of $T^{\tt L}$ in the same row is no greater than $y_1$. If  no entry of $T^{\tt L}$ is greater than $y$, we place $\boxed{y_1}$ next to $T^{\tt L}(1)$.

\item[($\mf s$-2)] Next, slide down $\boxed{y_2}$ until it is above $\boxed{y_1}$ and the entry of $T^{\tt L}$ in the same row is no greater than $y_2$. Repeat the same process with the other boxes $\boxed{y_3}, \boxed{y_4},\ldots$  until there is no moving down.

\item[($\mf s$-3)] Slide each box $\boxed{x}$ in $T^{\tt L}$ to the right if its right position is empty  (indeed the number of such boxes is $a$).

\item[($\mf s$-4)] Define ${}^{\tt R}T$ to be the tableau determined by the boxes $\boxed{y_i}$\,'s in $T^{\tt R}$ together with boxes $\boxed{x}$'s which have moved from $T^{\tt L}$ by ($\mf s$-3), and define ${}^{\tt L}T$ to be the tableau with the remaining boxes on the left.
\end{itemize}
For example, let $T$ be as in Example \ref{ex:R matrix}.
 
$$
\resizebox{.45\hsize}{!}{$
{\def\lr#1{\multicolumn{1}{|@{\hspace{.75ex}}c@{\hspace{.75ex}}|}{\raisebox{-.04ex}{$#1$}}}\raisebox{-.6ex}
{$\begin{array}{cc}
\cline{2-2}
 \cline{2-2}
 &\lr{\ov{9}}\\
\cline{1-1}\cline{2-2}
\lr{\ov{8}}&\lr{\ov{7}}\\
\cline{1-1}\cline{2-2}
\lr{\ov{6}}&\lr{\ov{2}}\\
\cline{1-1}\cline{2-2}
\lr{\ov{4}} \\
\cline{1-1}
\lr{\ov{2}} \\
\cline{1-1}
\lr{\ov{1}} \\
\cline{1-1}\\ 
\!\!\!\! T^{\tt L}\!\! & \!\! T^{\tt R}\!\!\!\!\!\!
\end{array}$}}
\ \ \  \rightarrow  \ \ \
{\def\lr#1{\multicolumn{1}{|@{\hspace{.75ex}}c@{\hspace{.75ex}}|}{\raisebox{-.04ex}{$#1$}}}\raisebox{-.6ex}
{$\begin{array}{cc}
\cline{2-2}
&\lr{\ov{9}}\\
\cline{1-1}\cline{2-2}
\lr{\ov{8}}& \lr{\ov{7}}\\
\cline{1-1}\cline{2-2}
\lr{\ov{6}}\\
\cline{1-1} 
\lr{\ov{4}}&  \\
\cline{1-1}\cline{2-2}
\lr{\ov{2}}&\lr{\ov{2}} \\
\cline{1-1}\cline{2-2}
\lr{\ov{1}}& \\
\cline{1-1} \\ \\
\end{array}$}}
\ \ \  \rightarrow  \ \ \
{\def\lr#1{\multicolumn{1}{|@{\hspace{.75ex}}c@{\hspace{.75ex}}|}{\raisebox{-.04ex}{$#1$}}}\raisebox{-.6ex}
{$\begin{array}{cc}
\cline{2-2}
&\lr{\ov{9}}\\
\cline{1-1}\cline{2-2}
\lr{\ov{8}}& \lr{\ov{7}}\\
\cline{1-1}\cline{2-2}
& \lr{\ov{6}}\\
\cline{2-2}
&\lr{\ov{4}} \\
\cline{1-1}\cline{2-2}
\lr{\ov{2}}&\lr{\ov{2}} \\
\cline{1-1}\cline{2-2}
&\lr{\ov{1}} \\
\cline{2-2} \\ \\
\end{array}$}}
\ \ \  \rightarrow  \ \ \
{\def\lr#1{\multicolumn{1}{|@{\hspace{.75ex}}c@{\hspace{.75ex}}|}{\raisebox{-.04ex}{$#1$}}}\raisebox{-.6ex}
{$\begin{array}{cc}
\cline{2-2}
&\lr{\ov{9}}\\
\cline{1-1}\cline{2-2}
\lr{\ov{8}}& \lr{\ov{7}}\\
\cline{1-1}\cline{2-2}
\lr{\ov{2}}& \lr{\ov{6}}\\
\cline{1-1}\cline{2-2}
&\lr{\ov{4}} \\
\cline{2-2}
&\lr{\ov{2}} \\
\cline{2-2}
&\lr{\ov{1}} \\
\cline{2-2}\\
\!\!\!\! {}^{\tt L}T\!\! & \!\! {}^{\tt R}T\!\!\!\!\!\!
\end{array}$}}$}
$$\vskip 2mm
\noindent Here the bottoms of the columns in each pair of $(T^{\tt L},{}^{\tt R}T)$ and $({}^{\tt L}T, T^{\tt R})$ are placed on the same horizontal lines, respectively.}
\end{rem}

\vskip 2mm

\begin{df}\label{def:psst}{\rm \mbox{}
\begin{itemize}
\item[(1)]
For $0\leq a_1\leq a_2<n$, we say that a pair  $(T_2,T_1)\in {\bf T}^{\mf g}(a_2)\times {\bf T}^{\mf g}(a_1)$ or ${\bf T}^{\mf g}(a_2)\times  {\bf T}^{\rm sp}$ is {\em admissible}, and write  $T_2\prec T_1$ if it satisfies

\begin{enumerate}
\item[(i)]  ${\rm ht}({T}^{\tt R}_2)\leq {\rm ht}(T^{\tt L}_1) -a_1$,

\item[(ii)] ${T}^{\tt R}_2(i)\leq  {}^{\tt L}T_1(i)$ for $i\geq 1$,

\item[(iii)] ${}^{\tt R}T_2(i+a_2-a_1)\leq {T}^{\tt L}_1(i)$ for $i\geq 1$,
\end{enumerate}
where we assume that $a_1=0$, $T_1=T_1^{\tt L}={}^{\tt L}T_1$ if $T_1\in {\bf T}^{\rm sp}$.


\item[(2)] For $\la\in \cP_n\cup\cP_n^{\rm sp}$, 
we define
\begin{equation*}
{\bf T}^{\mf g}_{\la}=
\begin{cases}
\{\,{\bf T}=(T_\ell\,\ldots,T_1)\in \widehat{\bf T}^{\mf g}_{\la}\ |\ \text{$T_{\ell}\prec \cdots \prec T_{1}$}\,\},& \text{if $\la\in\cP_n$},\\
\{\,{\bf T}=(T_\ell\,\ldots, T_1,T_0)\in \widehat{\bf T}^{\mf b_n}_{\la}\ |\ \text{$T_{\ell}\prec \cdots \prec T_{0}$}\,\},& \text{if $\la\in\cP_n^{\rm sp}$},\\
\end{cases}
\end{equation*}
where 
\begin{equation*}
\widehat{\bf T}^{\mf g}_{\la}=
\begin{cases}
{\bf T}^{\mf g}(a_\ell)\times \cdots \times {\bf T}^{\mf g}(a_1), & \text{if $\la\in\cP_n$}, \\
{\bf T}^{\mf b_n}(a_\ell)\times \cdots \times {\bf T}^{\mf b_n}(a_1)\times {\bf T}^{\rm sp}, & \text{if $\la\in\cP_n^{\rm sp}$}, \\
\end{cases}
\end{equation*}
and $a_\ell\geq \ldots \geq a_1$ is the sequence satisfying 
\begin{equation*}
\omega_\la=
\begin{cases}
\sum_{i=1}^\ell\omega_{n-a_i}, & \text{if $\la\in\cP_n$}, \\
\sum_{i=1}^\ell\omega_{n-a_i}+\omega_n, & \text{if $\la\in\cP_n^{\rm sp}$}. \\
\end{cases}
\end{equation*}
\end{itemize}
}
\end{df}

\begin{rem}\label{rem:vertical position}
{\rm
When we consider the admissibility of a pair $(T_2,T_1)$, it is convenient to assume that $T^{\tt body}_i$ and $T^{\tt tail}_i$ for $i=1,2$ are separated by a common horizontal line so that the vertical positions of the other column tableaux ${}^{\tt L}T_i$, ${}^{\tt R}T_i$ for $i=1,2$ are determined accordingly (see the dashed lines in Example \ref{ex:admissibility}).
Then the condition (1)(ii) and (1)(iii) are equivalent to saying that $(T^{\tt R}_2,{}^{\tt L}T_1)$ and $({}^{\tt R}T_2,T^{\tt L}_1)$ form semistandard tableaux with respect to this vertical position, respectively.
}
\end{rem}

\begin{ex}\label{ex:admissibility}
{\rm
For $T_2\in {\bf T}^{\mf c_7}(2)$ and $T_1\in {\bf T}^{\mf c_7}(1)$ below, we have

$$
\begin{array}{cccc}
  T_2 \quad\quad\quad &   {\mc E}^2T_2 \quad\quad\quad\quad &  T_1 \quad\quad\quad &  {\mc E}T_1\\
\end{array}
$$
$$\resizebox{.5\hsize}{!}
{\def\lr#1{\multicolumn{1}{|@{\hspace{.75ex}}c@{\hspace{.75ex}}|}{\raisebox{-.04ex}{$#1$}}}\raisebox{-.6ex}
{$\begin{array}{cccccccccccccccccc}
\cline{12-12}\cline{13-13}\cline{16-16}\cline{17-17}
& & & & & & & & & & & \lr{\color{blue}\ov{7}} &\lr{\ov{5}} & & & \lr{\color{red}\ov{6}} & \lr{\ov{7}} \\
\cline{2-2}\cline{3-3}\cline{6-6}\cline{7-7}\cline{12-12}\cline{13-13}\cline{16-16}\cline{17-17}
& \lr{\ov{6}}&\lr{\color{red}\ov{4}} & & & \lr{\ov{4}} &  \lr{\color{blue}\ov{6}} & & & & & \lr{\color{blue}\ov{6}}&\lr{\ov{2}} & & & \lr{\color{red}\ov{3}} &  \lr{\ov{5}}\\
\cline{2-2}\cline{3-3}\cline{6-6}\cline{7-7}\cline{12-12}\cline{13-13}\cline{16-16}\cline{17-17}
& \lr{\ov{5}} & \lr{\color{red}\ov{2}}& & & \lr{\ov{3}} & \lr{\color{blue}\ov{5}} & & & & & \lr{\color{blue}\ov{3}} & \lr{\ov{1}}& & & \lr{\color{red}\ov{2}} &  \lr{\ov{2}}\\
\cdashline{1-1}[0.5pt/1pt]\cline{2-2}\cline{3-3}\cdashline{4-5}[0.5pt/1pt]\cline{6-6}\cline{7-7}\cdashline{8-11}[0.5pt/1pt]
\cline{12-12}\cline{13-13}\cdashline{14-18}[0.5pt/1pt]\cline{16-16}\cline{17-17}
& \lr{\ov{4}} & & & & & \lr{\color{blue}\ov{4}} & & & & & \lr{\color{blue}\ov{2}} & & & & & \lr{\ov{1}}& \\
\cline{2-2}\cline{7-7}\cline{12-12}\cdashline{1-18}[0.5pt/1pt]\cline{17-17}
& \lr{\ov{3}} & & & & & \lr{\color{blue}\ov{2}}& \\
\cline{2-2}\cdashline{1-18}[0.5pt/1pt]\cline{7-7}
\end{array}$}} 
$$ \vskip 3mm
\noindent Then 
$({}^{\tt R}T_2,T^{\tt L}_1)$ (in blue) and $(T^{\tt R}_2,{}^{\tt L}T_1)$ (in red) form semistandard tableaux \vskip 2mm
$$ 
\resizebox{.06\hsize}{!}
{\def\lr#1{\multicolumn{1}{|@{\hspace{.75ex}}c@{\hspace{.75ex}}|}{\raisebox{-.04ex}{$#1$}}}\raisebox{-.6ex}
{$\begin{array}{cc}
 \cline{2-2}
  & \lr{\color{blue}\ov{7}} \\
\cline{1-1}\cline{2-2}
  \lr{\color{blue}\ov{6}} & \lr{\color{blue}\ov{6}}\\
\cline{1-1}\cline{2-2}
 \lr{\color{blue}\ov{5}}& \lr{\color{blue}\ov{3}}  \\
\cline{1-1}\cline{2-2}
\lr{\color{blue}\ov{4}} & \lr{\color{blue}\ov{2}} \\
\cline{1-1}\cline{2-2}
\lr{\color{blue}\ov{2}} \\ 
\cline{1-1}\\
\!\!\!\!\!\! {}^{\tt R}T_2\!\! & \!\! T^{\tt L}_1\!\!\!\!\!\!
\end{array}$}}
 \ \ \ \ \  \ \ \ \ \ 
\resizebox{.06\hsize}{!}
{\def\lr#1{\multicolumn{1}{|@{\hspace{.75ex}}c@{\hspace{.75ex}}|}{\raisebox{-.04ex}{$#1$}}}\raisebox{-.6ex}
{$\begin{array}{cc}
 \cline{2-2}
  & \lr{\color{red}\ov{6}} \\
\cline{1-1}\cline{2-2}
\lr{\color{red}{\ov{4}}} &  \lr{\color{red}\ov{3}}\\
\cline{1-1}\cline{2-2}
\lr{\color{red}\ov{2}} & \lr{\color{red}\ov{2}}  \\
\cline{1-1}\cline{2-2}\\
\\ \\
\!\!\!\! T^{\tt R}_2\!\! & \!\! {}^{\tt L}T_1\!\!\!\!\!\!
\end{array}$}}
$$\vskip 2mm
\noindent
which implies that $T_2\prec T_1$, and hence $(T_2,T_1)\in {\bf T}^{\mf c_7}_{\la}$, where $\la=(2^5,1)$.
}
\end{ex}

One may define a $\g$-crystal structure on $\widehat{\bf T}^{\mf g}_{\la}$ by identifying $(\,\ldots,T_2,T_1)$ with $T_1 \otimes T_2\otimes  \cdots $. Then ${\bf T}^{\mf g}_{\la}\cup\{{\bf 0}\}$ is invariant under $\te_i, \tf_i$ for $i\in I$, and ${\bf T}^{\mf g}_{\la}$ becomes a $\g$-crystal.
The following is proved in \cite[Theorem 7.4]{K15}.

\begin{thm}\label{thm:PSST iso B(la)} For $\la\in \cP_n\cup\cP_n^{\rm sp}$, we have
\begin{equation*}
{\bf T}^{\mf g}_{\la}\cong B(\omega_\la).
\end{equation*}
\end{thm}

Note that the highest weight element in ${\bf T}^{\mf g}_{\la}$ is given by
\begin{equation}\label{eq:hwv in PSST}
\begin{cases}
(H_{a_\ell}, \cdots , H_{a_1}), & \text{if $\la\in\cP_n$}, \\
(H_{a_\ell}, \cdots , H_{a_1},H_{\rm sp}), & \text{if $\la\in\cP_n^{\rm sp}$},
\end{cases}
\end{equation}
where $(a_1,\ldots,a_\ell)$ is the sequence given in Definition \ref{def:psst}(2), and $H_a$ (resp. $H_{\rm sp}$) is the highest weight element in ${\bf T}^\g(a)$ (resp. ${\bf T}^{\rm sp}$).


\subsection{Isomorphism between KN tableaux and spinor model}
In this subsection, we construct an explicit isomorphism from ${\bf KN}^{\g}_{\la}$ to ${\bf T}^{\g}_{\la}$ for $\la\in \cP_n \cup \cP_n^{\rm sp}$, which was briefly mentioned in \cite[Section 7.2]{K15} without proof. 
 
\subsubsection{Type $C_n$} 
For $0\leq a< n$ and $T\in {\bf T}^{\mf c_n}(a)$, let $\td{T}$ be the tableau in $SST_{\J_n^\times}(1^{n-a})$ defined as follows:
\begin{itemize}
\item[(1)] Let $\td{{}^{\tt R}T}$ be the tableau in $SST_{[n]}(1^m)$ with $m=n-{\rm ht}({}^{\tt R}T)$ such that $k$ appears in $\td{{}^{\tt R}T}$ if and only if $\ov{k}$ does not appear in ${}^{\tt R}T$ for each $k\in [n]$.

\item[(2)] Define $\td{T}$ to be the tableau in $SST_{\J_n^\times}(1^{n-a})$ given by putting ${}^{\tt L}T$ below $\td{{}^{\tt R}T}$.
\end{itemize}

\begin{ex}{\rm
$$
{\bf T}^{\mf c_5}(1) \ni T=\ 
\resizebox{.06\hsize}{!}{$
{\def\lr#1{\multicolumn{1}{|@{\hspace{.75ex}}c@{\hspace{.75ex}}|}{\raisebox{-.04ex}{$#1$}}}\raisebox{-.6ex}
{$\begin{array}{cc}
\\
\cline{1-1}\cline{2-2}
\lr{\overline{5}}&\lr{\overline{4}}\\ 
\cline{1-1}\cline{2-2}
\lr{\overline{3}}&\lr{\overline{2}}\\ 
\cline{1-1}\cline{2-2}
\lr{\overline{2}}\\ 
\cline{1-1} \\
\!\!\! T^{\tt L}\!\!\! & \!\! T^{\tt R}\!\!\!\!\!\!
\end{array}$}}$}\  \ \ \ \ \ \ \ \  
\resizebox{.06\hsize}{!}{$
{\def\lr#1{\multicolumn{1}{|@{\hspace{.75ex}}c@{\hspace{.75ex}}|}{\raisebox{-.04ex}{$#1$}}}\raisebox{-.6ex}
{$\begin{array}{cc}
\\
\cline{1-1}\cline{2-2}
\lr{\overline{5}}&\lr{\overline{4}}\\ 
\cline{1-1}\cline{2-2}
\lr{\overline{2}}&\lr{\overline{3}}\\ 
\cline{1-1}\cline{2-2}
&\lr{\overline{2}}\\ 
\cline{2-2}  \\
\!\!\! {}^{\tt L}T\!\!\! & \!\! {}^{\tt R}T\!\!\!\!\!\!

\end{array}$}}$}\ \ \ \ \ \  \td{{}^{\tt R}T} \ = \
\resizebox{.03\hsize}{!}{$
{\def\lr#1{\multicolumn{1}{|@{\hspace{.75ex}}c@{\hspace{.75ex}}|}{\raisebox{-.04ex}{$#1$}}}\raisebox{-.6ex}
{$\begin{array}{c}
\cline{1-1} 
\lr{1} \\ 
\cline{1-1} 
\lr{5} \\ 
\cline{1-1}   
\end{array}$}}$}\ \ \ \ \ \  \td{T} \ = \
\resizebox{.03\hsize}{!}{$
{\def\lr#1{\multicolumn{1}{|@{\hspace{.75ex}}c@{\hspace{.75ex}}|}{\raisebox{-.04ex}{$#1$}}}\raisebox{-.6ex}
{$\begin{array}{c}
\cline{1-1} 
\lr{1} \\ 
\cline{1-1} 
\lr{5} \\ 
\cline{1-1}   
\lr{\overline{5}}\\ 
\cline{1-1}
\lr{\overline{2}}\\ 
\cline{1-1}  
\end{array}$}}$}
$$}
\end{ex}

\begin{lem}\label{lem:psst_to_KN_C} The map $T \longmapsto \td{T}$ is an isomorphism of $\mf c_n$-crystals from ${\bf T}^{\mf c_n}(a)$ to ${\bf KN}^{\mf c_n}_{(1^{n-a})}$ for $0\leq a< n$.
\end{lem}
\pf By definition, it is clear that the map $T\longmapsto \td{T}$ is injective and preserves the weight. We first claim that $\td{T}\in {\bf KN}^{\mf c_n}_{(1^{n-a})}$. Note that $\td{T}\in SST_{\J_n^\times}(1^{n-a})$. Suppose that $\td{T}\not\in {\bf KN}^{\mf c_n}_{(1^{n-a})}$. By ($\mf c$-1) in Definition \ref{def:KN-C}, there exists $c\in [n]$ such $\td{T}(i_1)=\ov{c}$ and $\td{T}(i_2)=c$ for some $i_1<i_2$, and $i_1+(n-a-i_2+1)>c$. Put $x=n-a-i_2+1$ and $y=i_1$. Then $x$ is the number of entries in $\td{T}$ smaller than $c$, and $y$ is the number of entries in $\td{T}$ or ${}^{\tt L}T$ greater than or equal to $\ov{c}$. Note that there is no $\ov{c}$ in ${}^{\tt R}T$, and the number of entries in ${}^{\tt R}T$ greater than $\ov{c}$ is $c-x$. On the other hand, the pair $({}^{\tt L}T,{}^{\tt R}T)$ is obtained from $T$ by applying the steps ($\mf s$-1)-($\mf s$-4) in Section \ref{subsec:PSST}, and it is reversible. But if $y>c-x$, then the inverse steps from $({}^{\tt L}T,{}^{\tt R}T)$ to $T$ is not well-defined, which is a contradiction. This proves our claim. In particular, it implies that the map is a bijection from ${\bf T}^{\mf c_n}(a)$ to ${\bf KN}^{\mf c_n}_{(1^{n-a})}$ since both of them have the same cardinality as $B(\omega_{n-a})$. 

Next we claim that the map $T\longmapsto \td{T}$ commutes with $\te_i$ and $\tf_i$ for $i\in I$. Since ${\mc E}^a T$ is Knuth equivalent to $T$, we have $T\equiv_{\mf l} {}^{\tt R}T\otimes {}^{\tt L}T$, where $\equiv_{\mf l}$ means $\equiv$ as elements in $\mf l$-crystals. Also we have ${}^{\tt R}T\equiv_{\mf l} \td{{}^{\tt R}T}$ and $\td{T}\equiv_{\mf l} \td{{}^{\tt R}T}\otimes {}^{\tt L}T$. Hence $T\equiv_{\mf l} \td{T}$ and the map $T\longmapsto \td{T}$ commutes with $\te_i$ and $\tf_i$ for $i\in I\setminus\{n\}$. Now, consider $\te_n$ and $\tf_n$. It is not difficult to see that
\begin{equation*}
\begin{split}
\te_n T\neq {\bf 0}& \Longleftrightarrow \text{
$T$ has $\resizebox{.05\hsize}{!}{\def\lr#1{\multicolumn{1}{|@{\hspace{.6ex}}c@{\hspace{.6ex}}|}{\raisebox{-.25ex}{$#1$}}}\raisebox{-.65ex}
{$\begin{array}[b]{cc}
\cline{1-1}\cline{2-2}
\lr{\ov{n}}&\lr{\ov{n}}\\
\cline{1-1}\cline{2-2}
\end{array}$}}$ on its top
}\\
& \Longleftrightarrow \text{${}^{\tt L}T$ has $\ov{n}$ while $\td{{}^{\tt R}T}$ does not have $n$}\\
& \Longleftrightarrow \te_n \td{T}\neq {\bf 0},
\end{split}
\end{equation*}
and $\td{T'}=\te_n\td{T}$ for $T'=\td{e}_n T$ in this case. The proof for $\tf_n$ is similar. 
Therefore the map is an isomorphism of crystals.
\qed

For $0\leq a<n$ and $T\in {\bf KN}^{\mf c_n}_{(1^a)}$, we define
\begin{equation}\label{eq:KN_C_to_psst}
\Psi_a(T)={\mc F}^{n-a}(T_-,\td{T_+}),
\end{equation}
where $T_+$ (resp. $T_-$) is the subtableau of $T$ with entries in $[n]$ (resp. $[\ov{n}]$), and $\td{T_+}\in SST_{[\ov{n}]}(1^m)$ with $m=n-{\rm ht}(T_+)$ is such that $\ov k$ appears in $\td{T_+}$ if and only if $k$ does not appear in $T_+$ for $k\in [n]$. Here we regard $\Psi_a(T)$ as a two-columned skew tableau with $\mf r_{\Psi_a(T)}=0$.

\begin{lem}\label{lem:KN_C_to_psst}
The map $\Psi_a$ is an isomorphism of ${\mf c}_n$-crystals from ${\bf KN}^{\mf c_n}_{(1^a)}$ to ${\bf T}^{\mf c_n}(n-a)$ for $0< a \leq n$.
\end{lem}
\pf By \eqref{eq:LT RT}, we see immediately that $\Psi_a$ is the inverse of the map in Lemma \ref{lem:psst_to_KN_C}.
\qed
  
\begin{ex}{\rm  
$$
{\bf KN}_{(1^4)}^{\mf c_5}\ni {T} = \
\resizebox{.03\hsize}{!}{$
{\def\lr#1{\multicolumn{1}{|@{\hspace{.75ex}}c@{\hspace{.75ex}}|}{\raisebox{-.04ex}{$#1$}}}\raisebox{-.6ex}
{$\begin{array}{c}
\cline{1-1} 
\lr{1} \\ 
\cline{1-1} 
\lr{5} \\ 
\cline{1-1}   
\lr{\overline{5}}\\ 
\cline{1-1}
\lr{\overline{2}}\\ 
\cline{1-1}  
\end{array}$}}$}\ \ \ \ \ \ \ \ \ \ 
(T_-, \td{T_+}) =\
\resizebox{.06\hsize}{!}{$  
{\def\lr#1{\multicolumn{1}{|@{\hspace{.75ex}}c@{\hspace{.75ex}}|}{\raisebox{-.04ex}{$#1$}}}\raisebox{-.6ex}
{$\begin{array}{cc}
\cline{2-2}
\cdot &\lr{\overline{4}}\\ 
\cline{1-1}\cline{2-2}
\lr{\overline{5}}&\lr{\overline{3}}\\ 
\cline{1-1}\cline{2-2}
\lr{\overline{2}}&\lr{\overline{2}}\\ 
\cline{1-2}  
\end{array}$}}$}
\ \ \  \stackrel{{\mc F}}{\longrightarrow } \ \ \
\resizebox{.06\hsize}{!}{$
{\def\lr#1{\multicolumn{1}{|@{\hspace{.75ex}}c@{\hspace{.75ex}}|}{\raisebox{-.04ex}{$#1$}}}\raisebox{-.6ex}
{$\begin{array}{cc}
\cline{1-1}\cline{2-2}
\lr{\overline{5}}&\lr{\overline{4}}\\ 
\cline{1-1}\cline{2-2}
\lr{\overline{3}}&\lr{\overline{2}}\\ 
\cline{1-1}\cline{2-2}
\lr{\overline{2}}\\ 
\cline{1-1}  
\end{array}$}}
$}
\ = \Psi_4(T)  \in {\bf T}^{\mf c_5}(1)$$
}
\end{ex}

\begin{thm}\label{thm:KN_C_to_psst} 
Let $\la\in \cP_n$ with $\ell=\la_1$. The map
\begin{equation*}\label{eq:Psi_lambda:C}
\xymatrixcolsep{1pc}\xymatrixrowsep{0pc}\xymatrix{
\Psi_{\la} : {\bf KN}^{\mf c_n}_{\la}  \ar@{->}[r]  & \ {\bf T}^{\mf c_n}_{\la}  \\
\ \ \ T \   \ar@{|->}[r] &  (\Psi_{\la'_\ell}(T_l),\ldots,\Psi_{\la'_1}(T_1)) }
\end{equation*}
is an isomorphism of $\mf c_n$-crystals,
where $T_i$ denotes the $i$th column of $T$ from the right.
\end{thm}
\pf Note that $\Psi_\la(T)\in \widehat{\bf T}^{\mf c_n}_{\la}$ for $T\in {\bf KN}^{\mf c_n}_{\la}$. It is clear from the crystal structures on ${\bf KN}^{\mf c_n}_{\la}$ and $\widehat{\bf T}^{\mf c_n}_{\la}$, and Lemma \ref{lem:KN_C_to_psst} that $\Psi_\la : {\bf KN}^{\mf c_n}_{\la}\longrightarrow {\rm Im}\Psi_\la\subset \widehat{\bf T}^{\mf c_n}_{\la}$ is an isomorphism of $\mf c_n$-crystals.
Recall that if $H$ is the highest weight element in ${\bf KN}^{\mf c_n}_{\la}$, then the $i$th entry from the top in each column of $H$ is filled with $i$ for $1\leq i\leq n$. One can check that $H$ is mapped to the highest weight element in 
${\bf T}^{\mf c_n}_{\la}$ \eqref{eq:hwv in PSST} under $\Psi_\la$.
So ${\rm Im}\Psi_\la$ is the connected component of the highest weight element of ${\bf T}^{\mf c_n}_{\la}$. This implies that ${\rm Im}\Psi_\la={\bf T}^{\mf c_n}_{\la}$.
\qed

\begin{ex}\label{ex:embedding-1}
{\rm
Let  \vskip 2mm
$$
T=\quad
\resizebox{.12\hsize}{!}
{\def\lr#1{\multicolumn{1}{|@{\hspace{.75ex}}c@{\hspace{.75ex}}|}{\raisebox{-.04ex}{$#1$}}}\raisebox{-.6ex}
{$\begin{array}{cccc}
\cline{3-3}\cline{4-4}
& & \lr{1} &\lr{4}\\
\cline{2-2}\cline{3-3}\cline{4-4}
& \lr{3} & \lr{5} &\lr{\ov 4}\\
\cline{1-1}\cline{2-2}\cline{3-3}\cline{4-4}
\lr{5}&\lr{\ov 5} & \lr{\ov 5}&\lr{\ov 2}\\
\cline{1-1}\cline{2-2}\cline{3-3}\cline{4-4}
\lr{\ov 3}&\lr{\ov 2} & \lr{\ov 2}&\lr{\ov 1}\\
\cline{1-1}\cline{2-2}\cline{3-3}\cline{4-4}\\
\end{array}$}}\quad  \in {\bf KN}^{\mf c_5}_{\la},
$$ 
\noindent where $\la=(4,4,3,2,0)\in \cP_5$. Then $\Psi_{\la}(T)={\bf T}=(\Psi_2(T_4),\Psi_3(T_3),\Psi_4(T_2),\Psi_4(T_1))\in {\bf T}^{\mf c_5}_{\la}$ by \eqref{eq:KN_C_to_psst}, where 

$$\!\!\!\!
\resizebox{.43\hsize}{!}
{\def\lr#1{\multicolumn{1}{|@{\hspace{.75ex}}c@{\hspace{.75ex}}|}{\raisebox{-.04ex}{$#1$}}}\raisebox{-.6ex}
{$\begin{array}{ccccccccccccccccc}
\cline{10-10}\cline{14-14}
& & & & & & & & & \lr{1} & & & &\lr{4}\\
\cline{6-6}\cline{10-10}\cline{14-14}& & & & & \lr{3} & & & & \lr{5} & & & &\lr{\ov 4}\\
\cline{2-2}\cline{6-6}\cline{10-10}\cline{14-14}
& \lr{5} & & & & \lr{\ov 5} & & & & \lr{\ov 5} & & & & \lr{\ov 2}\\
\cline{2-2}\cline{6-6}\cline{10-10}\cline{14-14}
\ \, &\lr{\ov 3} & &  \, \ & & \lr{\ov 2} & &   \, \ & & \lr{\ov 2} & &  \, \ & & \lr{\ov 1}& \\
\cline{2-2}\cline{6-6}\cline{10-10}\cline{14-14}\cdashline{1-15}[0.5pt/1pt]\\ 
& \!\!\!\!\!\!\!\!\!T_4 \!\!\!\!\!\!\!\!\! & & &
& \!\!\!\!\!\!\!\!\!T_3 \!\!\!\!\!\!\!\!\! & & &
& \!\!\!\!\!\!\!\!\!T_2 \!\!\!\!\!\!\!\!\! & & & 
& \!\!\!\!\!\!\!\!\!T_1 \!\!\!\!\!\!\!\!\! & & & 
\end{array}$}}
$$
$$\resizebox{.42\hsize}{!}
{\def\lr#1{\multicolumn{1}{|@{\hspace{.75ex}}c@{\hspace{.75ex}}|}{\raisebox{-.04ex}{$#1$}}}\raisebox{-.6ex}
{$\begin{array}{ccccccccccccccccc}
\\
\cline{14-14}\cline{15-15}
& & & & & & & & & & & & & \lr{\ov{5}} & \lr{\ov{3}} \\
\cline{6-6}\cline{7-7}\cline{10-11}\cline{14-15}
& & & & & \lr{\ov{5}} &  \lr{\ov{5}} & &  & \lr{\ov{5}}&\lr{\ov{4}} & & & \lr{\ov{4}} &  \lr{\ov{2}}\\
\cline{2-3}\cline{6-7}\cline{10-11}\cline{14-15}
& \lr{\ov{4}} & \lr{\ov{3}}& & & \lr{\ov{4}} & \lr{\ov{2}} & & & \lr{\ov{3}} & \lr{\ov{2}}& & & \lr{\ov{2}} &  \lr{\ov{1}}& \\
\cline{2-3}\cdashline{1-5}[0.5pt/1pt]\cline{6-7}\cdashline{8-11}[0.5pt/1pt]
\cline{10-11}\cdashline{12-16}[0.5pt/1pt]\cline{14-15}
& \lr{\ov{3}} & & & & \lr{\ov{2}} & & &  & \lr{\ov{2}} & & & & \lr{\ov{1}} & \\
\cline{2-2}\cline{6-6}\cline{10-10}\cline{14-14}
& \lr{\ov{2}} & & & & \lr{\ov{1}} & \\
\cline{2-2}\cline{6-6}
& \lr{\ov{1}} \\
\cline{2-2}\\
\Psi_2(T_4)\!\!\!\!\!\!\!\!\!\!\!\!\!\!\!\!\!\! & & & & 
\Psi_3(T_3)\!\!\!\!\!\!\!\!\!\!\!\!\!\!\!\!\!\! & & & & 
\Psi_4(T_2)\!\!\!\!\!\!\!\!\!\!\!\!\!\!\!\!\!\! & & & & 
\Psi_4(T_1)\!\!\!\!\!\!\!\!\!\!\!\!\!\!\!\!\!\! 
\end{array}$}} \quad
$$
}
\end{ex}

\subsubsection{Type $B_n$}
As in the case of $C_n$, we first construct an explicit isomorphism from ${\bf KN}^{\mf b_n}_{\la}$ to ${\bf T}^{\mf b_n}_{\la}$ when $\omega_\la$ is a fundamental weight.

For $0\leq a< n$ and $T=(T^{\tt L}, T^{\tt R})\in {\bf T}^{\mf b_n}(a)$, let $\td{T}$ be the tableau in $SST_{\J_n}(1^{n-a})$ defined as follows:
\begin{itemize}
\item[(1)] Let $\td{{}^{\tt R}T}$ be the tableau in $SST_{[n]}(1^m)$ with $m=n-{\rm ht}({}^{\tt R}T)$ defined in the same way as in the case of $C_n$.

\item[(2)] Define $\td{T}$ to be the tableau in $SST_{\J_n}(1^{n-a})$ given by putting below $\td{{}^{\tt R}T}$ a single-columned tableau of height $a+{\rm ht}(T^{\tt R})-{\rm ht}(T^{\tt L})$ with $0$'s and then ${}^{\tt L}T$.
\end{itemize}
If $T\in {\bf T}^{\rm sp}$, then define $\td{T}$ to be the unique tableau in ${\bf KN}^{\rm sp}$ such that 
$\ov{k}$ appears in $T$ if and only if $\ov{k}$ appears in $\td{T}$.

\begin{ex}{\rm

$${\bf T}^{\mf b_5}(1) \ni T= \
\resizebox{.055\hsize}{!}{$
{\def\lr#1{\multicolumn{1}{|@{\hspace{.75ex}}c@{\hspace{.75ex}}|}{\raisebox{-.04ex}{$#1$}}}\raisebox{-.6ex}
{$\begin{array}{cc}
\\
\cline{2-2}
&\lr{\overline{5}}\\ 
\cline{2-2}
&\lr{\overline{4}}\\ 
\cline{1-1}\cline{2-2}
\lr{\overline{3}}&\lr{\overline{1}}\\ 
\cline{1-1}\cline{2-2}
\lr{\overline{1}}\\ 
\cline{1-1} \\
\!\!\! T^{\tt L}\!\!\! & \!\! T^{\tt R}\!\!\!\!\!\! 
\end{array}$}}$}\ \ \ \ \ \ \ \  
\resizebox{.055\hsize}{!}{$
{\def\lr#1{\multicolumn{1}{|@{\hspace{.75ex}}c@{\hspace{.75ex}}|}{\raisebox{-.04ex}{$#1$}}}\raisebox{-.6ex}
{$\begin{array}{cc}
\\
\cline{2-2}
&\lr{\overline{5}}\\ 
\cline{2-2}
&\lr{\overline{4}}\\ 
\cline{1-1}\cline{2-2}
\lr{\overline{1}}&\lr{\overline{3}}\\ 
\cline{1-1}\cline{2-2}
&\lr{\overline{1}}\\ 
\cline{2-2}  \\
\!\!\! {}^{\tt L}T\!\!\! & \!\! {}^{\tt R}T\!\!\!\!\!\! 
\end{array}$}}$}\ \ \ \ \ \ \td{{}^{\tt R}T} =  \
\resizebox{.028\hsize}{!}{$
{\def\lr#1{\multicolumn{1}{|@{\hspace{.75ex}}c@{\hspace{.75ex}}|}{\raisebox{-.04ex}{$#1$}}}\raisebox{-.0ex}
{$\begin{array}{c}
\cline{1-1} 
\lr{2} \\ 
\cline{1-1}   
\end{array}$}}$}\ \ \ \ \ \ \td{T} = \
\resizebox{.026\hsize}{!}{$
{\def\lr#1{\multicolumn{1}{|@{\hspace{.75ex}}c@{\hspace{.75ex}}|}{\raisebox{-.04ex}{$#1$}}}\raisebox{-.6ex}
{$\begin{array}{c}
\cline{1-1} 
\lr{2} \\ 
\cline{1-1} 
\lr{0} \\ 
\cline{1-1}   
\lr{0}\\ 
\cline{1-1}
\lr{\overline{1}}\\ 
\cline{1-1}  
\end{array}$}}$}
$$ \vskip 2mm
$$
{\bf T}^{\rm sp}\ni T=\ 
\resizebox{.028\hsize}{!}
{\def\lr#1{\multicolumn{1}{|@{\hspace{.75ex}}c@{\hspace{.75ex}}|}{\raisebox{-.04ex}{$#1$}}}\raisebox{-.6ex}
{$\begin{array}{c}
\cline{1-1} 
\lr{\ov 5} \\
\cline{1-1} 
\lr{\ov 4} \\
\cline{1-1}
\lr{\ov 2} \\
\cline{1-1} 
\lr{\ov 1} \\
\cline{1-1}
\end{array}$}}
\quad\quad\quad\quad  
\td{T}=\ 
\resizebox{.02\hsize}{!}
{\def\lr#1{\multicolumn{1}{|@{\hspace{.75ex}}c@{\hspace{.75ex}}|}{\raisebox{-.04ex}{$#1$}}}\raisebox{-.6ex}
{$\begin{array}{c}
\cline{1-1} 
\lr{\! 3\!} \\
\cline{1-1} 
\lr{\!\! \ov 5\!\!} \\
\cline{1-1} 
\lr{\!\!\ov 4\!\!} \\
\cline{1-1} 
\lr{\!\! \ov 2\!\!} \\
\cline{1-1}
\lr{\!\! \ov 1\!\!} \\
\cline{1-1}
\end{array}$}}
$$}
\end{ex}
\vskip 2mm

\begin{lem}\label{lem:psst_to_KN_B} The map $T \longmapsto \td{T}$ is an isomorphism of $\mf b_n$-crystals from ${\bf T}^{\mf b_n}(a)$ to ${\bf KN}^{\mf b_n}_{(1^{n-a})}$ for $0\leq a <n$, and from ${\bf T}^{\rm sp}$ to ${\bf KN}^{\rm sp}$, respectively.
\end{lem}
\pf The proof is almost the same as that of Lemma \ref{lem:psst_to_KN_C}.  
\qed\vskip 2mm

For $0< a\leq n$ and $T\in {\bf KN}^{\mf c_n}_{(1^a)}$ and $T'\in {\bf T}^{\rm sp}$, we define
\begin{equation}\label{eq:KN_B_to_psst}
\Psi_a(T)={\mc F}^{n-a}(T_-,\td{T_+}),\quad \Psi_{\rm sp}(T')=(T')_-,
\end{equation}
where $T_\pm$, $\td{T_+}$ are as in \eqref{eq:KN_C_to_psst}.
As in Lemma \ref{lem:KN_C_to_psst}, we see that \eqref{eq:KN_B_to_psst} is the inverse of the isomorphism in Lemma \ref{lem:psst_to_KN_B}. Hence we have

\begin{lem}\label{lem:KN_B_to_psst}
The map $\Psi_a$ and $\Psi_{\rm sp}$ are  isomorphisms of ${\mf b}_n$-crystals from ${\bf KN}^{\mf b_n}_{(1^a)}$ to ${\bf T}^{\mf b_n}(n-a)$ for $0< a \leq n$, and from ${\bf KN}^{\rm sp}$ to ${\bf T}^{\rm sp}$, respectively.
\end{lem}

\begin{thm}\label{thm:KN_B_to_psst} 
Let $\la\in \cP_n\cup \cP_n^{\rm sp}$. 
The map
\begin{equation*}\label{eq:Psi_lambda:B-1}
\xymatrixcolsep{3pc}\xymatrixrowsep{0pc}\xymatrix{
\Psi_{\la} : {\bf KN}^{\mf b_n}_{\la}  \ar@{->}[r]  & \ {\bf T}^{\mf b_n}_{\la}
}
\end{equation*}
given by 
\begin{equation*}
\Psi_\la(T)=
\begin{cases}
(\Psi_{\la'_\ell}(T_\ell),\ldots,\Psi_{\la'_1}(T_1)),&\!\! \text{if $\la\in \cP_n$ with $\ell=\la_1$}, \\
(\Psi_{\mu'_\ell}(T_\ell),\ldots,\Psi_{\mu'_1}(T_1),\Psi_{\rm sp}(T_0)),& \!\!
\text{if $\la=\mu+\sigma_n\in \cP_n^{\rm sp}$ with $\ell=\mu_1$},
\end{cases}
\end{equation*}
is an isomorphism of ${\mf b}_n$-crystals, where $T_i$ denotes the $i$th column of $T$ with width 1 from the right and $T_0$ denotes the column of $T$ with half width.
\end{thm}
\pf It follows from Lemma \ref{lem:KN_B_to_psst} and the same arguments as in Theorem \ref{thm:KN_C_to_psst}.
\qed

\section{Crystal of a maximal parabolic Verma module}\label{sec:parabolic}

In this section, we recall a combinatorial model for the crystal of a maximal parabolic Verma module \cite{K12} and then construct its embedding $\Phi_\la$ into the crystal of the $\bi$-Lusztig data. 

From now on $\epsilon$ denotes a constant such that $\epsilon=1$ for $\g=\mf b_n$ and $\epsilon=2$ for $\g=\mf c_n$.

\subsection{Crystal of a maximal parabolic Verma module}

Let
\begin{equation*}
\cP_{\mf g} =\{\,\epsilon\lambda=(\epsilon\lambda_i)_{i\geq 1}\,|\, \lambda=(\lambda_i)_{i\geq 1}\in\cP \,\}.
\end{equation*}
Put
\begin{equation*}
{\bf V}^{\g} = \bigsqcup_{\tau\in\cP_{\mf g}}SST_{[\ov{n}]}(\tau^\pi),\quad 
{\bf V}^\g_{k}=\bigsqcup_{\substack{\tau\in\cP_{\mf g} \\ \tau_1\leq \epsilon k}}SST_{[\ov{n}]}(\tau^\pi) \quad (k\geq 1).
\end{equation*}

Let $T\in {\bf V}^\g_k$ be given.
When $\g=\mf c_n$, let $T_i$ be the subtableau of $T$ consisting of the $(2i-1)$th and the $2i$-th columns of $T$ from the right for $1\leq i\leq k$. When $\g=\mf b_n$, let $T_i$ be the $i$-th column of $T$ from the right. Note that $T_i\in \B$ for $1\leq i\leq k$, where $\B={\bf T}^{\g}(0)$ for $\g={\mf c}_n$ and $\B={\bf T}^{\rm sp}$ for   $\g={\mf b}_n$.

We define a $\g$-crystal structure on ${\bf V}^\g_k$ by identifying $T\in {\bf V}^\g_k$ with $T_1\otimes\cdots\otimes T_k\otimes t_{-k\omega_n}\in {\bf B}^{\otimes k}\otimes T_{-k\omega_n}$. Then the inclusion ${\bf V}^\g_{k}\longrightarrow {\bf V}^\g_{k+1}$ is an embedding of $\g$-crystals, and this induces a $\g$-crystal structure on ${\bf V}^\g$ since $\bigcup_{k\geq 1}{\bf V}^\g_k={\bf V}^\g$. Moreover, we have
\begin{equation*}
{\bf V}^\g=\{\,\tf_{i_1}\cdots\tf_{i_r}\emptyset \,|\,r\geq 0,\, i_1,\ldots,i_r\in I\,\},
\end{equation*}
where $\emptyset$ is the empty tableau with ${\rm wt}(\emptyset)=0$ (see \cite[Section 3.1]{K12}).

\begin{rem}{\rm
Let $T\in {\bf V}^{\g}$ be given. Then $T\in {\bf V}^{\g}_k$ for some $k\geq 1$, and
$\te_nT$ is given by removing $\resizebox{.05\hsize}{!}{\def\lr#1{\multicolumn{1}{|@{\hspace{.6ex}}c@{\hspace{.6ex}}|}{\raisebox{-.25ex}{$#1$}}}\raisebox{-.65ex}
{$\begin{array}[b]{cc}
\cline{1-1}\cline{2-2}
\lr{\ov{n}}&\lr{\ov{n}}\\
\cline{1-1}\cline{2-2}
\end{array}$}}$
(resp. $\boxed{\ov{n}}$) when $\g=\mf c_n$ (resp. $\g=\mf b_n$) on top of a column in $T$ following tensor product rule of crystals in $\B^{\otimes k}$. Note that $\te_nT$ is uniquely determined for all sufficiently large $k$, and $\tf_nT$ can be described in a similar way. For example, \vskip 1mm
\begin{equation*}
\begin{split}
& \te_5 \left(\ \
\resizebox{.15\hsize}{!}{
\text{${\def\lr#1{\multicolumn{1}{|@{\hspace{.6ex}}c@{\hspace{.6ex}}|}{\raisebox{-.3ex}{$#1$}}}\raisebox{-.6ex}
{$\begin{array}{ccccc}
\cline{5-5}
& & & & \lr{\ov{5}} \\
\cline{3-5}
& & \lr{\ov{5}} & \lr{\ov{3}}& \lr{\ov{3}}\\
\cline{1-5}
\lr{\ov{5}} & \lr{\ov{4}} & \lr{\ov{3}} & \lr{\ov{2}}& \lr{\ov{1}}\\
\cline{1-5}
\end{array}$}}$}}\ \ \right) \quad =\quad
\resizebox{.17\hsize}{!}{
\text{${\def\lr#1{\multicolumn{1}{|@{\hspace{.6ex}}c@{\hspace{.6ex}}|}{\raisebox{-.3ex}{$#1$}}}\raisebox{-.6ex}
{$\begin{array}{cccccc}
 & & & & & \\
\cline{4-6}
 & & & \lr{\ov{5}} & \lr{\ov{3}}& \lr{\ov{3}}\\
\cline{2-6}
 &\lr{\ov{5}} & \lr{\ov{4}} & \lr{\ov{3}} & \lr{\ov{2}}& \lr{\ov{1}}\\
\cline{2-6}
\end{array}$}}$}}\quad 
\in {\bf V}^{\mf b_5}.
\end{split}
\end{equation*}
}
\end{rem}
\vskip 2mm

\begin{df}\label{def:parabolic Verma}{\rm
For $\lambda\in \cP_n\cup \cP_n^{\rm sp}$, let
\begin{equation*}
{\bf V}^{\mf g}_{\lambda}= SST_{[\ov{n}]}(\nu)\times {\bf V}^{\g},
\end{equation*}
where $\nu=(a_\ell,\ldots,a_1)'$ for the partition $(a_\ell,\ldots,a_1)$ given in Definition \ref{def:psst}(2).
}
\end{df}\vskip 2mm

Let us regard $SST_{[\ov{n}]}(\nu)$ as a $\g$-crystal, where the $\mf l$-crystal structure is defined as usual and $\varphi_n(T)=-\infty$ for all $T\in SST_{[\ov{n}]}(\nu)$.  
Then we define a $\g$-crystal structure on  ${\bf V}^\g_\la$  by identifying $(V_2,V_1)\in {\bf V}^\g_\la$ with $V_1\otimes V_2$.  
Let $H_\nu$ be the highest weight element in $SST_{[\ov{n}]}(\nu)$ as an $\mf l$-crystal. Then $(H_\nu,\emptyset)$ is the highest weight element in ${\bf V}^\g_\la$, and 
\begin{equation*}
{\bf V}^\g_\la=\{\,\tf_{i_1}\cdots\tf_{i_r}(H_\nu,\emptyset)\,|\,r\geq 0,\, i_1,\ldots,i_r\in I\,\},
\end{equation*}
\cite[Proposition 3.5]{K12}. The crystal ${\bf V}^\g_\la$ can be viewed as a crystal of a maximal parabolic Verma module in the following sense (see \cite[Section 3.2]{K12} for more details).

\begin{thm}
Let $\la\in\cP_n \cup \cP_n^{\rm sp}$. 
For each $k\in\Z_+$, there exists an embedding of $\g$-crystals 
\begin{equation*}
\xymatrixcolsep{2pc}\xymatrixrowsep{0pc}\xymatrix{
\theta_k : B(\omega_\la + k\omega_n)\otimes T_{-L\omega_n}\ \ar@{^{(}->}[r]  & \ {\bf V}^{\g}_{\la}}, 
\end{equation*}
with $L=\langle \omega_\la + k\omega_n, h_n\rangle$ such that ${\rm Im}\theta_k\subset {\rm Im}\theta_{k+1}$ and ${\bf V}^{\g}_{\la}=\bigcup_{k\in\Z_+}{\rm Im}\theta_k$.
\end{thm} 

\begin{rem}{\rm
An explicit characterization of ${\rm Im}\theta_k$ in ${\bf V}^\g_\la$ is also given in terms of a statistic $\Delta$ on ${\bf V}^\g_\la$, which is a combinatorial realization of $\varepsilon_n^\ast$ \cite[Definition 3.10, Theorem 3.11]{K12}.
}
\end{rem}

\subsection{Crystal structure on Lusztig data}

\subsubsection{Type $A_{2n-1}$}\label{subsec:type A crystal}  In this subsection, we consider the crystal of  Lusztig data of type $A_{2n-1}$ associated to a special family of reduced expressions, which will be used to describe the crystal of Lusztig data in types $B_n$ and $C_n$.

Let us fix some notations.
Let $\mathring{\mf g}=\gl_{2n}$ and $\mathring{I}=\{\,1,2,\ldots,2n-1\,\}$. We assume that $\mathring{P}=\bigoplus_{i=1}^{2n}\Z \upepsilon_i$ is the weight lattice, $\mathring{\Pi}=\{\,\upalpha_i= \upepsilon_i- \upepsilon_{i+1}\,|\,i\in \mathring{I}\,\}$ is the set of simple roots, and $\mathring{W}=\langle\, s_i\,|\,i\in \mathring{I}\,\rangle$ is the Weyl group of $\mathring{\g}$. Let $\mathring{w}_0$ denote  the longest element of length $M=(2n-1)n$ in $\mathring{W}$. Let $\sigma$ be the automorphism of the Dynkin diagram of type $A_{2n-1}$ given by $\sigma(\upalpha_i)=\upalpha_{2n-i}$ for $i\in \mathring{I}$.

Given ${\bi}=(i_1,\ldots,i_M)\in R(\mathring{w}_0)$, let $\mathring{\B}_{\bi}=\Z_+^{M}$ denote the crystal of $\bi$-Lusztig data for $U_q^-(\mathring{\g})$.
Let 
\begin{equation*}
\upbeta_1=\upalpha_{i_1},\  \upbeta_2=s_{i_1}(\upalpha_{i_2}), \ \ldots ,\ \upbeta_M=s_{i_1}\cdots s_{i_{M-1}}(\upalpha_{i_{M}}).
\end{equation*}
Given ${\bf c}=(c_1,\ldots,c_M)\in \mathring{\B}_{\bi}$ and $1\leq k\leq M$, let us also write 
\begin{center}
$c_k=c_{ij}$ if $\upbeta_k=\upepsilon_i-\upepsilon_j$ for some $1\leq i<j\leq 2n$.
\end{center} 
For $1\leq r<s\leq 2n$, let ${\bf 1}_{rs}=(c_{ij})\in \mathring{\B}_{{\bi}}$ be given by $c_{ij}=\delta_{ir}\delta_{js}$. 

Let $\Omega$ be a Dynkin quiver of type $A_{2n-1}$. We call a vertex $i\in \mathring{I}$ a sink of $\Omega$ if there is no arrow going out of $i$. For $i\in \mathring{I}$, let $s_i\Omega$ be the quiver given by reversing the arrows which end or start at $i$. 
We say that  ${\bf i}=(i_1,\ldots,i_M)\in R(\mathring{w}_0)$ is adapted to $\Omega$ if $i_1$ is a sink of $\Omega$, and $i_k$ is a sink of $s_{i_{k-1}}\cdots s_{i_2}s_{i_1}\Omega$ for $2\leq k\leq M$.

Suppose that ${\bi}\in R(\mathring{w}_0)$ is adapted to 
a quiver $\Omega_0$ of type $A_{2n-1}$, where 
\begin{equation}\label{eq:orientation with single sink}
\Omega_0 : \quad
\xymatrixcolsep{2pc}\xymatrixrowsep{0pc}\xymatrix{
\bullet \ar@{->}[r] &  \cdots \ar@{->}[r]  & \bullet  & \ar@{->}[l]  \cdots & \ar@{->}[l]  \bullet \\
{}_1 &  & {}_n  &  & {}_{n-1}}
\end{equation}
In this case of $\bi$, the crystal structure on $\mathring{\B}_{{\bi}}$ can be described explicitly as follows thanks to \cite[Theorem 7.1]{Re}.

Let ${\bf c}=(c_{ij})\in \mathring{\B}_{{\bi}}$ be given.
For $i=n$, we have
\begin{equation*}\label{eq:}
\begin{split}
\te_n{\bf c}&=
\begin{cases}
{\bf c}  - {\bf 1}_{n\,n+1}, & \text{if $c_{n\,n+1}>0$},\\
{\bf 0}, & \text{if $c_{n\,n+1}=0$},
\end{cases}\quad \quad
\tf_n{\bf c}= 
{\bf c}  + {\bf 1}_{n\,n+1}.
\end{split}
\end{equation*}
(Here ${\bf 0}$ is a formal symbol, not the zero vector in $\Z_+^M$.)

For $1\leq i\leq n-1$, if we put
\begin{equation*}\label{eq:crytal on B_Omega-1}
\begin{split}
c_k^{(i)}&= 
\begin{cases}
c_{i\,n+1}, & \text{if $k=1$},\\
c_1^{(i)}+\sum_{s=2}^{k} (c_{i\, n+s}-c_{i+1\, n+s-1}), & \text{if $2\leq k\leq n$},\\
c^{(i)}_{n}+ (c_{i\,n}-c_{i+1\,2n}), & \text{if $k=n+1$},\\
c^{(i)}_{n+1}+\sum_{s=1}^{k-n-1} (c_{i\,n-s}-c_{i+1\,n-s+1}),\!\! & \text{if $n+2\leq k\leq 2n-i$},
\end{cases}
\end{split}
\end{equation*}
and
{\allowdisplaybreaks
\begin{equation*}\label{eq:crytal on B_Omega-2}
\begin{split}
c^{(i)}& = \max\{\,c_k^{(i)}\,|\,1\leq k\leq 2n-i\,\},\\
k_0&=\min\{\,1\leq k\leq 2n-i\,|\,c^{(i)}_k=c^{(i)}\,\},\\
k_1&=\max\{\,1\leq k\leq 2n-i\,|\,c^{(i)}_k=c^{(i)}\,\},
\end{split}
\end{equation*}
then we have 
\begin{equation*}\label{eq:crytal on B_Omega-3}
\begin{split}
\te_i{\bf c}&=
\begin{cases}
{\bf c}  - {\bf 1}_{i\,k_0+n} + {\bf 1}_{i+1\,k_0+n}, & \text{if $c^{(i)}>0$ and $k_0\leq n$},\\
{\bf c}  - {\bf 1}_{i\,2n-k_0+1} + {\bf 1}_{i+1\,2n-k_0+1}, & \text{if  $n+1\leq k_0 \leq 2n-i-1$},\\
{\bf c}  - {\bf 1}_{i\,i+1}, & \text{if $k_0=2n-i$},\\
{\bf 0}, & \text{if $c^{(i)}=0$},
\end{cases}\\
\tf_i{\bf c}&=
\begin{cases}
{\bf c}  + {\bf 1}_{i\,k_1+n} - {\bf 1}_{i+1\,k_1+n}, & \text{if $k_1\leq n$},\\
{\bf c}  + {\bf 1}_{i\,2n-k_1+1} - {\bf 1}_{i+1\,2n-k_1+1}, & \text{if $n+1\leq k_1 \leq 2n-i-1$},\\
{\bf c}  + {\bf 1}_{i\,i+1}, & \text{if $k_1=2n-i$.}
\end{cases}\\\end{split}
\end{equation*}}
Note that if $c^{(i)}=0$, then we have $k_0=1$ with $c_{in+1}=0$, and if $k_0\geq n+1$, we have $c^{(i)}>0$.

The diagram automorphism $\sigma$ induces a bijection on $\mathring{\B}_{{\bi}}$, still denoted by $\sigma$, such that $\sigma({\bf c})=(c^\sigma_{ij})$ with 
\begin{equation}\label{eq:sigma on B(infty)}
\begin{split}
c^\sigma_{ij}=c_{2n-j+1\,2n-i+1}\quad (1\leq i<j\leq 2n),\\
\end{split}
\end{equation}
for ${\bf c}=(c_{ij})$.
For $n+1\leq i\leq 2n-1$, we have  
\begin{equation*}
\te_i = \sigma\circ \te_{2n-i}\circ \sigma,\quad \tf_i = \sigma\circ \tf_{2n-i}\circ \sigma.
\end{equation*}

\subsubsection{Types $C_n$ and $B_n$} \
The diagram automorphism $\sigma$ induces an automorphism on $\mathring{P}$ given by $\sigma(\upepsilon_i)=-\upepsilon_{2n-i+1}$ for $1\leq i\leq 2n$. Note that $\sigma\circ s_i \circ \sigma = s_{2n-i}$ for $i\in \mathring{I}$.
Let $- : P \longrightarrow \mathring{P}$ be an embedding given by $\ov{\epsilon_i}=\upepsilon_i-\upepsilon_{2n-i+1}$ for $1\leq i\leq n$ so that 
\begin{equation*}
\ov{\alpha_i} = 
\begin{cases}
\upalpha_i+\upalpha_{2n-i}, & \text{for $1\leq i\leq n-1$},\\
\epsilon\upalpha_n, & \text{for $i=n$}.
\end{cases}
\end{equation*}
Let $- \ : W \rightarrow \mathring{W}$ be an injective group homomorphism 
 given by 
\begin{equation*}
\ov{s_i}=
\begin{cases}
s_i s_{2n-i}, & \text{for $1\leq i<n$},\\
s_n, & \text{for $i=n$},
\end{cases}
\end{equation*}
(see \cite[Corollary 3.4]{FRS}).
By using induction on the length of $w$, one can check that for $w\in W$ and $\mu\in P$
\begin{equation}\label{eq:bar preserves w}
\ov{w}(\ov{\mu})=\ov{w(\mu)}.
\end{equation}

Let $B(\infty)$ and $\mathring{B}(\infty)$ be the crystals of $U^-_q(\g)$ and $U^-_q(\mathring{\g})$, respectively. 
We define the operators $\td{\tt e}_i$ and $\td{\tt f}_i$ for $i\in I$ on $\mathring{B}(\infty)$ by
\begin{equation*}
\td{\tt e}_i = 
\begin{cases}
\te_i\te_{2n-i}, & \text{for $1\leq i\leq n-1$},\\
\te_n^\epsilon, & \text{for $i=n$},
\end{cases}\quad
\td{\tt f}_i = 
\begin{cases}
\tf_i\tf_{2n-i}, & \text{for $1\leq i\leq n-1$},\\
\tf_n^\epsilon, & \text{for $i=n$},
\end{cases}
\end{equation*}
where $\te_i$ and $\tf_i$ are the Kashiwara operators on $\mathring{B}(\infty)$.
We define $\td{\tt e}^\ast_i=\ast \circ \td{\tt e}_i\circ \ast$ and 
$\td{\tt f}^\ast_i=\ast \circ \td{\tt f}_i\circ \ast$ for $i\in \mathring{I}$, where $\ast$ is the bijection on $\mathring{B}(\infty)$ induced from the involution $\ast$ on $U_q^-(\mathring{\g})$.
By \cite[Theorem 5.1]{Kas96}, \cite[Proposition 3.2]{NS03} and \cite[Theorem 2.3.1]{NS05}, one can characterize $B(\infty)$ as a subset of $\mathring{B}(\infty)$ as follows.

\begin{prop}\label{prop:similarity}
There exists a unique injective map 
\begin{equation*}
\xymatrixcolsep{2pc}\xymatrixrowsep{0pc}\xymatrix{
\chi : B(\infty) \ \ar@{^{(}->}[r]  & \ \mathring{B}(\infty)}, 
\end{equation*}
such that for $b\in B(\infty)$ and $i\in I$
\begin{equation*}
\ov{{\rm wt}(b)} = {\rm wt}(\chi(b)),\\
\end{equation*}
\begin{equation}\label{eq:similarity}
\begin{split}
&\chi(\te_i b) = \td{\tt e}_i \chi(b),
\quad 
\begin{cases}
\varepsilon_{i}(b)= \varepsilon_{i}(\chi(b))=\varepsilon_{2n-i}(\chi(b)), & \text{if $i\neq n$},\\
\epsilon\varepsilon_{n}(b)= \varepsilon_{n}(\chi(b)), &  
\end{cases} \\
&\chi(\tf_i b) = \td{\tt f}_i \chi(b)
\quad 
\begin{cases}
\varphi_{i}(b)= \varphi_{i}(\chi(b))=\varphi_{2n-i}(\chi(b)), & \text{if $i\neq n$},\\
\epsilon\varphi_{n}(b)= \varphi_{n}(\chi(b)). &  
\end{cases}\\
\end{split}
\end{equation}
Moreover, \eqref{eq:similarity} holds when $\te_i$, $\tf_i$, $\td{\tt e}_i$, and $\td{\tt f}_i$ are replaced by 
$\te^\ast_i$, $\tf^\ast_i$, $\td{\tt e}^\ast_i$, and $\td{\tt f}^\ast_i$, respectively.
Hence $\chi(B(\infty))=\{\,\td{\tt f}_{i_1}^{m_1}\cdots \td{\tt f}_{i_r}^{m_r} {\bf 1} \,|\,r\geq 0, i_1,\ldots,i_r\in I, m_1,\ldots, m_r\in \Z_+\,\}$ is isomorphic to $B(\infty)$ with respect to $\td{\tt e}_i$ and $\td{\tt f}_i$ for $i\in I$ as a $\g$-crystal, where ${\bf 1}$ denotes the highest weight element in $\mathring{B}(\infty)$.
\end{prop}

Let $w_0$ be the longest element of length $N=n^2$ in the Weyl group $W$ of $\g$.
Given ${\bi}=(i_1,\ldots,i_{{N}})\in R(w_0)$, consider $\B_{\bi}=\Z_+^N$ the crystal of $\bi$-Lusztig data for $U_q^-(\g)$.
Let 
\begin{equation*}
\beta_1=\alpha_{i_1},\  \beta_2=s_{i_1}(\alpha_{i_2}), \ \ldots ,\ \beta_N=s_{i_1}\cdots s_{i_{{N}-1}}(\alpha_{i_{{N}}}).
\end{equation*}
Let $\ov{\bi}$ be the reduced expression of 
$\mathring{w}_0=\ov{w_0}=\ov{s_{i_1}}\cdots \ov{s_{i_N}}\in \mathring{W}$.
Then the image of $\B_{\bi}\cong B(\infty)$ in $\mathring{\B}_{\ov{\bi}}\cong \mathring{B}(\infty)$ under the map $\chi$ in Proposition \ref{prop:similarity} can be characterized as follows.

\begin{thm}\label{thm:image of PBW under chi}
For $\bi\in R(w_0)$, we have
\begin{equation*}
\chi (\B_{\bi}) = \{\,{\bf c}\in \mathring{\B}_{\ov{\bi}}\,|\,\sigma({\bf c})={\bf c},\, \epsilon|c_{k\,2n-k+1}\,(1\leq k\leq n)\,\},
\end{equation*}
where $\sigma$ is the bijection on $\mathring{\B}_{\ov{\bi}}$ given in \eqref{eq:sigma on B(infty)}. 
Moreover, for ${\bf d}=(d_1,\ldots,d_N)\in \B_{\bi}$,
we have 
$\chi({\bf d})={\bf c}=(c_1,\ldots,c_M)$, where 
\begin{equation}\label{eq:image of d under xi}
\begin{split}
c_k=
\begin{cases}
d_j, & \text{if $\upbeta_k\neq \sigma(\upbeta_k)$ and $\upbeta_k+\sigma(\upbeta_k)=\ov{\beta_j}$ for some $1\leq j\leq N$},\\
\epsilon d_j, & \text{if $\upbeta_k = \sigma(\upbeta_k)$ and $\epsilon\upbeta_k=\ov{\beta_j}$ for some $1\leq j\leq N$},
\end{cases}
\end{split}
\end{equation} 
for $1\leq k\leq M$.
\end{thm}
\pf Let $\bi=(i_1,\ldots,i_N)$ be given. 
For convenience, let us write $\ov{i} = (i,2n-i)$ for $1\leq i\leq n-1$ and $\ov{n}=n$
so that 
\begin{equation*}
\ov{\bi}=(\ov{i_1},\ldots,\ov{i_N})=(i'_1,\ldots,i'_M)\in R(\mathring{w}_0).
\end{equation*}
Let $1\leq j\leq N$ be given. If $i_j\neq n$, then $\ov{i_j}=(i_j,2n-i_j)=(i'_k,i'_{k+1})$ in $\ov{\bi}$ for some $1\leq k<M$, and
\begin{equation*}
\begin{split}
\upbeta_{k}&=\ov{s_{i_1}}\cdots \ov{s_{i_{j-1}}}(\upalpha_{i_j}),\\ \quad
\upbeta_{{k+1}}&=\ov{s_{i_1}}\cdots \ov{s_{i_{j-1}}}s_{i_j}(\upalpha_{2n-i_j})
=\ov{s_{i_1}}\cdots \ov{s_{i_{j-1}}}(\upalpha_{2n-i_j}).
\end{split}
\end{equation*}
Since $\sigma\circ \ov{w}\circ \sigma =\ov{w}$ for $w\in W$, we have by \eqref{eq:bar preserves w}
\begin{equation}\label{eq:beta with 2 orbit}
\upbeta_{k}=\sigma(\upbeta_{{k+1}}),\quad
\upbeta_{k}+\sigma(\upbeta_{{k+1}})=\ov{\beta_{j}}.
\end{equation}
If $i_j=n$, then $\ov{i_j}=i'_k=n$ for some $1\leq k\leq M$ and 
\begin{equation}\label{eq:beta with 1 orbit}
\begin{split}
\epsilon\upbeta_{k}&=\epsilon\ov{s_{i_1}}\cdots \ov{s_{i_{j-1}}}(\upalpha_{n})=\ov{\beta_j}.
\end{split}
\end{equation}
  
Let ${\bf d}=(d_1,\ldots,d_N)\in \B_{\bi}$ be given. 
Define a sequence ${\bf d}_N,\ldots, {\bf d}_1$ in $\B_{\bi}$ inductively as follows:
\begin{equation}\label{eq:sequence d-1}
\begin{split}
{\bf d}_{j}&= \tf_{i_{j}}^{d_{j}}
\tf_{i_{j}}^{\ast\varphi_{i_{j}}} \te_{i_{j}}^{\varepsilon_{i_{j}}} {\bf d}_{j+1} \quad (1\leq j\leq N),\\
\end{split}
\end{equation}
where ${\bf d}_{N+1}={\bf 1}$, and 
$\varepsilon_{i_{j}} =\varepsilon_{i_{j}} ({\bf d}_{j+1})$,  
$\varphi_{i_{j}}=\varphi_{i_{j}}({\bf d}_{j+1})$.
By \cite[Proposition 3.4.7]{S94}, we have 
${\bf d}={\bf d}_1$.  
Applying $\chi$ to \eqref{eq:sequence d-1}, we have by Proposition \ref{prop:similarity}
\begin{equation}\label{eq:sequence d-2}
\begin{split}
\chi({\bf d}_{j})&= 
\td{\tt f}_{i_{j}}^{d_{j}}\td{\tt f}_{i_{j}}^{\ast\varphi_{i_{j}}} 
\td{\tt e}_{i_{j}}^{\varepsilon_{i_{j}}} \chi({\bf d}_{j+1})\\
&=
\begin{cases}
\left(\tf_{i_{j}}^{d_{j}}\tf_{i_{j}}^{\ast\varphi_{i_{j}}} \te_{i_{j}}^{\varepsilon_{i_{j}}}\right)
\left(\tf_{2n-i_{j}}^{d_{j}}\tf_{2n-i_{j}}^{\ast\varphi_{i_{j}}} \te_{2n-i_{j}}^{\varepsilon_{i_{j}}}\right) \chi({\bf d}_{j+1}), & \text{if $i_j\neq n$},\\
\tf_{i_{j}}^{\epsilon d_{j}}\tf_{i_{j}}^{\ast\epsilon\varphi_{i_{j}}} \te_{i_{j}}^{\epsilon\varepsilon_{i_{j}}} \chi({\bf d}_{j+1}),& \text{if $i_j=n$},
\end{cases}
\end{split}
\end{equation}
for $1\leq j\leq N$.
By \eqref{eq:similarity}, we have
\begin{equation*}
\begin{split}
&
\begin{cases}
\varepsilon_{i_{j}} 
=\varepsilon_{i_{j}} (\chi({\bf d}_{j+1}))
=\varepsilon_{2n-i_{j}} (\chi({\bf d}_{j+1})), \\
\varphi_{i_{j}}
=\varphi_{i_{j}}(\chi({\bf d}_{j+1}))
=\varphi_{2n-i_{j}}(\chi({\bf d}_{j+1})),
\end{cases}\quad \text{if $i_j\neq n$},\\
&
\begin{cases}
\epsilon\varepsilon_{i_{j}}
=\varepsilon_{i_{j}} (\chi({\bf d}_{j+1})), \\
\epsilon\varphi_{i_{j}} 
=\varphi_{i_{j}}(\chi({\bf d}_{j+1})),
\end{cases}
\hskip 3cm  \quad\text{if $i_j=n$}.
\end{split}
\end{equation*}

Finally it follows from \eqref{eq:beta with 2 orbit}, \eqref{eq:beta with 1 orbit}, \eqref{eq:sequence d-2}, and  \cite[Proposition 3.4.7]{S94} that ${\bf c}:=\chi({\bf d})$ satisfies \eqref{eq:image of d under xi}, which also implies that $\sigma({\bf c})={\bf c}$ and $\epsilon|c_{k\,2n-k+1}$ for $1\leq k\leq n$.  
Conversely, suppose that ${\bf c}\in \mathring{\B}_{\ov{\bi}}$ satisfies 
$\sigma({\bf c})={\bf c}$ and $\epsilon|c_{k\,2n-k+1}$ for $1\leq k\leq n$. 
By the above argument, there exists a unique ${\bf d}\in \B_{\bi}$ such that ${\bf c}=\chi({\bf d})$. This completes the proof.
\qed

\subsection{Embedding into the crystal of Lusztig data}\
Let $\bi\in R(\omega_0)$ be given. We assume that $\ov{\bi}$ is adapted to $\Omega_0$ in \eqref{eq:orientation with single sink}. For ${\bf d}=(d_1,\ldots,d_N)\in \B_{\bi}$, let us write
\begin{equation*}
d_k=
\begin{cases}
d^+_{ii}, & \text{if $\beta_k=\epsilon \epsilon_i$ for some $1\leq i\leq n$},\\
d^+_{ij}, & \text{if $\beta_k=\epsilon_i + \epsilon_j$ for some $1\leq i<j\leq n$},\\
d^-_{ij}, & \text{if $\beta_k=\epsilon_i - \epsilon_j$ for some $1\leq i<j\leq n$},
\end{cases}
\end{equation*}
for $1\leq k\leq N$, and  
\begin{equation*}
{\bf d}^+=(d^+_{ij})_{1\leq i\leq j\leq n},\quad {\bf d}^-=(d^-_{ij})_{1\leq i<j\leq n}.
\end{equation*}
By Theorem \ref{thm:image of PBW under chi}, we have an explicit crystal structure on $\B_{\bi}$ identifying ${\bf d}$ with $\chi({\bf d})={\bf c}=(c_{ij})\in \mathring{\B}_{\ov{\bi}}$, where  
\begin{equation*}\label{eq:d_ij}
\begin{cases}
d^+_{ii}=c_{i\, 2n-i+1}/\epsilon, & \text{for $1\leq i\leq n$},\\
d^+_{ij}=c_{i\,2n- j+1}, & \text{for $1\leq i< j\leq n$},\\
d^-_{ij}=c_{ij},& \text{for $1\leq i<j\leq n$}, 
\end{cases}
\end{equation*}
(see Section \ref{subsec:type A crystal}).
Let 
\begin{equation*}
\begin{split}
\B_{\bi}^+&=\{\,{\bf d}^+\,|\,{\bf d}\in \B_{\bi}\,\}=\Z_+^{\frac{n(n+1)}{2}},\quad 
\B_{\bi}^-=\{\,{\bf d}^-\,|\,{\bf d}\in \B_{\bi}\,\}=\Z_+^{\frac{n(n-1)}{2}}.
\end{split}
\end{equation*}
We regard $\B^\pm_{\bi}$ as $\g$-subcrystals of $\B_{\bi}$, where we assume that $\varphi_n({\bf d}^-)=-\infty$ for ${\bf d}^-\in \B_{\bi}^-$. Then we have the following.

\begin{lem}\label{lem:tensor decomp of B} 
Let $\bi\in R(w_0)$ such that $\ov{\bi}$ is adapted to $\Omega_0$. Then the map 
\begin{equation*}
\xymatrixcolsep{3pc}\xymatrixrowsep{0pc}\xymatrix{
\delta : \B_{\bi}  \ar@{->}[r]  & \ \B_{\bi}^+\otimes\B_{\bi}^- \\
\ \ \ {\bf d}  \ar@{|->}[r] & \  {\bf d}^+\otimes{\bf d}^-} 
\end{equation*}
is an isomorphism of $\g$-crystals.
\end{lem}
\pf Let us keep the notations for $\mathring{\B}_{\ov{\bi}}$ in Section \ref{subsec:type A crystal}.
For ${\bf c}=(c_{ij})\in \mathring{\B}_{\ov{\bi}}$, let
\begin{equation*}
{\bf c}^+=(c_{ij})_{(i,j)\in U},\quad {\bf c}^-=(c_{ij})_{(i,j)\in L},
\end{equation*}
where $U=\{\,(i,j)\,|\,1\leq i \red{\leq} n< j\leq 2n\,\}$ and $L=\{\,(i,j)\,|\,1\leq i < j\leq n, \text{ or } n\red{<} i<j\leq 2n\,\}$.
Let 
\begin{equation*}
\begin{split}
\mathring{\B}_{\bi}^+&=\{\,{\bf c}^+\,|\,{\bf c}\in \mathring{\B}_{\bi}\,\}=\Z_+^{n^2},\quad
\mathring{\B}_{\bi}^-=\{\,{\bf c}^-\,|\,{\bf c}\in \mathring{\B}_{\bi}\,\}=\Z_+^{n^2-n},
\end{split}
\end{equation*}
and regard $\mathring{\B}_{\ov{\bi}}^\pm$ as subcrystals of $\mathring{\B}_{\ov{\bi}}$, where we assume that 
$\varphi_n({\bf c}^-)=-\infty$ for ${\bf c}^-\in \mathring{\B}_{\ov{\bi}}^-$.
By \cite[Theorem 4.2]{K16-2}, the map 
\begin{equation*}
\xymatrixcolsep{3pc}\xymatrixrowsep{0pc}\xymatrix{
\mathring{\delta} : \mathring{\B}_{\ov{\bi}}  \ar@{->}[r]  & \ \mathring{\B}_{\ov{\bi}}^+\otimes \mathring{\B}_{\ov{\bi}}^- \\
\ \ \ {\bf c}  \ar@{|->}[r] &  \ {\bf c}^+\otimes{\bf c}^-} 
\end{equation*}
is an isomorphism of $\mathring{\g}$-crystals. 
By Theorem \ref{thm:image of PBW under chi}, $\chi$ induces a map $\chi : \B^\pm_{\bi} \longrightarrow \mathring{\B}^\pm_{\ov \bi}$ such that $\chi({\bf d}^\pm)=\chi({\bf d})^\pm$, for ${\bf d}\in \B_{\bi}$. 
This implies the following commutative diagram of bijections:
\begin{equation*}
\xymatrixcolsep{4pc}\xymatrixrowsep{2pc}\xymatrix{
\B_{\bi} \ar@{->}^{\!\!\!\!\!\!\delta}[r] \ar@{->}^{\chi}[d] &\ \B_{\bi}^+\otimes \B_{\bi}^- \ar@{->}^{\chi\otimes\chi}[d] \\  
\chi(\B_{\bi})  \ar@{->}^{\!\!\!\!\!\!\mathring{\delta}|_{\chi(\B_{\bi})}}[r] &\ \chi(\B_{\bi})^+ \otimes \chi(\B_{\bi})^-    }
\end{equation*}
Since $\mathring{\delta}$ is an isomorphism, we have by \eqref{eq:similarity}
\begin{equation*}
\begin{split}
\td{\tt e}_i(\chi({\bf d}^+)\otimes \chi({\bf d}^-))
&=
\begin{cases}
\td{\tt e}_i \chi({\bf d}^+) \otimes \chi({\bf d}^-), & 
\text{if $\varphi_i(\chi({\bf d}^+))\geq \varepsilon_i(\chi({\bf d}^-))$}, \\
\chi({\bf d}^+) \otimes \td{\tt e}_i \chi({\bf d}^-), & \text{if
$\varphi_i(\chi({\bf d}^+))<\varepsilon_i(\chi({\bf d}^-))$},
\end{cases}\\
\td{\tt f}_i(\chi({\bf d}^+)\otimes \chi({\bf d}^-))
&=
\begin{cases}
\td{\tt f}_i \chi({\bf d}^+) \otimes \chi({\bf d}^-), & 
\text{if $\varphi_i(\chi({\bf d}^+))> \varepsilon_i(\chi({\bf d}^-))$}, \\
\chi({\bf d}^+) \otimes \td{\tt f}_i \chi({\bf d}^-), & \text{if
$\varphi_i(\chi({\bf d}^+))\leq\varepsilon_i(\chi({\bf d}^-))$}.
\end{cases}\end{split}
\end{equation*}
This implies that $\delta$ is a morphism of $\g$-crystals, and hence an isomorphism.
\qed

\begin{rem}{\rm One can also derive the crystal structure on ${\B}_{\bi}$ presented here and its tensor product decomposition by using a recent work on combinatorial description of $\B_{\bi}$ \cite{SST}.}
\end{rem}

In order to describe an embedding of the crystal of a maximal parabolic Verma module into $\B_{\bi}$,
let us briefly recall some well-known results on combinatorics of semistandard tableaux (cf.\cite{Ful}).

Let $T\in SST_{[\ov{n}]}(\la^\pi)$ be given for $\la\in \cP_n$. For $\ov{i}\in [\ov{n}]$, we define $T \leftarrow \ov{i}$ to be the tableau obtained by applying the Schensted's column insertion of $\ov i$ into $T$ in a reverse way starting from the rightmost column of $T$ so that 
$(T\leftarrow \ov i)\in SST_{[\ov{n}]}(\mu^\pi)$ for some $\mu\supset \la$ obtained by adding a box in a corner of $\la$.  
Let ${\bf M}$ be the set of $[\ov{n}]\times [\ov{n}]$ matrices ${\bf m}=(m_{\ov i \ov j})$ such that $m_{\ov i \ov j} \in \Z_+$ and $\sum_{\ov i, \ov j}m_{\ov i \ov j}<\infty$.
Let ${\bf m}\in \bf M$ be given.
For $\ov i\in [\ov n]$, we identify the $\ov i$th row of ${\bf m}$ with a tableau $T\in SST_{[\ov n]}(k)$ for some $k\in \Z_+$ such that $m_{\ov i \ov j}$ is the number of $\ov j$'s in $T$, and let $w({\bf m})_{\ov i}=w(T)$. 
We define
$P({\bf m})= (((w_r \leftarrow w_{r-1})\cdots )\leftarrow w_1)$ if $w({\bf m})_{\ov n}\cdots w({\bf m})_{\ov 1}=w_1\ldots w_r$. 
By the well-known symmetry of RSK correspondence, we have $P({\bf m})=P({\bf m}^t)$ if and only if ${\bf m}={\bf m}^t$, where ${\bf m}^t$ is the transpose of the matrix ${\bf m}$. Hence the map sending ${\bf m}$ to $P({\bf m})$ is a bijection 
\begin{equation*}
\kappa : {\bf M}^{sym} \longrightarrow \bigsqcup_{\mu\in\cP_n}SST_{[\ov{n}]}(\mu^\pi),
\end{equation*}
where ${\bf M}^{sym}$ is the set of symmetric matrices in $\bf M$.

Now, consider 
${\bf V}^\g_\la= SST_{[\ov{n}]}(\nu)\times {\bf V}^{\g}$ for  $\la\in\cP_n \cup \cP_n^{\rm sp}$.
For ${\bf V}=(V_2,V_1)\in {\bf V}^\g_\la$, 
we define ${\bf d}_{\bf V}\in \B_{\bi}$ by  
\begin{equation}\label{eq:def of Phi}
\begin{cases}
({\bf d}_{\bf V})^+ =(d^+_{ij})_{1\leq i\leq j\leq n} & \text{with $d^+_{ij}= m_{\ov{i}\ov{j}}$ for $i< j$ and $d^+_{ii}=m_{\ov{i}\ov{i}}/\epsilon$},\\
({\bf d}_{\bf V})^- =(d^-_{ij})_{1\leq i< j\leq n} & \text{with $d^-_{ij}= v_{ij}$},
\end{cases}
\end{equation}
where $m_{\ov{i}\ov{j}}$ is the $(\ov{i},\ov{j})$-entry of ${\bf m}=\kappa^{-1}(V_1)\in {\bf M}^{sym}$, and $v_{ij}$ is the number of $\ov{i}$'s in the $(n-j+1)$th row of $V_2$ from the top. Note that when $\g=\mf c_n$, we have $V_1\in SST_{[\ov n]}(\mu^\pi)$ for some $\mu\in \cP_{\mf c}$ and hence $2 | m_{\ov i \ov i}$ for $1\leq i\leq n$, where ${\bf m}=(m_{\ov i \ov j})=\kappa^{-1}(V_1)$.

\begin{thm}\label{eq:Phi isomorphism}
Let $\la\in \cP_n\cup \cP_n^{\rm sp}$ and $\bi\in R(w_0)$ such that $\ov{\bi}$ is adapted to $\Omega_0$. Then the map 
\begin{equation*}\label{eq:Phi}
\xymatrixcolsep{4pc}\xymatrixrowsep{0pc}\xymatrix{
\Phi_\la: {\bf V}_\la^\g \otimes T_{r\omega_n} \ar@{->}[r]  & \  \B_{\bi}\otimes T_{\omega_\la} \\
\quad \quad {\bf V}\otimes t_{r\omega_n}  \ar@{|->}[r] & \ {\bf d}_{\bf V}\otimes t_{\omega_\la}}
\end{equation*}
is an embedding of $\g$-crystals, where $r=\langle \la,h_n \rangle$.
\end{thm}
\pf Let us first identify $\bf M$ with $\mathring{\B}^+_{\ov i}$, where ${\bf m}=(m_{\ov i \ov j})\in {\bf M}$ corresponds to ${\bf c}^+=(c_{ij})_{(i,j)\in U}\in \mathring{\B}^+_{\ov i}$ in such a way that
\begin{equation*}
m_{\ov i \ov j} = c_{i\, 2n-j+1}\quad (1\leq i\leq j\leq n).
\end{equation*}
Then ${\bf M}^{sym}_\epsilon:=\{\,{\bf m}\,|\,{\bf m}=(m_{\ov i \ov j})\in {\bf M}^{sym},\ \epsilon|m_{\ov i \ov i}\ (1\leq i\leq n) \,\}$ is equal to $\chi(\B^+_{\bi})$.
Under this identification, the $\mathring{\g}$-crystal (resp. ${\g}$-crystal) structure on $\mathring{\B}^+_{\ov i}$ (resp. ${\B}^+_{i}$) coincides with the one on $\bf M$ (resp. ${\bf M}^{sym}_\epsilon$) defined in \cite{K09}.

Let ${\bf V}=(V_2,V_1)\in {\bf V}^\g_\la$ be given.
 By \cite[Theorem 3.6]{K09} the map ${\bf V}^\g\longrightarrow \B_{\bi}^+$ sending $V_1$ to $({\bf d}_{\bf V})^+$ is an isomorphism of $\g$-crystals. 
Also the map $SST_{[\ov{n}]}(\nu)\longrightarrow \B_{\bi}^-\otimes T_{\omega_\la-r\omega_n}$ sending $V_2$ to $({\bf d}_{\bf V})^-\otimes t_{\omega_\la-r\omega_n}$ is an embedding of $\g$-crystals (see \cite[Proposition 5.2]{K16-2}). Therefore, the map sending ${\bf V}\otimes t_{r\omega_n}$ to ${\bf d}_{\bf V}\otimes t_{\omega_\la}$ is an embedding of $\g$-crystals by Lemma \ref{lem:tensor decomp of B}.
\qed

\begin{ex}\label{ex:embedding-2}
{\rm
Suppose that $\g=\mf c_5$ and $\la=(4,4,3,2,0)\in \cP_5$. By Definition \ref{def:psst}(2), we have $(a_4,a_3,a_2,a_1) =(3,2,1,1)$, and hence $\nu=(3,2,1,1)'=(4,2,1)$.
Let ${\bf V}=(V_2,V_1)\in {\bf V}_\la^{\mf c_5}$ be given by
$$
V_2=\
\resizebox{.13\hsize}{!}
{\def\lr#1{\multicolumn{1}{|@{\hspace{.75ex}}c@{\hspace{.75ex}}|}{\raisebox{-.04ex}{$#1$}}}\raisebox{-.6ex}
{$\begin{array}{ccccccccccccccccccc}
\cline{2-9}
& \lr{\ov{3}} & \lr{\ov{3}} & \lr{\ov 2} & \lr{\ov{1}}   \\
\cline{2-5}
& \lr{\ov{2}} &  \lr{\ov{1}} \\
\cline{2-3}
& \lr{\ov{1}} \\
\cline{2-2}
\end{array}$}}
\quad\quad
V_1=\
\resizebox{.27\hsize}{!}
{\def\lr#1{\multicolumn{1}{|@{\hspace{.75ex}}c@{\hspace{.75ex}}|}{\raisebox{-.04ex}{$#1$}}}\raisebox{-.6ex}
{$\begin{array}{ccccccccccccccccccc}
\cline{8-9}
& & & & & & & \lr{\ov{5}} & \lr{\ov{3}} &  \\
\cline{4-9}
& & & \lr{\ov{5}}  & \lr{\ov{5}}   & \lr{\ov{5}} & \lr{\ov{4}}  & \lr{\ov{4}}  & \lr{\ov{2}}\\
\cline{2-9}
& \lr{\ov{4}} & \lr{\ov{4}}  & \lr{\ov{3}} & \lr{\ov{2}} & \lr{\ov{2}} & \lr{\ov{2}} & \lr{\ov{2}} &  \lr{\ov{1}}\\
\cline{2-9}
\end{array}$}}.
$$ \vskip 2mm 
Since 
\begin{equation*}
{\bf m}
=
\begin{bmatrix}
m_{\ov 5 \ov 5} & m_{\ov 5\ov 4} & m_{\ov 5\ov 3} & m_{\ov 5 \ov 2} & m_{\ov 5\ov 1} \\
m_{\ov 4 \ov 5} & m_{\ov 4\ov 4} & m_{\ov 4\ov 3} & m_{\ov 4\ov 2} & m_{\ov 4\ov 1} \\
m_{\ov 3\ov 5} & m_{\ov 3\ov 4} & m_{\ov 3\ov 3} & m_{\ov 3\ov 2} & m_{\ov 3\ov 1} \\
m_{\ov 2\ov 5} & m_{\ov 2\ov 4} & m_{\ov 2\ov 3} & m_{\ov 2\ov 2} & m_{\ov 2\ov 1} \\
m_{\ov 1\ov 5} & m_{\ov 1\ov 4} & m_{\ov 1\ov 3} & m_{\ov 1\ov 2} & m_{\ov 1\ov 1}
\end{bmatrix}
=
\begin{bmatrix}
2 & 1 & 1 & 0 & 0 \\
1 & 0 & 0 & 3 & 0 \\
1 & 0 & 0 & 0 & 1 \\
0 & 3 & 0 & 2 & 0 \\
0 & 0 & 1 & 0 & 0 \\
\end{bmatrix}\ \in {\bf M}^{sym}
\end{equation*}
satisfies $\kappa({\bf m})=P({\bf m})=V_1$, we have 
\begin{equation*}
{\bf d}^+_{\bf V} =
\begin{bmatrix}
d^+_{55} & d^+_{45} & d^+_{35} & d^+_{25} & d^+_{15} \\
  & d^+_{44} & d^+_{34} & d^+_{24} & d^+_{14} \\
 &   & d^+_{33} & d^+_{32} & d^+_{13} \\
 &   &   & d^+_{22} & d^+_{12} \\
  &   &   &   & d^+_{11} \\
\end{bmatrix}=
\begin{bmatrix}
1 & 1 & 1 & 0 & 0 \\
  & 0 & 0 & 3 & 0 \\
  &   & 0 & 0 & 1 \\
  &   &   & 1 & 0 \\
  &   &   &   & 0 \\
\end{bmatrix}, \\
\end{equation*}
by \eqref{eq:def of Phi}.
On the other hand, we have 
\begin{equation*}
\begin{split}
{\bf d}^-_{\bf V}&=
\begin{bmatrix}
d^-_{12} & d^-_{13} & d^-_{14} & d^-_{15}   \\
         & d^-_{23} & d^-_{24} & d^-_{25}   \\
         &          & d^-_{34} & d^-_{35}   \\
         &          &          & d^-_{45}   \\
\end{bmatrix}=
\begin{bmatrix}
0 & 1 & 1 & 1   \\
  & 0 & 1 & 1   \\
  &   & 0 & 2  \\
  &   &   & 0   \\
\end{bmatrix},
\end{split}
\end{equation*}
where $d^-_{ij}$ is the number of $\ov i$'s in the $(6-j)$th row in $V_2$ for $1\leq i<j\leq 5$.
}
\end{ex}

\section{Embedding of KN tableaux into Lusztig data}\label{sec:main}
Let us construct an embedding $\Theta_\la$ of the spinor model into the crystal of a maximal parabolic Verma module, which completes the description of an embedding of the crystal of KN tableaux into that of Lusztig data in \eqref{eq:factorization of embedding}.

\subsection{Separation algorithm}
 
For $r\geq 2$, let 
\begin{equation*}
{\bf E}^r=\bigsqcup_{(u_r,\ldots,u_1)\in \Z_+^r} SST_{[\ov n]}(1^{u_r})
\times \cdots\times SST_{[\ov n]}(1^{u_1}).
\end{equation*}
For $(U_r,\ldots,U_1)\in {\bf E}^r$ and $1\leq i\leq r-1$, we define
\begin{equation}\label{eq:Ei and Fi via jdt}
\begin{split}
\mc{X}_i(U_r,\ldots,U_1) &= 
\begin{cases}
(U_r,\ldots, \mc{X}(U_{i+1},U_i),\ldots,U_1),& \text{if $\mc{X}(U_{i+1},U_i)\neq {\bf 0}$},\\
{\bf 0}, &\text{if $\mc{X}(U_{i+1},U_i)={\bf 0}$},
\end{cases}
\end{split}
\end{equation} 
where $\mc{X}(\, \cdot\, ,\, \cdot \, )$ for $\mc{X}=\mc{E},\mc{F}$ is defined in \eqref{eq:E and F via jdt}.

\begin{lem}\label{lem:regularity of E}
${\bf E}^r$ is a regular $\mf{sl}_r$-crystal with respect to $\mc{E}_i$ and $\mc{F}_i$ for $1\leq i\leq r-1$.
\end{lem}
\pf Let ${\bf N}$ be the set of $[\ov n]\times [r]$ matrices ${\bf m}=(m_{\ov i j})_{\ov i\in [\ov n], j\in [r]}$ such that $m_{{\ov i} j}\in \{0,1\}$ and $\sum m_{{\ov i} j}<\infty$. One may identify ${\bf N}$ with ${\bf E}^r$, where the $j$th column of ${\bf m}$ corresponds to $T\in SST_{[\ov n]}(1^{u_j})$ with $m_{{\ov i} j}$ the number of occurrence of $\ov i$ in $T$. We recall that there exists a $(\mf{gl}_n,\mf{gl}_r)$-bicrystal structure on ${\bf N}$ (see for example, \cite[Section 4]{K07}). Then it is not difficult to see that $\mc{E}_i$ and $\mc{F}_i$ for $1\leq i\leq r-1$ in \eqref{eq:Ei and Fi via jdt} coincide with the Kashiwara operators for the (regular) $\gl_r$-crystal structure on ${\bf N}$. This proves the lemma.
\qed

Let $\la\in \cP_n$ be given with $\ell=\la'_1$ and $(a_\ell,\ldots,a_1)$ as in Definition \ref{def:psst}(2). Consider an embedding of sets given by
\begin{equation*}
\xymatrixcolsep{4pc}\xymatrixrowsep{0pc}\xymatrix{
{\bf T}_\la^\g  \ar@{->}[r]  & \  {\bf E}^{2\ell} \\
(T_\ell,\ldots,T_1)  \ar@{|->}[r] & (\underbrace{T_\ell^{\tt L},T_\ell^{\tt R},\ldots,T_1^{\tt L},T_1^{\tt R}}_{2\ell})}.
\end{equation*}
Let ${\bf T}=(T_\ell,\ldots,T_1)\in {\bf T}_\la^\g$ given. By abuse of notation, we identify ${\bf T}$ with its image in ${\bf E}^{2\ell}$.
We define $\ov{\bf T}=(\ov T_{2\ell},\ldots,\ov T_1)\in {\bf E}^{2\ell}$ by
\begin{equation}\label{eq:definition of separation}
\begin{split}
\ov{\bf T}=
\left(\mc{F}^{a_{1}}_{\ell}\cdots \mc{F}^{a_{1}}_{2}\right)\cdots
\left(\mc{F}^{a_{\ell-3}}_{2\ell-4}\, \mc{F}^{a_{\ell-3}}_{2\ell-5}\,\mc{F}^{a_{\ell-3}}_{2\ell-6}\right)
\left(\mc{F}^{a_{\ell-2}}_{2\ell-3}\, \mc{F}^{a_{\ell-2}}_{2\ell-4}\right)
\mc{F}^{a_{\ell-1}}_{2\ell-2} \ {\bf T}.
\end{split}
\end{equation}

Recall that ${\rm ht}(T^{\tt tail}_i)=a_i$ for $1\leq i\leq \ell$.
By definition of $\mc{F}_i$ in \eqref{eq:Ei and Fi via jdt}, we first move $a_{\ell-1}$ boxes in $T^{\tt tail}_{\ell-1}$ to the $(2\ell-1)$th column (component) via jeu de taquin, next move $a_{\ell-2}$ boxes in $T^{\tt tail}_{\ell-2}$ to the $(2\ell-2)$th column and so on. 
Hence we have
\begin{equation}\label{eq:height of ovT}
\begin{split}
{\rm ht}(\ov T_i)=
\begin{cases}
b_{i},&\text{for $1\leq i\leq \ell$},\\
b_{i}+a_{i-\ell}, & \text{for $\ell<i\leq 2\ell$}, 
\end{cases}
\end{split}
\end{equation}
where $b_{2i-1}={\rm ht}(T^{\tt R}_i)$ and $b_{2i}={\rm ht}(T^{\tt L}_i)-a_i$ for $1\leq i\leq \ell$.
Put 
\begin{equation*}
\begin{split}
\nu&=(a_\ell,\ldots,a_1)'\in \cP_n,\\
\tau&=(b_1,\ldots,b_{2\ell})'\in\cP_\g.
\end{split}
\end{equation*}
As in Remark \ref{rem:vertical position}, we assume that 
$T^{\tt body}_i$ and $T^{\tt tail}_i$ ($1\leq i\leq \ell$) are separated by a common horizontal line, say $L$. According to this convention, we may assume that the bottoms of $\ov{T}_i$ ($1\leq i\leq \ell$) are placed on $L$, and $\ov{T}_i$ ($\ell<i\leq 2\ell$) are placed with $b_i$ boxes  above $L$ and $a_{i-\ell}$ boxes below $L$. Hence $\ov{\bf T}$ can be viewed as a tableau of the following shape

\begin{equation}\label{eq:shape after separation}
\raisebox{.9ex}{$\eta=$} 
\resizebox{.38\hsize}{!}{
\def\lr#1{\multicolumn{1}{|@{\hspace{.75ex}}c@{\hspace{.75ex}}|}{\raisebox{-.04ex}{$#1$}}}
\def\l#1{\multicolumn{1}{|@{\hspace{.75ex}}c@{\hspace{.75ex}}}{\raisebox{-.04ex}{$#1$}}}
\def\r#1{\multicolumn{1}{@{\hspace{.75ex}}c@{\hspace{.75ex}}|}{\raisebox{-.04ex}{$#1$}}}
\raisebox{-.6ex}
{$\begin{array}{cccccccccc}
\cline{8-9}
& & & & & & &\l{\ \ } & \r{ } & \\
\cline{6-7}
 & & & & & \l{\ \ } & & & \r{ }\\
\cline{5-5} 
& & & & \l{\ \ } &  &  & & \r{\!\!\!\!\!\!\!\!\!\!\!\!\!\!\!\!\!\!\!\!\!\!\!\!\!\!\tau^\pi } \\
\cline{3-4}
& & \l{\ \ } & & & & & & \r{ } &  \\
\cline{2-2}\cdashline{1-10}[0.5pt/1pt]\cline{6-9}
& \l{\ \ } &  & & \r{\ \ }  &  & &  \\ 
\cline{4-5}
& \l{\ \ } & \r{\!\!\!\!\!\!\!\!\!\nu} & & &  & &  \\
& \l{\ \ } & \r{\ \ } & & & & &  \\  
\cline{3-3} 
& \lr{\ \ } & & & &  & &  \\  
\cline{2-2} 
& \l{\leftarrow} & & \!\!\!\!\ell &  \r{\rightarrow} & \l{\leftarrow} & & \!\!\!\!\ell & \r{\rightarrow}
\end{array}$}}\ \ \raisebox{.7ex}{$L$}
\end{equation}
\vskip 2mm
\noindent such that $\ov{T}_i$ is the $i$th column from the right.

\begin{lem}\label{lem:separation} Under the above hypothesis, we have
\begin{itemize}
\item[(1)]  $\ov{\bf T}\in SST_{[\ov{n}]}(\eta)$, where $\eta$ is given in \eqref{eq:shape after separation},

\item[(2)] $\ov{\bf T}\equiv_{\mf l} {\bf T}$, where $b\equiv_{\mf l} b'$ means that $b$ is equivalent to $b'$ as elements of $\mf l$-crystals.

\end{itemize}
\end{lem}
\pf 
Let $\ov{\bf T}=(\ov T_{2\ell},\ldots,\ov T_1)$ be in \eqref{eq:definition of separation}, which is a tableau of shape $\eta$ \eqref{eq:shape after separation}.
It is clear that (2) holds since $\mc{F}_i$ preserves Knuth equivalence. So it suffices to show (1).

For simplicity, we say that $(\ov{T}_{i+1},\ov{T}_i)$ is semistandard when the subtableau of $\ov{\bf T}$ consisting of the columns $(\ov{T}_{i+1},\ov{T}_i)$ is $[\ov{n}]$-semistandard. We will show that $(\ov{T}_{i+1},\ov{T}_i)$ is semistandard for $1\leq i\leq 2\ell-1$. 
We use induction on $\ell$ and may assume that $\ell>1$ since it is clear when $\ell=1$. 

{\em Case 1.} $\ell=2$. Let ${\bf T}=(T_2,T_1)$ be given. Recall that ${\rm ht}(T^{\tt R}_1)=b_1$, ${\rm ht}(T^{\tt L}_1)=a_1+b_2$, ${\rm ht}(T^{\tt R}_2)=b_3$, and ${\rm ht}(T^{\tt L}_2)=a_2+b_4$ (cf.~\eqref{eq:height of ovT}).
We have
\begin{equation*}
\ov{\bf T}=(\ov{T}_4,\ov{T}_3,\ov{T}_2,\ov{T}_1)
=\mc{F}_2^{a_1}(T_2^{\tt L},T_2^{\tt R},T_1^{\tt L},T_1^{\tt R})
=(T_2^{\tt L},\mc{F}^{a_1}(T_2^{\tt R},T_1^{\tt L}),T_1^{\tt R}).
\end{equation*}


Let $1\leq j_1<j_2<\cdots<j_{b_2}\leq a_1+b_2$ be a sequence determined inductively as follows:
\begin{equation}\label{eq:aux-2}
\begin{split}
& j_1=\min\{\,j\,|\,1\leq j \leq a_1+b_2,\ T_1^{\tt L}(j)\leq T_1^{\tt R}(1)\,\},  \\
& j_k=\min\{\,j\,|\,j_{k-1} < j \leq a_1+b_2,\ T_1^{\tt L}(j)\leq T_1^{\tt R}(k)\,\}\quad (2\leq k\leq b_2).
\end{split}
\end{equation}
Then ${}^{\tt L}T_1$ is the tableau with entries $T_1^{\tt L}(k)$ for $k\in \{j_1,\ldots,j_{b_2}\}$ (see Remark \ref{rem:equivalent algorithm}).
Next, we define a sequence $1\leq i_1< i_2< \cdots < i_{b_3}\leq a_1+b_3$ inductively as follows:
\begin{equation}\label{eq:aux-1}
\begin{split}
& i_{b_3} =\max\{\,i\,|\,1\leq i\leq a_1+b_3,\ T_2^{\tt R}(b_3)\leq T_1^{\tt L}(i) \,\}, \\
& i_k =\max\{\,i\,|\, 1\leq i< i_{k+1},\ T_2^{\tt R}(k)\leq T_1^{\tt L}(i) \,\}\quad (1\leq k\leq b_3-1).
\end{split}
\end{equation}
Note that such $i_1,\ldots,i_{b_3}$ do exist since $T_2^{\tt R}(k)\leq {}^{\tt L}T_1(k)=T_1^{\tt L}(j_k)$ for $1\leq k\leq b_3$ (see Definition \ref{def:psst} (1)(ii)), and hence $i_k\geq j_k$ for for $1\leq k\leq b_3$. 

Let $X=\{i_1,\ldots,i_{b_3}\}\cup \{a_1+b_3+1, a_1+b_3+2,\ldots, a_1+b_2 \}$. 
Then by definition of $\mc{F}^{a_1}$, it is not difficult to see that $\ov{T}_2$ is the tableau with the entries $T_1^{\tt L}(k)$ for $k\in X$, while $\ov{T}_3$ is the tableau with the entries of $T_2^{\tt R}$ and $T_1^{\tt L}(k)$ for $1\leq k\leq a_1+b_3$ with $k\not\in \{i_1,\ldots,i_{b_3}\}$. It is clear that $(\ov{T}_3,\ov{T}_2)$ is semistandard.

Since $i_k\geq j_k$ for $1\leq k\leq b_3$, we have $\ov{T}_2(k)=T^{\tt L}_1(i_k)\leq T^{\tt L}_1(j_k)\leq T^{\tt R}_1(k)=\ov{T}_1(k)$ 
for $1\leq k\leq b_3$. We also have $\ov{T}_2(k)=T_1^{\tt L}(k+a_1)\leq T_1^{\tt R}(k)=\ov{T}_1(k)$ for $b_3<k\leq b_2$. So $(\ov{T}_2,\ov{T}_1)$ is semistandard.

Let $1\leq i'_1< i'_2< \cdots < i'_{b_4}\leq a_2+b_4$ be a sequence given by
\begin{equation}
\begin{split}
& i'_1=\min\{\,i\,|\,1\leq i \leq a_2+b_4,\ T_2^{\tt L}(i)\leq T_2^{\tt R}(1)\,\},  \\
& i'_k=\min\{\,i\,|\,i'_{k-1} < i \leq a_2+b_4,\ T_2^{\tt L}(i)\leq T_2^{\tt R}(k)\,\}\quad (2\leq k\leq b_4).
\end{split}
\end{equation}
Recall that ${}^{\tt R}T_2$ is the tableau with entries in $T^{\tt R}_2$ and $T_2^{\tt L}(k)$ for $1\leq k\leq a_2+b_4$ with $k\not\in \{\, i'_1, i'_2, \cdots , i'_{b_4}\,\}$.

By Definition \ref{def:psst} (1)(iii) and our choice of $i_1,\ldots,i_{b_3}$, we have
\begin{equation}\label{eq:i prime}
i'_k \leq i_k + (a_2-a_1)
\end{equation}
for $1\leq k\leq b_4$.
Note that for $1\leq l\leq a_1+b_4$
\begin{equation}\label{eq:ov T 4}
\ov{T}_4(l+a_2-a_1)=T^{\tt L}_2(l + a_2-a_1)\leq {}^{\tt R}T_2(l + a_2-a_1)\leq T_1^{\tt L}(l).
\end{equation}

Suppose that $i_{k}< l < i_{k+1}$ for $0\leq k\leq b_3$ and $l\leq a_1+b_4$, where we assume that $i_0=0$ and $i_{b_3+1}=a_1+b_3+1$. Then $\ov{T}_3(l)=T_1^{\tt L}(l)$, and hence $\ov{T}_4(l+a_2-a_1)\leq \ov{T}_3(l)$ by \eqref{eq:ov T 4}.
Next, suppose that $l=i_k$ for $1\leq k\leq b_3$ and $l\leq a_1+b_4$. Then $\ov{T}_3(l)=\ov{T}_3(i_k)=T_2^{\tt R}(k)$. 
On the other hand, if $i_k +a_2-a_1 = i'_{k'}$ for some $k'$ (which implies that $k'\geq k$ by \eqref{eq:i prime}), then ${}^{\tt R}T_2(i_k + a_2-a_1)=T_2^{\tt R}(k')\leq T_2^{\tt R}(k)$.
Otherwise ${}^{\tt R}T_2(i_k + a_2-a_1) < {}^{\tt R}T_2(i'_k)=T_2^{\tt R}(k)$ by \eqref{eq:i prime}. 
Hence we have $\ov{T}_4(i_k+a_2-a_1)\leq \ov{T}_3(i_k)$.

It follows that $(\ov{T}_4,\ov{T}_3)$ is semistandard. 
Hence $\ov{\bf T}$ is semistandard.

{\em Case 2.} $\ell=3$.  Let ${\bf T}=(T_3,T_2,T_1)$ be given. Note that
\begin{equation*}
\ov{\bf T} = \mc{F}^{a_1}_3\mc{F}^{a_1}_2\mc{F}^{a_2}_4 {\bf T}
\end{equation*}
First consider 
\begin{equation*}\label{eq:ov T1}
(\ov{T}_6,\ov{T}_5,U_4,U_3,\ov{T}_2,\ov{T}_1) = \mc{F}^{a_1}_2\mc{F}^{a_2}_4 {\bf T}. 
\end{equation*}
We have by {\em Case 1} that both $(\ov{T}_2,\ov{T}_1)$ and $(\ov{T}_6,\ov{T}_5)$ are semistandard. Let $1\leq i_1<i_2<\cdots<i_{b_3}\leq a_1+b_3$ be as in \eqref{eq:aux-1}.  
Note that $U_3$ and $\ov{T}_2$ are the same as $\ov{T}_3$ and $\ov{T}_2$ in {\em Case 1}, respectively. 

Let $1\leq j_1<j_2<\cdots<j_{b_5}\leq a_2+b_5$ be a sequence defined inductively as follows:
\begin{equation}\label{eq:aux-3}
\begin{split}
& j_{b_5} =\max\{\,j\,|\,1\leq j\leq a_2+b_5,\ T_3^{\tt R}(b_5)\leq T_2^{\tt L}(j) \,\}, \\
& j_k =\max\{\,j\,|\, 1\leq j< j_{k+1},\ T_3^{\tt R}(k)\leq T_2^{\tt L}(j) \,\}\quad (1\leq k\leq b_5-1).
\end{split}
\end{equation}
Put $Y=\{j_1,\ldots,j_{b_5}\}\cup \{a_2+b_5+1, a_2+b_5+2,\ldots, a_2+b_4 \}$. 
By the same arguments for $(\ov{T}_3, \ov{T}_2)$ in {\it Case 1},
$U_4$ is the tableau with the entries $T_2^{\tt L}(k)$ for $k\in Y$, while $\ov{T}_5$ is the tableau with the entries of $T_3^{\tt R}$ and $T_2^{\tt L}(k)$ for $1\leq k\leq a_2+b_5$ with $k\not\in Y$.

Next consider 
\begin{equation*}\label{eq:ov T2}
\ov{\bf T} = \mc{F}^{a_1}_3 (\ov{T}_6,\ov{T}_5,U_4,U_3,\ov{T}_2,\ov{T}_1)=(\ov{T}_6,\ov{T}_5,\mc{F}^{a_1}(U_4,U_3),\ov{T}_2,\ov{T}_1).
\end{equation*}
Note that
\begin{equation}\label{eq:aux-4}
U_4(k)\leq T^{\tt R}_2(k)\quad (1\leq k\leq b_4),
\end{equation}
by the same argument for the semistandardness of $(\ov{T}_2,\ov{T}_1)$ in {\em Case 1}. Then we first have
$\ov{T}_3(k)\leq T^{\tt R}_2(k)$ for $1\leq k\leq b_4$ by \eqref{eq:aux-4} and definition of $\mc{F}^{a_1}(U_4,U_3)$, and hence $\ov{T}_3(k)\leq T^{\tt R}_2(k)\leq T_1^{\tt L}(i_k)=\ov{T}_2(k)$ for $1\leq k\leq b_4$. 
We also have $\ov{T}_3(k)=U_3(k+a_1)\leq \ov{T}_2(k)$ for $b_4<k\leq b_3$ by the semistandardness of $(\ov{T}_3,\ov{T}_2)$ in {\em Case 1}.
So $(\ov{T}_3,\ov{T}_2)$ is semistandard. It is clear that $(\ov{T}_4,\ov{T}_3)=\mc{F}^{a_1}(U_4,U_3)$ is semistandard by \eqref{eq:aux-4} and definition of $\mc{F}$.

So it remains to show that $(\ov{T}_5, \ov{T}_4)$ is semistandard.
Applying the argument in {\em Case 1} to $(T_2,T_1)$, $(T_2^{\tt L},U_3)$ is semistandard, that is, 
\begin{equation}\label{eq:aux-5}
T_2^{\tt L}(k+a_2-a_1)\leq U_3(k)\quad (1\leq k\leq a_1+b_4).
\end{equation}
Let $1\leq j'_1<j'_2<\cdots<j'_{b_4}\leq a_1+b_4$ be a sequence defined inductively as follows:
\begin{equation}\label{eq:aux-5-2}
\begin{split}
& j'_{b_4} =\max\{\,j\,|\,1\leq j\leq a_1+b_4,\ U_4(b_4)=T^{\tt L}_2(j_{b_4})\leq U_3(j) \,\}, \\
& j'_k =\max\{\,j\,|\, 1\leq j< j'_{k+1},\ U_4(k)=T^{\tt L}_2(j_{k})\leq U_3(j) \,\}\quad (1\leq k\leq b_4-1).\\
\end{split}
\end{equation}
Such $j'_1,\ldots,j'_{b_4}$ exist by \eqref{eq:aux-4}. 
Let $Y'=\{j'_1,\ldots,j'_{b_4}\}\cup \{a_1+b_4+1, a_1+b_4+2,\ldots, a_1+b_3 \}$. 
By definition of $\mc{F}^{a_1}$,
$\ov{T}_3$ is the tableau with the entries $U_3(k)$ for $k\in Y'$, while $\ov{T}_4$ is the tableau with the entries of $U_4$ and $U_3(k)$ for $1\leq k\leq a_1+a_4$ with $k\not\in Y'$.  

For $1\leq k\leq b_4$, we have $U_4(k)=T_2^{\tt L}(j_k)\leq U_3(j_k-(a_2-a_1))$ by \eqref{eq:aux-5}, which implies from \eqref{eq:aux-5-2} that
\begin{equation}\label{eq:aux-6}
j_k\leq (a_2-a_1) + j'_k. 
\end{equation}
Note that for $a_2-a_1< l\leq a_2+b_5$, we have
\begin{equation}\label{eq:aux-8}
\ov{T}_4(l-(a_2-a_1))=
\begin{cases}
T_2^{\tt L}(j_{k'}), & \text{if $l-(a_2-a_1)=j'_{k'}$ for some $k'$},\\
U_3(l-(a_2-a_1)), & \text{otherwise}.
\end{cases}
\end{equation}

Suppose that $l=j_k$ for $1\leq k\leq b_5$ with $j_k> a_2-a_1$. Then we have 
\begin{equation}\label{eq:aux-7}
\ov{T}_5(j_k)=T^{\tt R}_3(k)\leq T^{\tt L}_2(j_k)\leq U_3(j_k-(a_2-a_1))
\end{equation}
by \eqref{eq:aux-3} and \eqref{eq:aux-5}. 
If $j_k-(a_2-a_1)=j'_{k'}$ in \eqref{eq:aux-8}, then we have $k'\leq k$ by \eqref{eq:aux-6}.
By \eqref{eq:aux-8} and \eqref{eq:aux-7}, we have $\ov{T}_5(j_k)\leq \ov{T}_4(j_k-(a_2-a_1))$. 

Next, suppose that $l> a_2-a_1$ and $j_{k}<l<j_{k+1}$ with $0\leq k\leq b_5$, where we assume that  $j_0=0$ and $j_{b_5+1}=a_2+b_5+1$. 
If $l-(a_2-a_1)=j'_{k'}$ for some $k'$, then we also have $k'\leq k$ by \eqref{eq:aux-6}, and 
\begin{equation*}
\ov{T}_5(l)=T^{\tt L}_2(l)< T^{\tt L}_2(j_k)\leq T^{\tt L}_2(j_{k'})=\ov{T}_4(l-(a_2-a_1)).  
\end{equation*}
Otherwise, we have
\begin{equation*}
\ov{T}_5(l)=T^{\tt L}_2(l)\leq U_3(l-(a_2-a_1))=\ov{T}_4(l-(a_2-a_1))
\end{equation*}
So $(\ov{T}_5,\ov{T}_4)$ is semistandard, and hence $\ov{\bf T}$ is semistandard.

{\em Case 3.} $\ell>3$. 
Let
\begin{equation}\label{eq:operator X}
\begin{split}
\mc{Y}_1&= 
\left( \mc{F}_3^{a_1}\mc{F}_2^{a_1} \right)\cdots 
\left( \mc{F}_{2\ell-5}^{a_{\ell-3}}\mc{F}_{2\ell-6}^{a_{\ell-3}} \right)
\left( \mc{F}_{2\ell-3}^{a_{\ell-2}}\mc{F}_{2\ell-4}^{a_{\ell-2}} \right) \mc{F}_{2\ell-2}^{a_{\ell-1}},\\
\mc{Y}_2&= 
\left(\mc{F}^{a_{1}}_{\ell}\cdots \mc{F}^{a_{1}}_{4}\right)\cdots
\left(\mc{F}^{a_{\ell-4}}_{2\ell-5}\, \mc{F}^{a_{\ell-4}}_{2\ell-6}\right)
\mc{F}^{a_{\ell-3}}_{2\ell-4}.
\\
\end{split}
\end{equation}
We have 
\begin{equation}\label{eq:aux-9}
\ov{\bf T}=\mc{Y}_2\mc{Y}_1{\bf T}\in {\bf E}^{2\ell},
\end{equation}
since $\mc{F}_i\mc{F}_j=\mc{F}_j\mc{F}_i$ for $|i-j|>1$ by Lemma \ref{lem:regularity of E}. 
Let 
\begin{equation}\label{eq:aux-11}
\mc{Y}_1 {\bf T} =(U_{2\ell},\ldots,U_1).
\end{equation}
It is clear that $(U_{2\ell},U_{2\ell-1})=(\ov{T}_{2\ell},\ov{T}_{2\ell-1})$ and 
$(U_{2},U_{1})=(\ov{T}_{2},\ov{T}_{1})$, and both pairs are semistandard by {\em Case 1}. 
For $2\leq k\leq \ell-1$, we have
\begin{itemize}
\item[$\cdot$] $S_k:=(U_{2k},U_{2k-1})$ is semistandard of shape $\la(a_{k-1},b_{2k-1}-b_{2k},b_{2k})$ by the proof of the semistandardness of $(\ov{T}_4,\ov{T}_3)$ in {\em Case 2}, 

\item[$\cdot$] $S_k\in {\bf T}^{\g}(a_{k-1})$ since $r_{S_k}=0$ by applying Lemma \ref{lem:cal E and F} 
to ${\mc F}_{2k}^{a_k}$ and ${\mc F}_{2k-1}^{a_{k-1}}{\mc F}_{2k-2}^{a_{k-1}}$,

\item[$\cdot$] $S_{k}\prec S_{k-1}$ since $(S_{k}^{\tt R},{}^{\tt L}S_{k-1})$ (resp. $({}^{\tt R}S_{k},S_{k-1}^{\tt L})$) is semistandard by the proof of the semistandardness of $(\ov{T}_5,\ov{T}_4)$ (resp. $(\ov{T}_3,\ov{T}_2)$) in {\em Case 2}.
\end{itemize}
Hence we have
\begin{equation}\label{eq:S}
{\bf S}=(S_{\ell-1},\ldots,S_2)\in {\bf T}^{\g}(\zeta),
\end{equation}
where $\omega_\zeta = \omega_\la - \omega_{n-a_\ell}- \omega_{n-a_1}$, and by induction hypothesis,
\begin{equation}\label{eq:aux-10}
\ov{\bf S} =\mc{Y}_2 {\bf S} =(\ov{T}_{2\ell-2},\ldots,\ov{T}_3)\in {\bf E}^{2\ell-4}
\end{equation}
is semistanard. Finally, $(\ov{T}_{2\ell-1},\ov{T}_{2\ell-2})$ and $(\ov{T}_{3},\ov{T}_{2})$ are semistandard 
by the proofs of the semistandardness of $(\ov{T}_5,\ov{T}_4)$ and $(\ov{T}_3,\ov{T}_2)$ in {\em Case 2}, respectively. Therefore $\ov{T}$ is semistandard by \eqref{eq:aux-9} and \eqref{eq:aux-10}. This completes the induction.
\qed\vskip 2mm

Let $\la=\mu+\sigma_n\in \cP_n^{\rm sp}$ be given with $\ell=\mu_1$. For ${\bf T}=(T_\ell,\ldots,T_1,T_0)\in {\bf T}_\la^\g$, we define 
\begin{equation}\label{eq:definition of separation-spin}
\ov{\bf T}=(\ov T_{2\ell},\ldots,\ov T_1,\ov T_0)\in {\bf E}^{2\ell+1},
\end{equation}
where $(\ov T_{2\ell},\ldots,\ov T_1)$ is obtained as in \eqref{eq:definition of separation} with $(a_\ell,\ldots,a_1)$ given in Definition \ref{def:psst}(2), and $\ov{T}_0=T_0$. 
Since $T_1\prec T_0$, $\ov{\bf T}$ is also semistandard of shape $\eta$ and $\ov{\bf T}\equiv_{\mf l} {\bf T}$ by Lemma \ref{lem:separation}. In this case $\tau=(b_0,b_1,\ldots,b_{2\ell})'\in\cP_\g$ with $b_0={\rm ht}(T_0)$.

\begin{df}
Let $\la\in \cP_n\cup\cP_n^{\rm sp}$.  
For ${\bf T}\in {\bf T}_\la^\g$, let $\ov{\bf T}$ be the tableau of shape $\eta$ given in \eqref{eq:definition of separation} and \eqref{eq:definition of separation-spin}.
We define ${\bf T}^{\tt tail}$  and ${\bf T}^{\tt body}$ to be the tableaux of shape $\nu$ and $\tau^\pi$ given by the subtableaux of $\ov{\bf T}$ located below and above $L$ in \eqref{eq:shape after separation}, respectively.
\end{df}

\begin{cor} Under the above hypothesis, we have
\begin{itemize}
\item[(1)]  $({\bf T}^{\tt tail}, {\bf T}^{\tt body})\in SST_{[\ov{n}]}(\nu)\times SST_{[\ov{n}]}(\tau^\pi)$,

\item[(2)] ${\bf T}^{\tt body}\otimes {\bf T}^{\tt tail} \equiv_{\mf l} {\bf T}$.
\end{itemize}
\end{cor}
\pf (1) follows from Lemma \ref{lem:separation}(1). (2) follows from the fact that the crystal structure on the connected component of $\ov{\bf T}$ does not depend on the choice of reading entries (column reading or row reading) (cf.\cite[Theorem 7.3.6]{HK}).
\qed

\begin{ex}\label{ex:embedding-3}
{\rm
Let ${\bf T}=(T_4,T_3,T_2,T_1)\in {\bf T}^{\mf c_5}_\la$ with $\la=(4,4,3,2,0)\in \cP_5$ 
given by
$$\resizebox{.41\hsize}{!}
{\def\lr#1{\multicolumn{1}{|@{\hspace{.75ex}}c@{\hspace{.75ex}}|}{\raisebox{-.04ex}{$#1$}}}\raisebox{-.6ex}
{$\begin{array}{ccccccccccccccccc}
\\
\cline{14-14}\cline{15-15}
& & & & & & & & & & & & & \lr{\ov{5}} & \lr{\ov{3}} \\
\cline{6-6}\cline{7-7}\cline{10-11}\cline{14-15}
& & & & & \lr{\ov{5}} &  \lr{\ov{5}} & &  & \lr{\ov{5}}&\lr{\ov{4}} & & & \lr{\ov{4}} &  \lr{\ov{2}}\\
\cline{2-3}\cline{6-7}\cline{10-11}\cline{14-15}
& \lr{\ov{4}} & \lr{\ov{3}}& & & \lr{\ov{4}} & \lr{\ov{2}} & & & \lr{\ov{3}} & \lr{\ov{2}}& & & \lr{\ov{2}} &  \lr{\ov{1}}& \\
\cline{2-3}\cdashline{1-5}[0.5pt/1pt]\cline{6-7}\cdashline{8-11}[0.5pt/1pt]
\cline{10-11}\cdashline{12-16}[0.5pt/1pt]\cline{14-15}
& \lr{\ov{3}} & & & & \lr{\ov{2}} & & &  & \lr{\ov{2}} & & & & \lr{\ov{1}} & \\
\cline{2-2}\cline{6-6}\cline{10-10}\cline{14-14}
& \lr{\ov{2}} & & & & \lr{\ov{1}} & \\
\cline{2-2}\cline{6-6}
& \lr{\ov{1}} \\
\cline{2-2}\\
T_4\!\!\!\!\!\!\!\!\!\!\!\!\!\!\!\!\!\! & & & & 
T_3\!\!\!\!\!\!\!\!\!\!\!\!\!\!\!\!\!\! & & & & 
T_2\!\!\!\!\!\!\!\!\!\!\!\!\!\!\!\!\!\! & & & & 
T_1\!\!\!\!\!\!\!\!\!\!\!\!\!\!\!\!\!\! 
\end{array}$}} \quad .
$$ 
Regarding 
$$
{\bf T}=\quad 
\resizebox{.43\hsize}{!}
{\def\lr#1{\multicolumn{1}{|@{\hspace{.75ex}}c@{\hspace{.75ex}}|}{\raisebox{-.04ex}{$#1$}}}\raisebox{-.6ex}
{$\begin{array}{ccccccccccccccccccccc}
\cline{14-14}\cline{16-16}
& & & & & & & & & & & & & \lr{\ov{5}} & & \lr{\ov{3}} \\
\cline{6-6}\cline{8-8}\cline{10-10}\cline{12-12}\cline{14-14}\cline{16-16}
& & & & & \lr{\ov{5}} & & \lr{\ov{5}} &  & \lr{\ov{5}}& & \lr{\ov{4}} & & \lr{\ov{4}} & & \lr{\ov{2}}\\
\cline{2-2}\cline{4-4}\cline{6-6}\cline{8-8}\cline{10-10}\cline{12-12}\cline{14-14}\cline{16-16}
&\lr{\ov{4}} & & \lr{\ov{3}}  & & \lr{\ov{4}} & & \lr{\ov{2}} & & \lr{\ov{3}}& & \lr{\ov{2}}& & \lr{\ov{2}} & &  \lr{\ov{1}}&\\
\cline{2-2}\cline{4-4}\cline{6-6}\cline{8-8}
\cline{10-10}\cline{12-12}\cdashline{1-17}[0.5pt/1pt]\cline{14-14}\cline{16-16}
&\lr{\ov{3}} & & & & \lr{\ov{2}} & & &  & \lr{\ov{2}} & & & & \lr{\ov{1}} & \\
\cline{2-2}\cline{6-6}\cline{10-10}\cline{14-14}
&\lr{\ov{2}} & & & & \lr{\ov{1}} & \\
\cline{2-2}\cline{6-6}
&\lr{\ov{1}} \\
\cline{2-2}\\
&\text{\em 8\!} & & \text{\em 7\!} & &\text{\em 6\!} & &\text{\em 5\!} & &\text{\em 4\!} & &\text{\em 3\!} & &\text{\em 2\!} & &\text{\em 1\!} 
\end{array}$}} \quad  \in {\bf E}^8,
$$ 
\noindent where the numbers below each column tableau denote the indices for the components in ${\bf E}^8$,
we have by \eqref{eq:definition of separation}
$$
\!\!\!\!\!\!\!\!\!\!\!\!\!\!\!\!\!\!\!\!\!\!\!\!\!\!\!\!\!\!
{\bf T}\quad  \stackrel{\mc{F}_6^2}{\longrightarrow} \quad 
\resizebox{.43\hsize}{!}
{\def\lr#1{\multicolumn{1}{|@{\hspace{.75ex}}c@{\hspace{.75ex}}|}{\raisebox{-.04ex}{$#1$}}}\raisebox{-.6ex}
{$\begin{array}{ccccccccccccccccccccc}
\cline{14-14}\cline{16-16}
& & & & & & & & & & & & & \lr{\ov{5}} & & \lr{\ov{3}} \\
\cline{6-6}\cline{8-8}\cline{10-10}\cline{12-12}\cline{14-14}\cline{16-16}
& & & & & \lr{\ov{5}} & & \lr{\ov{5}} &  & \lr{\ov{5}}& & \lr{\ov{4}} & & \lr{\ov{4}} & & \lr{\ov{2}}\\
\cline{2-2}\cline{4-4}\cline{6-6}\cline{8-8}\cline{10-10}\cline{12-12}\cline{14-14}\cline{16-16}
&\lr{\ov{4}} & & \lr{\ov{4}}  & & \lr{\ov{2}} & & \lr{\ov{2}} & & \lr{\ov{3}}& & \lr{\ov{2}}& & \lr{\ov{2}} & &  \lr{\ov{1}}& \\
\cline{2-2}\cline{4-4}\cline{6-6}\cline{8-8}
\cline{10-10}\cline{12-12}\cdashline{1-17}[0.5pt/1pt]\cline{14-14}\cline{16-16}
&\lr{\ov{3}} &  & \lr{\ov{3}} & &  & & &  & \lr{\ov{2}} & & & & \lr{\ov{1}} & \\
\cline{2-2}\cline{4-4}\cline{10-10}\cline{14-14}
&\lr{\ov{2}} & & \lr{\ov{1}} & &   & \\
\cline{2-2}\cline{4-4}
&\lr{\ov{1}} \\
\cline{2-2}
\end{array}$}} 
$$ 
$$
\quad \quad  \!\!\!\!\!\!\!\!\!\!\!\!\!\!\!\!\!\!\!\!\!\!\!\!\!\!\!\!\!\!
\stackrel{\mc{F}_5\mc{F}_4}{\longrightarrow} \quad 
\resizebox{.43\hsize}{!}
{\def\lr#1{\multicolumn{1}{|@{\hspace{.75ex}}c@{\hspace{.75ex}}|}{\raisebox{-.04ex}{$#1$}}}\raisebox{-.6ex}
{$\begin{array}{ccccccccccccccccccccc}
\cline{14-14}\cline{16-16}
& & & & & & & & & & & & & \lr{\ov{5}} & & \lr{\ov{3}} \\
\cline{6-6}\cline{8-8}\cline{10-10}\cline{12-12}\cline{14-14}\cline{16-16}
& & & & & \lr{\ov{5}} & & \lr{\ov{5}} &  & \lr{\ov{5}}& & \lr{\ov{4}} & & \lr{\ov{4}} & & \lr{\ov{2}}\\
\cline{2-2}\cline{4-4}\cline{6-6}\cline{8-8}\cline{10-10}\cline{12-12}\cline{14-14}\cline{16-16}
&\lr{\ov{4}} & & \lr{\ov{4}}  & & \lr{\ov{3}} & & \lr{\ov{2}} & & \lr{\ov{2}}& & \lr{\ov{2}}& & \lr{\ov{2}} & &  \lr{\ov{1}}&\\
\cline{2-2}\cline{4-4}\cline{6-6}\cline{8-8}
\cline{10-10}\cline{12-12}\cdashline{1-17}[0.5pt/1pt]\cline{14-14}\cline{16-16}
&\lr{\ov{3}} &  & \lr{\ov{3}} & & \lr{\ov 2} & & &  & & & & & \lr{\ov{1}} & \\
\cline{2-2}\cline{4-4}\cline{6-6}\cline{14-14}
&\lr{\ov{2}} & & \lr{\ov{1}} & &   & \\
\cline{2-2}\cline{4-4}
&\lr{\ov{1}} \\
\cline{2-2}
\end{array}$}}
$$ 
$$
\ \! \quad  \!\!\!\!\!\!\!\!\!\!\!\!\!\!\!\!\!\!\!\!\!\!\!\!\!\!\!\!\!\!
\stackrel{\mc{F}_4\mc{F}_3\mc{F}_2}{\longrightarrow} \quad 
\resizebox{.43\hsize}{!}
{\def\lr#1{\multicolumn{1}{|@{\hspace{.75ex}}c@{\hspace{.75ex}}|}{\raisebox{-.04ex}{$#1$}}}\raisebox{-.6ex}
{$\begin{array}{ccccccccccccccccccccc}
\cline{14-14}\cline{16-16}
& & & & & & & & & & & & & \lr{\ov{5}} & & \lr{\ov{3}} \\
\cline{6-6}\cline{8-8}\cline{10-10}\cline{12-12}\cline{14-14}\cline{16-16}
& & & & & \lr{\ov{5}} & & \lr{\ov{5}} &  & \lr{\ov{5}}& & \lr{\ov{4}} & & \lr{\ov{4}} & & \lr{\ov{2}}\\
\cline{2-2}\cline{4-4}\cline{6-6}\cline{8-8}\cline{10-10}\cline{12-12}\cline{14-14}\cline{16-16}
&\lr{\ov{4}} & & \lr{\ov{4}}  & & \lr{\ov{3}} & & \lr{\ov{2}} & & \lr{\ov{2}}& & \lr{\ov{2}}& & \lr{\ov{2}} & &  \lr{\ov{1}}&\\
\cline{2-2}\cline{4-4}\cline{6-6}\cline{8-8}
\cline{10-10}\cline{12-12}\cdashline{1-17}[0.5pt/1pt]\cline{14-14}\cline{16-16}
&\lr{\ov{3}} &  & \lr{\ov{3}} & & \lr{\ov 2} & & \lr{\ov{1}} &  & & & & & & \\
\cline{2-2}\cline{4-4}\cline{6-6}\cline{8-8}
&\lr{\ov{2}} & & \lr{\ov{1}} & &   & \\
\cline{2-2}\cline{4-4}
&\lr{\ov{1}} \\
\cline{2-2}
\end{array}$}}\quad .
$$ 
Therefore, we get 
$$
\ov{\bf T}=\quad 
\resizebox{.27\hsize}{!}
{\def\lr#1{\multicolumn{1}{|@{\hspace{.75ex}}c@{\hspace{.75ex}}|}{\raisebox{-.04ex}{$#1$}}}\raisebox{-.6ex}
{$\begin{array}{ccccccccccccccccccc}
\cline{8-9}
& & & & & & & \lr{\ov{5}} & \lr{\ov{3}} &  \\
\cline{4-9}
& & & \lr{\ov{5}}  & \lr{\ov{5}}   & \lr{\ov{5}} & \lr{\ov{4}}  & \lr{\ov{4}}  & \lr{\ov{2}}\\
\cline{2-9}
& \lr{\ov{4}} & \lr{\ov{4}}  & \lr{\ov{3}} & \lr{\ov{2}} & \lr{\ov{2}} & \lr{\ov{2}} & \lr{\ov{2}} &  \lr{\ov{1}}\\
\cline{2-9}\cdashline{1-11}[0.5pt/1pt]
& \lr{\ov{3}} & \lr{\ov{3}} & \lr{\ov 2} & \lr{\ov{1}}   \\
\cline{2-5}
& \lr{\ov{2}} &  \lr{\ov{1}} \\
\cline{2-3}
& \lr{\ov{1}} \\
\cline{2-2}
\end{array}$}}
$$ \vskip 2mm
\noindent and hence \vskip 2mm
$$
{\bf T}^{\tt tail}=\quad 
\resizebox{.13\hsize}{!}
{\def\lr#1{\multicolumn{1}{|@{\hspace{.75ex}}c@{\hspace{.75ex}}|}{\raisebox{-.04ex}{$#1$}}}\raisebox{-.6ex}
{$\begin{array}{ccccccccccccccccccc}
\cline{2-9}
& \lr{\ov{3}} & \lr{\ov{3}} & \lr{\ov 2} & \lr{\ov{1}}   \\
\cline{2-5}
& \lr{\ov{2}} &  \lr{\ov{1}} \\
\cline{2-3}
& \lr{\ov{1}} \\
\cline{2-2}
\end{array}$}}
\quad\quad
{\bf T}^{\tt body}=\quad 
\resizebox{.27\hsize}{!}
{\def\lr#1{\multicolumn{1}{|@{\hspace{.75ex}}c@{\hspace{.75ex}}|}{\raisebox{-.04ex}{$#1$}}}\raisebox{-.6ex}
{$\begin{array}{ccccccccccccccccccc}
\cline{8-9}
& & & & & & & \lr{\ov{5}} & \lr{\ov{3}} &  \\
\cline{4-9}
& & & \lr{\ov{5}}  & \lr{\ov{5}}   & \lr{\ov{5}} & \lr{\ov{4}}  & \lr{\ov{4}}  & \lr{\ov{2}}\\
\cline{2-9}
& \lr{\ov{4}} & \lr{\ov{4}}  & \lr{\ov{3}} & \lr{\ov{2}} & \lr{\ov{2}} & \lr{\ov{2}} & \lr{\ov{2}} &  \lr{\ov{1}}\\
\cline{2-9}
\end{array}$}}\quad .
$$ 
}
\end{ex}\vskip 5mm 

\begin{rem}{\rm
The $\ell$-tuples $(T^{\tt tail}_\ell,\ldots,T^{\tt tail}_1)$ and $(T^{\tt body}_\ell,\ldots,T^{\tt body}_1)$ form tableaux of shape $\nu$ and $\tau^\pi$, respectively, which are not necessarily $[\ov n]$-semistandard. 
But they become $[\ov{n}]$-semistandard when ${\bf T}$ is equivalent to an $[\ov{n}]$- semistandard tableau of a relatively small shape
as an element of $\mf l$-crystals. 
More precisely, if ${\bf T}\equiv_{\mf l} S$ for some $S\in SST_{[\ov{n}]}(\gamma)$ with $\gamma_1<  \ell$, then  
$(T^{\tt tail}_\ell,\ldots,T^{\tt tail}_1)$ and $(T^{\tt body}_\ell,\ldots,T^{\tt body}_1)$ are $[\ov{n}]$-semistandard \cite[Lemma 4.2]{K15-1}, which implies that $(T^{\tt tail}_\ell,\ldots,T^{\tt tail}_1)={\bf T}^{\tt tail}$ and $(T^{\tt body}_\ell,\ldots,T^{\tt body}_1)={\bf T}^{\tt body}$. This phenomenon is closely related to so-called a {\em stable branching rule}, which is discussed in \cite{K15-1} with more details.
}
\end{rem}

\subsection{Main result}
Let us establish an embedding of ${\bf T}^\g_\la$ into ${\bf V}^\g_\la\otimes T_{r\omega_n}$, which is the last step for the construction of an embedding of ${\bf KN}^\g_\la$ into $\B_{\bi}\otimes T_{\omega_\la}$. 

\begin{thm}\label{thm:Theta isomorphism}
For $\la\in \cP_n\cup \cP_n^{\rm sp}$, the map 
\begin{equation*}\label{eq:Theta}
\xymatrixcolsep{2pc}\xymatrixrowsep{0pc}\xymatrix{
\Theta_\la:  {\bf T}_\la^\g \ar@{->}[r]  & \  {\bf V}_\la^\g \otimes T_{r\omega_n} \\
\quad \quad {\bf T}  \ar@{|->}[r] & \ ({\bf T}^{\tt tail}, {\bf T}^{\tt body})\otimes t_{r\omega_n}}
\end{equation*}
is an embedding of $\g$-crystals, where $r=\langle \la,h_n \rangle$.
\end{thm}
\pf Note that $\Theta_\la$ is injective since the map ${\bf T}\mapsto \ov{\bf T}$ is reversible.
Let ${\bf T}\in {\bf T}_\la^\g$ be given. By Lemma \ref{lem:separation}(2), we have 
$\Theta_\la(\tilde{x}_i{\bf T})=\tilde{x}_i\Theta_\la({\bf T})$ for $x=e, f$ and $i\in I\setminus \{n\}$. So it remains to show that 
\begin{equation}\label{eq:condition for embedding}
\text{if $\tilde{x}_n{\bf T}\neq {\bf 0}$, then $\Theta_\la(\tilde{x}_n{\bf T})=\tilde{x}_n\Theta_\la({\bf T})$.}
\end{equation}
Let us prove \eqref{eq:condition for embedding} when $\g=\mf c_n$ and $\la\in \cP_n$. The proof for $\g=\mf b_n$ and $\la\in \cP_n\cup\cP_n^{\rm sp}$ is very similar, which we leave it to the reader.

Let ${\bf T}=(T_\ell,\ldots,T_1)$ and $\ov{\bf T}=(\ov T_{2\ell},\ldots,\ov T_1)$. First, we associate a sequence $\sigma=(\sigma_\ell,\ldots,\sigma_1)$ given by
\begin{equation}
\sigma_i=
\begin{cases}
+\ ,& \text{if $\tf_n T_i\neq {\bf 0}$ or  $T_i$ has 
$\resizebox{.05\hsize}{!}{\def\lr#1{\multicolumn{1}{|@{\hspace{.6ex}}c@{\hspace{.6ex}}|}{\raisebox{.2ex}{$#1$}}}\raisebox{-.65ex}
{$\begin{array}[b]{cc}
\cline{1-1}\cline{2-2}
\lr{{x}}&\lr{{y}}\\
\cline{1-1}\cline{2-2}
\end{array}$}}$ 
on its top \ $(\ov n < x \leq  y\leq \ov 1)$},\\
-\ ,& 
\text{if $\te_n T_i\neq {\bf 0}$ or $T_i$ has 
$\resizebox{.05\hsize}{!}{\def\lr#1{\multicolumn{1}{|@{\hspace{.6ex}}c@{\hspace{.6ex}}|}{\raisebox{.1ex}{$#1$}}}\raisebox{-.65ex}
{$\begin{array}[b]{cc}
\cline{1-1}\cline{2-2}
\lr{\ov{n}}&\lr{\ov{n}}\\
\cline{1-1}\cline{2-2}
\end{array}$}}$ 
on its top},\\
\ \cdot\ \, ,& \text{otherwise}.
\end{cases}
\end{equation}
In $\sigma=(\sigma_\ell,\ldots,\sigma_1)$, we replace a
pair $(\sigma_{j},\sigma_{i})=(-,+)$, where $j>i$ and
$\sigma_k=\,\cdot\,$ for $j>k>i$, with $(\,\cdot\,,\,\cdot\,)$, and repeat
this process as far as possible until we get a sequence with no $-$
placed to the left of $+$, which we denote by 
$\sigma^{\rm red}=(\sigma^{\rm red}_\ell,\ldots,\sigma^{\rm red}_1)$.
Recall from tensor product rule of crystals that if $i$ is the largest (resp. smallest) index such that $\sigma^{\rm red}_i=-$ (resp. $\sigma^{\rm red}_i=+$), then $\te_n {\bf T}=(T_\ell,\ldots,\te_nT_i,\ldots,T_1)$ (resp. $\tf_n {\bf T}=(T_\ell,\ldots,\tf_nT_i,\ldots,T_1)$).

Next let $\ov{\sigma}=(\ov{\sigma}_\ell,\ldots,\ov{\sigma}_1)$ be a sequence associated to $\ov{\bf T}$ given by
\begin{equation}\label{eq:sign for separation}
\ov{\sigma}_i=
\begin{cases}
+\ ,& \text{if $(\ov{T}_{2i},\ov{T}_{2i-1})$ has 
$\resizebox{.05\hsize}{!}{\def\lr#1{\multicolumn{1}{|@{\hspace{.6ex}}c@{\hspace{.6ex}}|}{\raisebox{.2ex}{$#1$}}}\raisebox{-.65ex}
{$\begin{array}[b]{cc}
\cline{1-1}\cline{2-2}
\lr{{x}}&\lr{{y}}\\
\cline{1-1}\cline{2-2}
\end{array}$}}$ 
on its top \ $(\ov n < x \leq  y\leq \ov 1)$},\\
-\ ,& 
\text{if $(\ov{T}_{2i},\ov{T}_{2i-1})$ has 
$\resizebox{.05\hsize}{!}{\def\lr#1{\multicolumn{1}{|@{\hspace{.6ex}}c@{\hspace{.6ex}}|}{\raisebox{-.0ex}{$#1$}}}\raisebox{-.65ex}
{$\begin{array}[b]{cc}
\cline{1-1}\cline{2-2}
\lr{\ov{n}}&\lr{\ov{n}}\\
\cline{1-1}\cline{2-2}
\end{array}$}}$
on its top},\\
\ \cdot\ \,,& \text{otherwise},
\end{cases}
\end{equation}
and define $\ov{\sigma}^{\rm red}=(\ov{\sigma}^{\rm red}_\ell,\ldots,\ov{\sigma}^{\rm red}_1)$ in the same way as $\sigma^{\rm red}$.

We will prove that ${\sigma}^{\rm red}=\ov{\sigma}^{\rm red}$, which implies \eqref{eq:condition for embedding}. We assume that $\ell\geq 2$ and use induction on $\ell$. We separate the cases as in the proof of Lemma \ref{lem:separation}(1) and keep the notations there. Note that $b_{2i-1}=b_{2i}$ for $1\leq i\leq \ell$ since $\g=\mf{c}_n$.
We denote by 
$\resizebox{.07\hsize}{!}{\def\lr#1{\multicolumn{1}{|@{\hspace{.6ex}}c@{\hspace{.6ex}}|}{\raisebox{.1ex}{$#1$}}}\raisebox{-.65ex}
{$\begin{array}[b]{cc}
\cline{1-1}\cline{2-2}
\lr{{x_i}}&\lr{{y_i}}\\
\cline{1-1}\cline{2-2}
\end{array}$}}$
the domino on top of $T_i$ for $1\leq i\leq\ell$, and by
$\resizebox{.07\hsize}{!}{\def\lr#1{\multicolumn{1}{|@{\hspace{.6ex}}c@{\hspace{.6ex}}|}{\raisebox{.1ex}{$#1$}}}\raisebox{-.65ex}
{$\begin{array}[b]{cc}
\cline{1-1}\cline{2-2}
\lr{{z_i}}&\lr{{w_i}}\\
\cline{1-1}\cline{2-2}
\end{array}$}}$
the domino on top of $(\ov T_{2i},\ov{T}_{2i-1})$ for $1\leq i\leq\ell$.\vskip 2mm

{\em Case 1.} $\ell=2$.  
Suppose that $\sigma_2=+$, that is, $\ov n < x_2\leq y_2$. 
If $b_3< b_2$, then 
$x_2=z_2$, $x_2\leq w_2\leq y_2$, and $x_1=z_1$, $y_1=w_1$, which implies that $\sigma_2=\ov\sigma_2=+$ and $\sigma_1=\ov\sigma_1$.
If $b_3=b_2$, then $x_2=z_2\leq x_1$, $x_2\leq w_2\leq y_2$, and $x_1\leq z_1$, $y_1=w_1$, which implies that $\sigma_i=\ov\sigma_i=+$ for $i=1,2$.

Suppose that $\sigma_2=-$, that is, $x_2 = y_2=\ov n$. 
In this case, we have $x_2=y_2=z_2=w_2=\ov n$ and $x_1=z_1$, $y_1=w_1$, which implies that 
$\sigma_2=\ov\sigma_2=-$ and $\sigma_1=\ov\sigma_1$.

Suppose that $\sigma_2=\ \cdot\ $, that is, $x_2 =\ov n<y_2$.
If $b_3<b_3$, then it is clear that $x_2=z_2=\ov n<w_2\leq y_2$ and $x_1=z_1$, $y_1=w_1$, which implies that $\sigma_i=\ov\sigma_i$ for $i=1,2$. 
Assume that $b_3=b_2$.
If $\ov n< x_1$, then $x_2=z_2=\ov n<w_2\leq y_2$, $x_1\leq z_1$, and $y_1=w_1$, which implies that $\sigma_2=\ov\sigma_2=\ \cdot\ $ and $\sigma_1=\ov\sigma_1=+$. 
If $\ov n=x_1$, then $\ov n < y_1$ (since $\ov n< y_2$ and $T_2\prec T_1$), $x_2=z_2=w_2=\ov n$, and $\ov n< z_1\leq w_1$, which implies that $\sigma=(\ \cdot\ , \ \cdot\ )$ and $\ov\sigma = (-,+)$ and hence $\sigma^{\rm red}=\ov{\sigma}^{\rm red}$.
\vskip 2mm

{\em Case 2.} $\ell=3$. Let $\sigma'=(\sigma'_3,\sigma'_2,\sigma'_1)$ be the sequence associated to 
\begin{equation*}
{\bf U}=(\ov{T}_6,\ov{T}_5,U_4,U_3,U_2,\ov{T}_1) = \mc{F}^{a_2}_4 {\bf T},
\end{equation*}
as in \eqref{eq:sign for separation}. 
Note that $\ov{T}_1={T}^{\tt R}_1$, $U_2={T}^{\tt L}_1$, $U_3=T^{\tt R}_2$, $\ov{T}_6=T^{\tt L}_3$, and $\sigma'_3=\ov{\sigma}_3$, $\sigma'_1={\sigma}_1$. 
By {\em Case 1}, we have 
$\sigma^{\rm red}={\sigma'}^{\rm red}$. 
Recall that $\ov{\bf T}=\mc{F}^{a_1}_3\mc{F}^{a_1}_2{\bf U}$.
Let
$\resizebox{.07\hsize}{!}{\def\lr#1{\multicolumn{1}{|@{\hspace{.6ex}}c@{\hspace{.6ex}}|}{\raisebox{.0ex}{$#1$}}}\raisebox{-.65ex}
{$\begin{array}[b]{cc}
\cline{1-1}\cline{2-2}
\lr{u_2}&\lr{v_2}\\
\cline{1-1}\cline{2-2}
\end{array}$}}$
be the domino on top of $(U_4,U_3)$.

Suppose that $\sigma'_2=+$, that is, $\ov n < u_2\leq v_2$. 
If $b_3<b_2$, then $\ov n < z_2\leq w_2$, $x_1=z_1$, and $y_1=w_1$, which implies  
that $\sigma'_2=\ov\sigma_2=+$ and $\sigma'_1=\ov\sigma_1$. 
Assume that $b_3=b_2$. Note that $\ov n< y_1$ since $\ov n< y_2$ and $T_2\prec T_1$.
If $\ov n< x_1$, then $\ov n < z_2\leq w_2$, $x_1\leq z_1$, and $y_1= w_1$, which implies that $\sigma'_i=\ov\sigma_i=+$ for $i=1,2$.
If $x_1=\ov n$, then we should have $b_5=b_3$, $x_3=x_2=\ov n$, $\ov{n}<y_3$ and hence $z_3=w_3=\ov{n}$ since $T_3\prec T_2\prec T_1$. 
So $\sigma_i=\ \cdot \ $ for $1\leq i\leq 3$ and $\sigma'_3=-$. 
This implies that $z_2=\ov{n}<w_2$, $x_1< z_1$, and $y_1= w_1$.
So we have $\ov{\sigma}=(-,\ \cdot\ , +)$, while $\sigma'=(-, +,\ \cdot\ )$ and $\sigma=(\ \cdot\ ,\ \cdot\ ,\ \cdot\ )$.

Suppose that $\sigma'_2=-$, that is, $u_2 = v_2=\ov n$. 
By the same argument as in {\em Case 1} (the case when $\sigma_2=-$), we have $u_2=v_2=z_2=w_2=\ov n$ and $x_1=z_1$, $y_1=w_1$, which implies that 
$\sigma'_2=\ov\sigma_2=-$ and $\sigma'_1=\ov\sigma_1$.

Suppose that $\sigma'_2=\ \cdot\ $, that is, $u_2 =\ov n<v_2$.
By the same argument as in {\em Case 1} (the case when $\sigma_2=\ \cdot\ $), we have either $(\sigma'_2,\sigma'_1)=(\ov\sigma_2,\ov\sigma_1)$ or 
$(\sigma'_2,\sigma'_1)=(\ \cdot\ ,\ \cdot\ )$ and $(\ov\sigma_2,\ov\sigma_1)=(-,+)$. 

We have $\sigma'^{\rm red}=\ov\sigma^{\rm red}$ in any case and hence $\sigma^{\rm red}=\ov{\sigma}^{\rm red}$.

\vskip 2mm

{\em Case 3.} $\ell>3$. For $1\leq i\leq \ell$, define ${\bf T}^{(i)}$ inductively as follows; 
\begin{equation}
\begin{split}
{\bf T}^{(1)} &= \mc{F}_{2\ell-2}^{a_{\ell-1}} {\bf T},\\
{\bf T}^{(i)} &= \left( \mc{F}_{2\ell-2i+1}^{a_{\ell-i}}\mc{F}_{2\ell-2i}^{a_{\ell-i}} \right){\bf T}^{(i-1)}\quad (2\leq i\leq \ell-1),\\
{\bf T}^{(\ell)} &= \mc{Y}_2{\bf T}^{(\ell-1)}.
\end{split}
\end{equation}
Note that ${\bf T}^{(\ell-1)}=\mc{Y}_1{\bf T}$ and 
${\bf T}^{(\ell)}=\ov{\bf T}$
(see \eqref{eq:operator X}). Let $\sigma^{(i)}=(\sigma^{(i)}_\ell,\ldots,\sigma^{(i)}_1)$ be the sequence associated to ${\bf T}^{(i)}$ defined as in \eqref{eq:sign for separation} for $1\leq i\leq \ell$.

By {\em Case 1}, we have 
$\sigma^{\rm red}=(\sigma^{(1)})^{\rm red}$, where 
$\sigma^{\rm red}_k=(\sigma^{(1)})^{\rm red}_k$ for $\ell-2\leq k\leq 1$. 
By the arguments for $\sigma'^{\rm red}=\ov\sigma^{\rm red}$ in {\em Case 2}, we have 
$(\sigma^{(i)})^{\rm red}=(\sigma^{(i-1)})^{\rm red}$, where 
$(\sigma^{(i)})^{\rm red}_k=(\sigma^{(i-1)})^{\rm red}_k$ for 
$k\not\in \{\ell-i,\ell-i+1\}$ for $2\leq i\leq \ell-1$. 
Hence $\sigma^{\rm red}=(\sigma^{(\ell-1)})^{\rm red}$.

On the other hand, applying induction hypothesis to ${\bf S}$ in \eqref{eq:S} we have
\begin{equation*}
\left(\sigma^{(\ell-1)}_{\ell-1},\ldots,\sigma^{(\ell-1)}_2\right)^{\rm red}
=\left(\sigma^{(\ell)}_{\ell-1},\ldots,\sigma^{(\ell)}_2\right)^{\rm red}.
\end{equation*}
Since $\sigma^{(\ell-1)}_{k}=\sigma^{(\ell)}_{k}$ for $k=1,\ell$, we have 
\begin{equation}\label{eq:signature}
\sigma^{\rm red}=(\sigma^{(\ell-1)})^{\rm red}=(\sigma^{(\ell)})^{\rm red}=\ov\sigma^{\rm red}.
\end{equation}
\qed

\begin{rem}{\rm For the proof of Theorem \ref{thm:Theta isomorphism} when $\g=\mf b_n$, we have to consider a sequence $\sigma=(\sigma_{2\ell},\ldots,\sigma_1)$ or 
$(\sigma_{2\ell},\ldots,\sigma_1,\sigma_0)$ for ${\bf T}\in {\bf T}^{\mf b_n}_\la$ where $\sigma_{2i}$ (resp. $\sigma_{2i-1}$) is associated to $T^{\tt L}_i$ (resp. $T^{\tt R}_i$), and $\sigma_0$ to $T_0\in {\bf T}^{\rm sp}$. We define $\sigma_k= +$ and $-$ if the top entry of the column is greater than $\ov{n}$ and equal to $\ov{n}$, respectively. The other arguments are almost the same.  
}
\end{rem}

Now, combining Theorems \ref{thm:KN_C_to_psst}, \ref{thm:KN_B_to_psst}, \ref{eq:Phi isomorphism}, and \ref{thm:Theta isomorphism}, we have the following.

\begin{thm}\label{thm:main}
Let $\la\in \cP_n\cup \cP_n^{\rm sp}$ and $\bi\in R(w_0)$ such that $\ov{\bi}$ is adapted to $\Omega_0$. The map 
$\Xi_\la:=\Phi_\la\circ\Theta_\la\circ\Psi_\la$
\begin{equation*}
\xymatrixcolsep{4pc}\xymatrixrowsep{3pc}\xymatrix{
{\bf KN}_\la^\g \ar@{->}^{\Xi_\la }[r] \ar@{->}^{\Psi_\la }[d] &\ \B_{\bi}\otimes T_{\omega_\la} \\  
{\bf T}^{\g}_\la  \ar@{->}^{\Theta_\la }[r] &\ {\bf V}^\g_\la \otimes T_{r\omega_n} \ar@{->}^{\Phi_\la }[u]   }
\end{equation*}
is an embedding of $\g$-crystals, where $r=\langle \la,h_n \rangle$. 
\end{thm}

\begin{rem}\label{rem:related works}
{\rm
In \cite{HL}, it is proved that the set of certain KN tableaux of classical type so called {\em marginally large tableaux} is isomorphic to $B(\infty)$ (see also \cite{C}). 
One may define an embedding of ${\bf KN}_\la^\g$ into the crystal of marginally large tableaux by mapping $T$ to $T'\equiv T\otimes H_\gamma$, where $H_\gamma$ is a highest element in ${\bf KN}_\mu^\g$ with $\langle \omega_\gamma,h_i \rangle \gg 0$ for all $i$. 
Unlike in  type $A$, it does not seem to be easy in general to describe $T'$ since the insertion scheme in types $BCD$ \cite{Le02,Le03} is more involved. 
On the other hand, the embedding $\Xi_\la$ in Theorem \ref{thm:main} uses only the Sch\"{u}tzenberger's jeu de taquin sliding and RSK. }
\end{rem}

\subsection{Examples}\label{sec:examples}
Suppose that $\bi\in R(w_0)$ for $\mf c_5$ or $\mf b_5$, and $\ov{\bi}$ is adapted to $\Omega_0$.

(1) Let  \vskip 2mm
$$
T=\quad
\resizebox{.12\hsize}{!}
{\def\lr#1{\multicolumn{1}{|@{\hspace{.75ex}}c@{\hspace{.75ex}}|}{\raisebox{-.04ex}{$#1$}}}\raisebox{-.6ex}
{$\begin{array}{cccc}
\cline{3-3}\cline{4-4}
& & \lr{1} &\lr{4}\\
\cline{2-2}\cline{3-3}\cline{4-4}
& \lr{3} & \lr{5} &\lr{\ov 4}\\
\cline{1-1}\cline{2-2}\cline{3-3}\cline{4-4}
\lr{5}&\lr{\ov 5} & \lr{\ov 5}&\lr{\ov 2}\\
\cline{1-1}\cline{2-2}\cline{3-3}\cline{4-4}
\lr{\ov 3}&\lr{\ov 2} & \lr{\ov 2}&\lr{\ov 1}\\
\cline{1-1}\cline{2-2}\cline{3-3}\cline{4-4}
\end{array}$}}\quad  \in {\bf KN}^{\mf c_5}_{\la}, 
$$ \vskip 3mm
\noindent where $\la=(4,4,3,2,0)\in \cP_5$. 
By Example \ref{ex:embedding-1}, ${\bf T}=\Psi_\la(T)$ is given by
$${\bf T} =\quad  
\resizebox{.42\hsize}{!}
{\def\lr#1{\multicolumn{1}{|@{\hspace{.75ex}}c@{\hspace{.75ex}}|}{\raisebox{-.04ex}{$#1$}}}\raisebox{-.6ex}
{$\begin{array}{ccccccccccccccccc}
\\
\cline{14-14}\cline{15-15}
& & & & & & & & & & & & & \lr{\ov{5}} & \lr{\ov{3}} \\
\cline{6-6}\cline{7-7}\cline{10-11}\cline{14-15}
& & & & & \lr{\ov{5}} &  \lr{\ov{5}} & &  & \lr{\ov{5}}&\lr{\ov{4}} & & & \lr{\ov{4}} &  \lr{\ov{2}}\\
\cline{2-3}\cline{6-7}\cline{10-11}\cline{14-15}
& \lr{\ov{4}} & \lr{\ov{3}}& & & \lr{\ov{4}} & \lr{\ov{2}} & & & \lr{\ov{3}} & \lr{\ov{2}}& & & \lr{\ov{2}} &  \lr{\ov{1}}& \\
\cline{2-3}\cdashline{1-5}[0.5pt/1pt]\cline{6-7}\cdashline{8-11}[0.5pt/1pt]
\cline{10-11}\cdashline{12-16}[0.5pt/1pt]\cline{14-15}
& \lr{\ov{3}} & & & & \lr{\ov{2}} & & &  & \lr{\ov{2}} & & & & \lr{\ov{1}} & \\
\cline{2-2}\cline{6-6}\cline{10-10}\cline{14-14}
& \lr{\ov{2}} & & & & \lr{\ov{1}} & \\
\cline{2-2}\cline{6-6}
& \lr{\ov{1}} \\
\cline{2-2}\\
\end{array}$}} \quad .
$$ 
\noindent By Example \ref{ex:embedding-3}, we have 
$$
\ov{\bf T}=\quad 
\resizebox{.27\hsize}{!}
{\def\lr#1{\multicolumn{1}{|@{\hspace{.75ex}}c@{\hspace{.75ex}}|}{\raisebox{-.04ex}{$#1$}}}\raisebox{-.6ex}
{$\begin{array}{ccccccccccccccccccc}
\cline{8-9}
& & & & & & & \lr{\ov{5}} & \lr{\ov{3}} &  \\
\cline{4-9}
& & & \lr{\ov{5}}  & \lr{\ov{5}}   & \lr{\ov{5}} & \lr{\ov{4}}  & \lr{\ov{4}}  & \lr{\ov{2}}\\
\cline{2-9}
& \lr{\ov{4}} & \lr{\ov{4}}  & \lr{\ov{3}} & \lr{\ov{2}} & \lr{\ov{2}} & \lr{\ov{2}} & \lr{\ov{2}} &  \lr{\ov{1}}\\
\cline{2-9}\cdashline{1-11}[0.5pt/1pt]
& \lr{\ov{3}} & \lr{\ov{3}} & \lr{\ov 2} & \lr{\ov{1}}   \\
\cline{2-5}
& \lr{\ov{2}} &  \lr{\ov{1}} \\
\cline{2-3}
& \lr{\ov{1}} \\
\cline{2-2}
\end{array}$}}
$$ \vskip 2mm
\noindent with
$$
{\bf T}^{\tt tail}=\quad 
\resizebox{.13\hsize}{!}
{\def\lr#1{\multicolumn{1}{|@{\hspace{.75ex}}c@{\hspace{.75ex}}|}{\raisebox{-.04ex}{$#1$}}}\raisebox{-.6ex}
{$\begin{array}{ccccccccccccccccccc}
\cline{2-9}
& \lr{\ov{3}} & \lr{\ov{3}} & \lr{\ov 2} & \lr{\ov{1}}   \\
\cline{2-5}
& \lr{\ov{2}} &  \lr{\ov{1}} \\
\cline{2-3}
& \lr{\ov{1}} \\
\cline{2-2}
\end{array}$}}
\quad\quad
{\bf T}^{\tt body}=\quad 
\resizebox{.27\hsize}{!}
{\def\lr#1{\multicolumn{1}{|@{\hspace{.75ex}}c@{\hspace{.75ex}}|}{\raisebox{-.04ex}{$#1$}}}\raisebox{-.6ex}
{$\begin{array}{ccccccccccccccccccc}
\cline{8-9}
& & & & & & & \lr{\ov{5}} & \lr{\ov{3}} &  \\
\cline{4-9}
& & & \lr{\ov{5}}  & \lr{\ov{5}}   & \lr{\ov{5}} & \lr{\ov{4}}  & \lr{\ov{4}}  & \lr{\ov{2}}\\
\cline{2-9}
& \lr{\ov{4}} & \lr{\ov{4}}  & \lr{\ov{3}} & \lr{\ov{2}} & \lr{\ov{2}} & \lr{\ov{2}} & \lr{\ov{2}} &  \lr{\ov{1}}\\
\cline{2-9}
\end{array}$}}.
$$ \vskip 2mm 
\noindent  
Hence $\Theta_\la({\bf T})={\bf V}\otimes t_{4\omega_5}$ with ${\bf V}=({\bf T}^{\tt tail}, {\bf T}^{\tt body})$. 
Finally, by Example \ref{ex:embedding-2} we have $\Phi_\la({\bf V}\otimes t_{4\omega_5})={\bf d}\otimes t_{\omega_\la}\in \B_{\bi}\otimes T_{\omega_\la}$, where

\begin{equation*}
\begin{split}
{\bf d}^+ &=
\begin{bmatrix}
d^+_{55} & d^+_{45} & d^+_{35} & d^+_{25} & d^+_{15} \\
  & d^+_{44} & d^+_{34} & d^+_{24} & d^+_{14} \\
 &   & d^+_{33} & d^+_{32} & d^+_{13} \\
 &   &   & d^+_{22} & d^+_{12} \\
  &   &   &   & d^+_{11} \\
\end{bmatrix}=
\begin{bmatrix}
2 & 1 & 1 & 0 & 0 \\
  & 0 & 0 & 3 & 0 \\
  &   & 0 & 0 & 1 \\
  &   &   & 2 & 0 \\
  &   &   &   & 0 \\
\end{bmatrix}, \\
{\bf d}^-&=
\begin{bmatrix}
d^-_{12} & d^-_{13} & d^-_{14} & d^-_{15}   \\
         & d^-_{23} & d^-_{24} & d^-_{25}   \\
         &          & d^-_{34} & d^-_{35}   \\
         &          &          & d^-_{45}   \\
\end{bmatrix}=
\begin{bmatrix}
0 & 1 & 1 & 1   \\
  & 0 & 1 & 1   \\
  &   & 0 & 2  \\
  &   &   & 0   \\
\end{bmatrix}.
\end{split}
\end{equation*}
\noindent
\vskip 5mm

(2) Let $\la=(\frac{7}{2},\frac{7}{2},\frac{7}{2},\frac{5}{2},\frac{1}{2})\in \cP_5^{\rm sp}$ and  
$$
T=\quad
\resizebox{.11\hsize}{!}
{\def\lr#1{\multicolumn{1}{|@{\hspace{.75ex}}c@{\hspace{.75ex}}|}{\raisebox{-.04ex}{$#1$}}}\raisebox{-.6ex}
{$\begin{array}{cccc}
\cline{4-4}
& & & \lr{\!3\!} \\
\cline{2-4}
& \lr{2} & \lr{0} &\lr{\!\ov 5\!}\\
\cline{1-4}
\lr{1} & \lr{3} & \lr{\ov 5} &\lr{\!\ov 4\!}\\
\cline{1-1}\cline{2-2}\cline{3-3}\cline{4-4}
\lr{4}&\lr{0} & \lr{\ov 3}&\lr{\!\ov 2\!}\\
\cline{1-1}\cline{2-2}\cline{3-3}\cline{4-4}
\lr{\ov 4}&\lr{\ov 1} & \lr{\ov 1}&\lr{\!\ov 1\!}\\
\cline{1-1}\cline{2-2}\cline{3-3}\cline{4-4} 
\end{array}$}}\quad  \in {\bf KN}^{\mf b_5}_{\la}.
$$ 
\noindent Then
$\Psi_{\la}(T)={\bf T}=
(\Psi_3(T_3),\Psi_4(T_2),\Psi_4(T_1),\Psi_{\rm sp}(T_0))
\in {\bf T}^{\mf b_5}_{\la}$, where
$$\ 
\resizebox{.42\hsize}{!}
{\def\lr#1{\multicolumn{1}{|@{\hspace{.75ex}}c@{\hspace{.75ex}}|}{\raisebox{-.04ex}{$#1$}}}\raisebox{-.6ex}
{$\begin{array}{ccccccccccccccccc}
\cline{14-14}
& & & & & & & & & & & & &\lr{\!3\!}\\
\cline{6-6}\cline{10-10}\cline{14-14}
&  & & & & \lr{2} & & & & \lr{0} & & & &\lr{\!\ov 5\!}\\
\cline{2-2}\cline{6-6}\cline{10-10}\cline{14-14}
& \lr{1} & & & & \lr{3} & & & & \lr{\ov 5} & & & &\lr{\!\ov 4\!}\\
\cline{2-2}\cline{6-6}\cline{10-10}\cline{14-14}
& \lr{4} & & & & \lr{0} & & & & \lr{\ov 3} & & & & \lr{\!\ov 2\!}\\
\cline{2-2}\cline{6-6}\cline{10-10}\cline{14-14}
\ \, &\lr{\ov 4} & &  \, \ & & \lr{\ov 1} & &   \, \ & & \lr{\ov 1} & &  \, \ & & \lr{\!\ov 1\!}& \\
\cline{2-2}\cline{6-6}\cline{10-10}\cline{14-14}\cdashline{1-15}[0.5pt/1pt]\\ 
\ \  & \!\!\!\!\!\!\!\!\!T_3 \!\!\!\!\!\!\!\!\! & & &
& \!\!\!\!\!\!\!\!\!T_2 \!\!\!\!\!\!\!\!\! & & &
& \!\!\!\!\!\!\!\!\!T_1 \!\!\!\!\!\!\!\!\! & & & 
& \!\!\!\!\!\!\!\!\!T_0 \!\!\!\!\!\!\!\!\! & & & 
\end{array}$}}
$$
$$
\resizebox{.4\hsize}{!}
{\def\lr#1{\multicolumn{1}{|@{\hspace{.75ex}}c@{\hspace{.75ex}}|}{\raisebox{-.04ex}{$#1$}}}\raisebox{-.6ex}
{$\begin{array}{cccccccccccccccc}
\\
\cline{11-11}\cline{14-14}
& & & & & & & & &  & \lr{\ov{5}} & & & \lr{\ov{5}} &   \\
\cline{10-11}\cline{14-14}
& & & & & & & & & \lr{\ov{5}} & \lr{\ov{4}} & & & \lr{\ov{4}} &   \\
\cline{7-7}\cline{10-10}\cline{11-11}\cline{14-14}
& & & & & &  \lr{\ov{5}} & &  & \lr{\ov{3}}&\lr{\ov{3}} & & & \lr{\ov{2}} & \\
\cline{2-3}\cline{6-7}\cline{10-11}\cline{14-14} 
&\lr{\ov{5}} & \lr{\ov{3}}& & & \lr{\ov{4}} & \lr{\ov{1}} & & & \lr{\ov{2}} & \lr{\ov{1}}& & & \lr{\ov{1}} &   \\
\cline{2-3}\cdashline{1-5}[0.5pt/1pt]\cline{6-7}\cdashline{8-15}[0.5pt/1pt]
\cline{10-11}\cdashline{12-14}[0.5pt/1pt]\cline{14-14}
&\lr{\ov{4}} & & & & \lr{\ov{1}} & & &  & \lr{\ov{1}} & & & & & \\
\cline{2-2}\cline{6-6}\cline{10-10}
&\lr{\ov{2}} & & & & & \\
\cline{2-2}
\\
\Psi_3(T_3)\!\!\!\!\!\!\!\!\!\!\!\!\!\!\!\!\!\! & & & & 
\Psi_4(T_2)\!\!\!\!\!\!\!\!\!\!\!\!\!\!\!\!\!\! & & & & 
\Psi_4(T_1)\!\!\!\!\!\!\!\!\!\!\!\!\!\!\!\!\!\! & & & & 
\Psi_{\rm sp}(T_0)\!\!\!\!\!\!\!\!\!\!\!\!\!\!\!\!\!\! 
\end{array}$}} \quad .
$$ 
Since
$$
{\bf T}=\quad 
\resizebox{.4\hsize}{!}
{\def\lr#1{\multicolumn{1}{|@{\hspace{.75ex}}c@{\hspace{.75ex}}|}{\raisebox{-.04ex}{$#1$}}}\raisebox{-.6ex}
{$\begin{array}{cccccccccccccccccccc}
\cline{12-12}\cline{14-14} 
& & & & & & & & & & & \lr{\ov{5}} & & \lr{\ov{5}} &   \\
\cline{10-10}\cline{12-12}\cline{14-14} 
& & & & & & & & & \lr{\ov 5} & & \lr{\ov 4} & & \lr{\ov{4}} &   \\
\cline{8-8}\cline{10-10}\cline{12-12}\cline{14-14} 
& & & & & & & \lr{\ov{5}} &  & \lr{\ov{3}}& & \lr{\ov{3}} & & \lr{\ov{2}} &   \\
\cline{2-2}\cline{4-4}\cline{6-6}\cline{8-8}\cline{10-10}\cline{12-12}\cline{14-14} 
&\lr{\ov{5}} & & \lr{\ov{3}}  & & \lr{\ov{4}} & & \lr{\ov{1}} & & \lr{\ov{2}}& & \lr{\ov{1}}& & \lr{\ov{1}} &   \\
\cline{2-2}\cline{4-4}\cline{6-6}\cline{8-8}
\cline{10-10}\cline{12-12}\cdashline{1-15}[0.5pt/1pt]\cline{14-14} 
&\lr{\ov{4}} & & & & \lr{\ov{1}} & & &  & \lr{\ov{1}} & & & &   \\
\cline{2-2}\cline{6-6}\cline{10-10} 
&\lr{\ov{2}} & & & &   \\
\cline{2-2}
&\!\!\text{\em 7\!} & &\text{\em 6\!} & &\text{\em 5\!} & &\text{\em 4\!} & &\text{\em 3\!} & &\text{\em 2\!} & &\text{\em 1\!} 
\end{array}$}} \quad  \in {\bf E}^7,
$$ 
we have 
$$
\ov{\bf T}=(\mc{F}_4\mc{F}_3)\mc{F}_5 {\bf T}=\quad 
\resizebox{.24\hsize}{!}
{\def\lr#1{\multicolumn{1}{|@{\hspace{.75ex}}c@{\hspace{.75ex}}|}{\raisebox{-.04ex}{$#1$}}}\raisebox{-.6ex}
{$\begin{array}{ccccccccccccccccccc}
\cline{7-8}
& & & & &   & \lr{\ov{5}} & \lr{\ov{5}} &  \\
\cline{6-8}
& & & & &  \lr{\ov 5} & \lr{\ov{4}} & \lr{\ov{4}} &  \\
\cline{5-8}
& & & & \lr{\ov{5}} & \lr{\ov{3}} & \lr{\ov{3}}  & \lr{\ov{2}}  &  \\
\cline{2-8}
& \lr{\ov{5}} & \lr{\ov{4}}  & \lr{\ov{2}} & \lr{\ov{1}} & \lr{\ov{1}} & \lr{\ov{1}} & \lr{\ov{1}} &  \\
\cline{2-8}\cdashline{1-11}[0.5pt/1pt]
& \lr{\ov{4}} & \lr{\ov{3}} & \lr{\ov 1} &    \\
\cline{2-4}
& \lr{\ov{2}} &  \\
\cline{2-2}
\end{array}$}}\quad ,
$$\vskip 2mm
\noindent and hence  
$\Theta_\la({\bf T})={\bf V}\otimes t_{4\omega_5}=({\bf T}^{\tt tail}, {\bf T}^{\tt body})\otimes t_{4\omega_5}$, where \vskip 2mm
$$
{\bf T}^{\tt tail}=\quad 
\resizebox{.13\hsize}{!}
{\def\lr#1{\multicolumn{1}{|@{\hspace{.75ex}}c@{\hspace{.75ex}}|}{\raisebox{-.04ex}{$#1$}}}\raisebox{-.6ex}
{$\begin{array}{ccccccccccccccccccc}
\cline{2-4} 
& \lr{\ov{4}} & \lr{\ov{3}} & \lr{\ov 1} &    \\
\cline{2-4}
& \lr{\ov{2}} &  \\
\cline{2-2}
\end{array}$}}\quad\quad
{\bf T}^{\tt body}=\quad 
\resizebox{.24\hsize}{!}
{\def\lr#1{\multicolumn{1}{|@{\hspace{.75ex}}c@{\hspace{.75ex}}|}{\raisebox{-.04ex}{$#1$}}}\raisebox{-.6ex}
{$\begin{array}{ccccccccccccccccccc}
\cline{7-8}
& & & & &   & \lr{\ov{5}} & \lr{\ov{5}} &  \\
\cline{6-8}
& & & & &  \lr{\ov 5} & \lr{\ov{4}} & \lr{\ov{4}} &  \\
\cline{5-8}
& & & & \lr{\ov{5}} & \lr{\ov{3}} & \lr{\ov{3}}  & \lr{\ov{2}}  &  \\
\cline{2-8}
& \lr{\ov{5}} & \lr{\ov{4}}  & \lr{\ov{2}} & \lr{\ov{1}} & \lr{\ov{1}} & \lr{\ov{1}} & \lr{\ov{1}} &  \\
\cline{2-8}
\end{array}$}}\quad .
$$ \vskip 3mm
\noindent Therefore, we have $\Phi_\la({\bf V}\otimes t_{4\omega_5})={\bf d}\otimes t_{\omega_\la}\in \B_{\bi}\otimes T_{\omega_\la}$, where
\begin{equation*}
\begin{split}
{\bf d}^+ &=
\begin{bmatrix}
d^+_{55} & d^+_{45} & d^+_{35} & d^+_{25} & d^+_{15} \\
  & d^+_{44} & d^+_{34} & d^+_{24} & d^+_{14} \\
 &   & d^+_{33} & d^+_{32} & d^+_{13} \\
 &   &   & d^+_{22} & d^+_{12} \\
  &   &   &   & d^+_{11} \\
\end{bmatrix}=
\begin{bmatrix}
4 & 0 & 0 & 0 & 1 \\
  & 2 & 0 & 0 & 1 \\
  &   & 1 & 1 & 0 \\
  &   &   & 0 & 1 \\
  &   &   &   & 0 \\
\end{bmatrix},\\
{\bf d}^-&=
\begin{bmatrix}
d^-_{12} & d^-_{13} & d^-_{14} & d^-_{15}   \\
         & d^-_{23} & d^-_{24} & d^-_{25}   \\
         &          & d^-_{34} & d^-_{35}   \\
         &          &          & d^-_{45}   \\
\end{bmatrix}=
\begin{bmatrix}
0 & 0 & 0 & 1   \\
  & 0 & 1 & 0   \\
  &   & 0 & 1  \\
  &   &   & 1   \\
\end{bmatrix}. \\
\end{split}
\end{equation*}


\vskip 5mm

{\small
\begin{thebibliography}{HK}

 
\bibitem{BZ-1}
A. Berenstein, A. Zelevinsky, {\em Tensor product multiplicities and convex polytopes in
partition space}, J. of Geom. and Physics \textbf{5} (1988)  453--472.  
  
  
\bibitem{C}
G. Cliff, {\em Crystal bases and Young tableaux}, J. Algebra \textbf{202} (1998) 10--35.  
  
\bibitem{HK}
J. Hong, S.-J. Kang, {\em Introduction to Quantum Groups and Crystal Bases}, Graduate Studies in
Mathematics 42,  Amer. Math. Soc., 2002.

\bibitem{FOS}
G. Fourier, M. Okado, A. Schilling, {\em Kirillov-Reshetikhin crystals for nonexceptional types}, Adv. Math. \textbf{222} (2009) 1080--1116.

\bibitem{FRS}
J. Fuchs, U. Ray, and C. Schweigert, {\em Some automorphisms of generalized Kac Moody
algebras}, J. Algebra \textbf{191} (1997) 518--540.

\bibitem{Ful}
W. Fulton, {\em Young tableaux, with Application to Representation
theory and Geometry}, Cambridge Univ. Press, 1997.

\bibitem{Ha}
T. Hayashi, {\em $q$-analogues of Clifford and Weyl algebras- spinor and oscillator representations of quantum enveloping algebras}, Comm. Math. Phys. \text{127} (1990) 129--144.

\bibitem{HL}
J. Hong, H.-M. Lee, {\em Young tableaux and crystal $B(\infty)$ for finite simple Lie algebras}, J. Algebra \textbf{320} (2008) 3680--3693.

\bibitem{Kam}
J. Kamnitzer, {\em  The crystal structure on the set of Mirkovi\'{c}-Vilonen polytopes},  Adv. Math. \textbf{215} (2007) 66--93.

\bibitem{Kas91}
M. Kashiwara, {\em On crystal bases of the $q$-analogue of universal enveloping algebras}, Duke Math. J. \textbf{63} (1991) 465--516.

\bibitem{Kas93}
M. Kashiwara, {\em The crystal base and Littelmann's refined Demazure character formula}, Duke Math. J. \textbf{71} (1993) 839--858.

\bibitem{Kas95}
M. Kashiwara, {\em On crystal bases}, Representations of groups,
CMS Conf. Proc., 16, Amer. Math. Soc., Providence, RI, (1995) 155--197.

\bibitem{Kas96}
M. Kashiwara, {\em Similarity of crystal bases}, Contemp. Math., \textbf{194} (1996) 177--186.

\bibitem{KashNaka}
M. Kashiwara, T. Nakashima, {\em Crystal graphs for representations of the $q$-analogue of classical Lie algebras}, J. Algebra \textbf{165} (1994) 295--345.

\bibitem{K07}
J.-H. Kwon, {\em Crystal graphs for Lie superalgebras and Cauchy decomposition}, J. Algebraic Combin. \textbf{25} (2007) 57--100.

\bibitem{K09}
J.-H. Kwon,
{\em Demazure crystals of generalized Verma modules and a flagged RSK correspondence}, 
J. Algebra \textbf{322} (2009) 2150--2179.

\bibitem{K12}
J.-H. Kwon, {\em Littlewood identity and crystal bases}, Adv. Math. \textbf{230} (2012) 699--745.

\bibitem{K15-1}
J.-H. Kwon, {\em Combinatorial extension of stable branching rules for classical groups},  preprint (2015) arXiv:1512.01877, to appear in Trans. Amer. Math. Soc.

\bibitem{K15}
J.-H. Kwon,  {\em Super duality and crystal bases for quantum orthosymplectic
superalgebras}, Int. Math. Res. Not. \textbf{23} (2015) 12620--12677.

\bibitem{K16}
J.-H. Kwon,  {\em Super duality and crystal bases for quantum orthosymplectic
superalgebras II}, J. Algebraic Combin. \textbf{43} (2016) 553--588.

\bibitem{K16-2}
J.-H. Kwon, {\em A crystal embedding into Lusztig data of type A},  J. Combin. Theory Ser. A \textbf{154} (2018) 422--443.

\bibitem{La}
A. Lascoux, {\em Double crystal graphs}, in: Studies in Memory of Issai Schur, in: Progr. Math., vol. 210, Birkh\"{a}user, 2003, pp. 95--114.

\bibitem{Le02} 
C. Lecouvey, {\em Schensted-type correspondence, plactic monoid, and jeu de taquin for type $C_n$}, J. Algebra \textbf{247} (2002) 295--331.

\bibitem{Le03}
C. Lecouvey, {\em Schensted-type correspondences and plactic monoids for types $B_n$ and $D_n$}, J. Algebraic Combin. \textbf{18} (2003) 99--133.

\bibitem{LOS}
C. Lecouvey, M. Okado, M. Shimozono, {\em Affine crystals, one-dimensional sums and parabolic Lusztig q-analogues}, Math. Z. \textbf{271} (2012) 819--865.

\bibitem{Lu90}
G. Lusztig, {\em Canonical bases arising from quantized universal enveloping algebras},  J. Amer. Math. Soc. \textbf{3}  (1990)  447--498.

\bibitem{Lu90-2}
G. Lusztig, {\em Canonical bases arising from quantized enveloping algebras II}, Progr. Theor. Phys. Suppl., \textbf{102} (1990)  175--201.

\bibitem{Lu93}
G. Lusztig, {\em Introduction to quantum groups}, Progress in Math. \textbf{110}, Birkh\"{a}user, 1993.

\bibitem{Na}
T. Nakashima, {\em Crystal base and a generalization of the LR rule for the classical Lie algebras}, Comm. Math. Phys. \textbf{154} (1993) 215--243.

\bibitem{NS03}
S. Naito, D. Sagaki, {\em Path model for a level-zero extremal weight module over a quantum affine algebra}, Int. Math. Res. Not. (2003) 1731--1754.

\bibitem{NS05}
S. Naito, D. Sagaki, {\em Crystal base elements of an extremal weight module fixed by a diagram automorphism}, Algebr. Represent. Theory \textbf{8} (2005) 689--707.

\bibitem{Re}
M. Reineke,
{\em On the coloured graph structure of Lusztig's canonical basis}, Math. Ann. \textbf{307} (1997) 705--723.

\bibitem{S94}
Y. Saito, {\em  PBW basis of quantized universal enveloping algebras}, Publ. Res. Inst. Math. Sci. \textbf{30} (1994) 209--232.

\bibitem{SST16}
B. Salisbury, A. Schultze, P. Tingley, {\em PBW bases and marginally large tableaux in type D}, preprint (2016), arXiv:1606.02517

\bibitem{SST17}
J. Criswell, B. Salisbury, P. Tingley, {\em PBW bases and marginally large tableaux in types B and C}, preprint (2017), arXiv:1708.04311.


\bibitem{SST}
B. Salisbury, A. Schultze, P. Tingley, {\em Combinatorial descriptions of the crystal structure on certain PBW bases}, Transformation Groups (2017). https://doi.org/10.1007/s00031-017-9434-9.

\end{thebibliography}
}
\end{document}